\documentclass[11pt,twoside]{article}

\usepackage{epsfig,amsfonts,color}
\usepackage[final]{pdfpages}
\usepackage{amsmath}
\usepackage{blkarray}
\usepackage{bm}
\usepackage{algorithm}
\usepackage{algpseudocode}
\usepackage{multirow}
\renewcommand{\algorithmicfunction}{\textbf{Function}}
\algdef{SE}[TASK]{Task}{EndTask}%
   [2]{\algorithmicfunction\ \tt{#1}\ifthenelse{\equal{#2}{}}{}{(#2)}}%
   {\algorithmicend\ \algorithmicfunction}%
\usepackage{graphicx}
\usepackage[margin = 1cm]{caption}
\usepackage{subcaption}
\usepackage{placeins}

\usepackage{amssymb, palatino, geometry,url}
\usepackage[colorlinks=true,linkcolor=blue,citecolor=blue,urlcolor=blue]{hyperref}
\usepackage[titletoc,title]{appendix}

\usepackage{multirow}
\usepackage{mathtools}
\usepackage{fancyhdr}
\newcommand{\be}{\begin{equation}}
\newcommand{\ee}{\end{equation}}
\newcommand{\bea}{\begin{eqnarray}}
\newcommand{\eea}{\end{eqnarray}}

\newcommand{\bvec}{\left(\begin{array}{c}}
\newcommand{\evec}{\end{array}\right)}
\newcommand{\bsub}{\begin{subequations}}
\newcommand{\esub}{\end{subequations}}

\newcommand{\mn}{\mathcal{N}}
\newcommand{\me}{\mathcal{E}}
\newcommand{\mg}{\mathcal{G}}

\newcommand{\mc}{\mathcal{C}}

\newcommand{\msg}{\mathcal{SG}}
\newcommand{\ms}{\mathcal{S}}
\newcommand{\bx}{\boldsymbol{x}}

\usepackage{accsupp}

\usepackage{listings}
\DeclareCaptionFont{white}{\color{white}}
\DeclareCaptionFormat{listing}{{\parbox{\dimexpr\textwidth-2\fboxsep\relax}{#1#2#3}}}
\usepackage[dvipsnames]{xcolor}

\lstdefinelanguage{julia}
{
keywordsprefix=\@,
morekeywords={
exit,whos,edit,load,is,isa,isequal,typeof,tuple,ntuple,uid,hash,finalizer,convert,promote,
subtype,typemin,typemax,realmin,realmax,sizeof,eps,promote_type,method_exists,applicable,
invoke,dlopen,dlsym,system,error,throw,assert,new,Inf,Nan,pi,im,begin,while,for,in,return,
break,continue,macro,quote,let,if,elseif,else,try,catch,end,bitstype,ccall,do,using,module,
import,export,importall,baremodule,immutable,local,global,const,Bool,Int,Int8,Int16,Int32,
Int64,Uint,Uint8,Uint16,Uint32,Uint64,Float32,Float64,Complex64,Complex128,Any,Nothing,None,
function,type,typealias,abstract,get_node,add_edge,create_estimation_model,set_solution,
solve, get_solution, solve_ss_problem, create_estimation_problem, addnode, Partition, 
assemble_optigraph, local_subgraphs, apply_partition, OptiGraph, optimizer_with_attributes,
BendersAlgorithm, optimize, aggregate, source_graph, aggregate_to_depth, add_subgraph, fill, set_to_node_objectives, set_optimizer
},
morekeywords = [2]{triggered_by,compute_time,trigger_during_busy,send_on,delay,send_wait,start},
sensitive=true,
morecomment=[l]{\#},
morestring=[b]',
morestring=[b]"
}

\usepackage{authblk}
\usepackage{amsthm}
\theoremstyle{plain}

\theoremstyle{definition}

\begin{document}

\title{Graph-Based Modeling and Decomposition of\\ Hierarchical Optimization Problems}

\author{David L. Cole$^{1}$, Filippo Pecci$^{2,3}$, Omar J. Guerra$^{4}$\\ Harsha Gangammanavar$^{5}$, Jesse D. Jenkins$^{6,7}$, and Victor M. Zavala$^{1, 8 *}$\\
 {\small $^{1}$ Department of Chemical and Biological Engineering}\\
 {\small University of Wisconsin-Madison, Madison, WI 53706, USA}\\
 {\small $^{2}$ CMCC Foundation - Euro-Mediterranean Center on Climate Change}\\
 {\small Lecce, Italy 73100}\\
 {\small $^{3}$ RFF-CMCC European Institute on Economic and the Environment}\\
 {\small Lecce, Italy 73100}\\
 {\small $^{4}$ Grid Planning and Analysis Center}\\
 {\small National Renewable Energy Laboratory, Golden, CO 80401, USA}\\
 {\small $^{5}$ Operations Research and Engineering Management}\\
 {\small Southern Methodist University, Dallas, TX 75205, USA}\\
 {\small $^{6}$ Andlinger Center for Energy and Environment}\\
 {\small Princeton University, Princeton, NJ 08540, USA}\\
 {\small $^{7}$ Department of Mechanical and Aerospace Engineering}\\
 {\small Princeton University, Princeton, NJ 08540, USA}\\
 {\small $^{8}$ Mathematics and Computer Science Division}\\
 {\small Argonne National Laboratory, Lemont, IL 60439, USA}\\
 {\small $^*$ Corresponding author email address: victor.zavala@wisc.edu}
}

\date{}

\maketitle

\abstract{We present a graph-theoretic modeling approach for hierarchical optimization that leverages the OptiGraph abstraction implemented in the {\tt Julia} package {\tt Plasmo.jl}. We show that the abstraction is flexible and can effectively capture complex hierarchical connectivity that arises from decision-making over multiple spatial and temporal scales (e.g., integration of planning, scheduling, and operations in manufacturing and infrastructures). We also show that the graph abstraction facilitates the conceptualization and implementation of decomposition and approximation schemes. Specifically, we propose a graph-based Benders decomposition (gBD) framework that enables the exploitation of hierarchical (nested) structures and that uses graph aggregation/partitioning procedures to discover such structures. In addition, we provide a {\tt Julia} implementation of gBD, which we call {\tt PlasmoBenders.jl}. We illustrate the capabilities using examples arising in the context of energy and power systems.}
\\

\noindent {\bf Keywords}: Optimization, Graph Theory, Hierarchical, Multi-Scale, Decomposition, Julia.

\section{Introduction}

Hierarchical optimization problems often involve connectivity across different decision levels and scales (e.g., physical assets of a system are coupled across different temporal or spatial scales) \cite{anandalingam1992hierarchical,biegler2014multi}. These types of problems are also often high-dimensional and their solution requires approximation schemes (e.g., aggregation/reduction schemes) or decomposition schemes (e.g., Benders or Lagrangian decomposition) \cite{chaieb2015hierarchical,conejo2006decomposition,mitrai2021efficient}. Hierarchical optimization problems often arise from the integration of decisions in supply chains, planning, scheduling, and operation/control levels in industrial systems such as energy systems \cite{braun2018hierarchical,guo2017hierarchical,jacobson2024computationally,kong2017hierarchical,tian2015hierarchical}, chemical processes \cite{abaecherli2017hierarchical,van2024hierarchical,wu2007hierarchical}, agricultural systems \cite{Albornoz2015,liu2021hierarchical}, water infrastructure networks \cite{ocampo2010model}, and healthcare systems \cite{jemai2013home}. 
\\

To handle the complexity of hierarchical optimization problems, it is necessary to navigate and exploit their unique {\it connectivity structure}.  For example, Benders Decomposition (BD) \cite{jacobson2024computationally,pecci2024regularized,rahmaniani2017benders,zhang2024decomposition}, Lagrangian Decomposition (LD)   \cite{li2010production,zhang2024decomposition}, and the Alternating Direct Method of Multipliers (ADMM) \cite{braun2018hierarchical,khaki2018hierarchical,shin2019hierarchical} can be used to decompose and exploit the structure of these problems. Other approaches involve computing approximate solutions by solving hierarchical layers (or components of each layer) in a sequential manner \cite{guerra2019integrated,long2016hierarchical}; this approximation strategy is used, for instance, for handling integrated unit commitment and economic dispatch problems in power systems \cite{abujarad2017recent,atakan2022towards,lara2024powersimulations} and for solving dynamic programs over long horizons \cite{shin2021diffusing}. However, to our knowledge, there is currently no modeling framework for expressing hierarchical connectivity in a unified manner. As a result, approximation/decomposition approaches are often problem-specific, and successful applications typically rely on specialized domain knowledge.
\\

In this work, we present a general, graph-theoretic modeling approach for optimization that facilitates the navigation, discovery, and manipulation of hierarchical structures. This approach is based on the OptiGraph abstraction implemented in the {\tt Julia} package {\tt Plasmo.jl} \cite{jalving2022graph}. Under this abstraction, the structure of an optimization problem is expressed as a graph in terms of nodes and (hyper)edges. Nodes embed optimization problems with their own variables, constraints, objective functions, and data; hyperedges contain linking constraints that couple variables of optimization problems from different nodes. Importantly, the node of an OptiGraph can embed another OptiGraph, thus enabling the modeling of hierarchical structures. 
\\

The OptiGraph representation is abstract and can be applied to problems arising in distinct {\it application domains}. The OptiGraph abstraction is also agnostic to the meaning of {\it its modeling elements}; for example, a node can capture either a spatial domain (e.g., a geographical region), a temporal domain (e.g., a day or month), a probabilistic scenario (e.g., for stochastic programs), or a physical asset (e.g., a battery of a power plant). Furthermore, graph structures can be manipulated via aggregation and partitioning, which can facilitate the visualization/discovery of structures and the re-casting of problem structure in a way that is amenable to tailored solution algorithms. For instance, an algorithm that requires an acyclic graph structure (e.g., nested BD requires a tree structure) could potentially be applied to a general graph by eliminating cycles via node aggregation.
\\

We illustrate how the graph abstraction can be used to enable the implementation of approximation and decomposition schemes. Specifically, we introduce a graph-based Benders decomposition scheme (gBD) that generalizes a nested Benders decomposition algorithm \cite{birge1985decomposition,pereira1985stochastic,pereira1991multi} to general graph structures. In doing so, we are extending the work of Brunaud \cite{brunaud2019models} who originally proposed using {\tt Plasmo.jl}'s OptiGraph abstraction as a basis for BD. We also show how graph operations facilitate the application of BD to a variety of problem classes and structures. For example, nested BD (also known as dual dynamic programming) is traditionally applied to decompose problems over a time domain, but the graph structure enables its use in spatial and spatio-temporal domains (e.g., infrastructure networks). We also introduce a {\tt Julia} implementation of the gBD scheme that we call {\tt PlasmoBenders.jl}.
\\

Several applications of graph-based modeling frameworks have been reported in the literature. For example, the work of Daoutidis and co-workers represents optimization problems as graphs \cite{allman2019decode,daoutidis2019decomposition,mitrai2020decomposition,mitrai2022stochastic,mitrai2023graph,tang2018optimal,tang2023resolving} but consider a fine level of granularity where individual variables and constraints are represented as nodes/edges. Mitrai and Daoutidis have used this representation together with a network analysis approach for detecting hierarchical structures \cite{mitrai2021efficient}. The OptiGraph abstraction presented here works at a different level of granularity, in which nodes and edges represent entire optimization subproblems. This level of representation is also similar to that of Marcucci and co-workers \cite{marcucci2024graphs,marcucci2024shortest,morozov2025mixed} in which they represent optimization problems in the nodes of a graph with edges representing linking constraints. They solve a mixed integer program to determine a subgraph of the original graph (e.g., a path between two nodes) and solve the resulting subgraph for control applications. Their representation of the original graph could be effectively captured by the OptiGraph abstraction.
\\

This work advances and distinguishes itself from the existing literature in several ways. i) While building on the OptiGraph abstraction introduced in \cite{jalving2022graph}, we formalize the mathematical notation for modeling with hierarchical graphs and manipulating their graph structure. ii) We formalize the BD algorithm for the graph-based problem (gBD). While Brunaud \cite{brunaud2019models} originally introduced Benders decomposition on the nodes of {\tt Plasmo.jl}, they did not formalize the algorithm or discuss the implications of generalizing the algorithm to the graph structure. iii) We introduce an open-source software tool for implementing gBD. While Brunaud \cite{brunaud2019models} initially developed a software tool for BD with {\tt Plasmo.jl}, this tool is not maintained. Further, it relies on an implementation of {\tt Plasmo.jl} which is more node-centric (every node is an independent modeling object in {\tt Julia}), whereas {\tt Plasmo.jl} under more recent versions has been made subgraph-centric (every subgraph is an independent modeling object in {\tt Julia}). The latter approach is more memory efficient and follows the mathematical formalization. In addition, {\tt PlasmoBenders.jl} also differs from the {\tt SDDP.jl} \cite{dowson2021sddp} package. This latter package is focused on stochastic problem instances---with special functionality for working with probabilities and scenarios---and is implemented with {\tt JuMP.jl} as its primary backend, whereas {\tt PlasmoBenders.jl} uses the OptiGraph abstraction and is more general in its modeling and partitioning capabilities. iv) We present many important implications of generalizing the hierarchical optimization problems to the graph structure. These include highlighting how approximation schemes can be implemented and tested (see Case Study 1) as well as how the graph structure can be flexibly exploited to improve or enhance algorithm performance. v) Finally, the implementation of gBD as outlined in this paper provides a framework for defining structure-exploiting algorithms to graph-based problems in a flexible and agnostic way. This framework could be extended to other decomposition approaches (e.g., ADMM). 
\\

The paper is structured as follows. Section \ref{sec:abstraction} shows how the OptiGraph abstraction can be used to model hierarchical optimization problems and provides a tutorial overview of {\tt Plasmo.jl} functionalities. Section \ref{sec:Decomposition} shows how approximation and decomposition schemes can be applied to graph structures. Section \ref{sec:case_studies} provides challenging case studies, including a tri-level energy market problem showing how graphs can capture hierarchical coupling and can be used to obtain approximate solutions, as well as a capacity expansion model and a production cost model for power systems that are solved using gBD. Section \ref{sec:future_work} provides conclusions and discusses directions of future work. 

\section{Hierarchical Optimization}\label{sec:abstraction}

Hierarchical optimization problems can be defined in several different ways \cite{anandalingam1992hierarchical,bennett2023hierarchical,bracken1973mathematical,ishikawa2015effective}. For our purposes, a hierarchical optimization problem exhibits a {\em single} (centralized/composite) objective function, and its constraint set has a nested algebraic structure \cite{ishikawa2015effective}. This is not to be confused with multi-level optimization problems, where lower-level optimization problems (with their own objective) are embedded in the constraint set of a higher-level problem \cite{bracken1973mathematical}. Multi-level problems can also be expressed as graphs, but they require specialized solution approaches (e.g., replacing low-level optimization problems with their optimality conditions) that are not considered here and that are left as future work. 

\subsection{OptiGraph Abstraction}

The OptiGraph abstraction was introduced by Jalving and co-workers \cite{jalving2022graph}, and we follow their notation and definitions and expand on them for representing hierarchical problems. An OptiGraph is a hierarchical hypergraph containing sets of OptiNodes $\mathcal{N}$ and OptiEdges $\mathcal{E}$. An OptiNode is a modeling object that embeds an optimization problem with its own objective, variables, constraints, and data; an OptiEdge is a modeling object that embeds constraints that link variables of multiple OptiNodes (thus generating connectivity). The OptiEdges can connect {\it more} than two nodes (they are hyperedges); for simplicity, we refer to hyperedges as edges. Importantly, an {\em OptiNode can also embed an entire OptiGraph}, thus providing the ability to represent hierarchical structures. 
\\

When solving the full problem, the individual objectives of the OptiNodes are combined into a single scalar composite objective function (e.g., via weighted summation). For simplicity, we will refer to OptiNodes, OptiEdges, and OptiGraphs as nodes, edges, and graphs (or subgraphs in the hierarchical case) when it is clear from context that we are referencing the OptiGraph abstraction. An OptiGraph for an optimization problem is denoted as $\mg$ with the set of nodes and edges denoted as $\mn(\mg)$ and $\me(\mg)$, respectively. We also denote edges incident to node $n$ as $\me(n)$ and denote nodes connected by edge $e$ as $\mn(e)$.

An optimization problem modeled under the OptiGraph abstraction can be expressed in the following abstract form:
\begin{subequations}\label{eq:optigraph}
    \begin{align}
        \min_{\{\bx_n\}_{n \in \mn(\mg)}} &\; f\left(\{\bx_n \}_{n \in \mn(\mg)}\right) & (\textrm{Objective}) \label{eq:optigraph_objective} \\
        \textrm{s.t.} &\; \bx_n \in \mathcal{X}_n, \quad n \in \mn(\mg), \quad & (\textrm{Node Constraints}) \label{eq:optigraph_nodes} \\
        &\; g_e(\{\bx_n\}_{n \in \mn(e)})\geq 0, \quad e \in \me(\mg). \quad & (\textrm{Link Constraints}) \label{eq:optigraph_edges}
    \end{align}
\end{subequations}
\noindent Here, $\{\bx_n \}_{n \in \mn(\mg)}$ is the set of decision variables over all the nodes in the graph (i.e., $\bx_n$ is the set of variables of node $n$). The constraints on the node $n$ are represented by the set $\mathcal{X}_n$. The objective function \eqref{eq:optigraph_objective} is a scalar composite function of any variables contained on nodes in the graph. In \eqref{eq:optigraph}, there is no restriction on variables being continuous or discrete or on objectives or constraints being linear or nonlinear; as such, from a modeling standpoint, the OptiGraph abstraction can express structures of linear programs (LPs), mixed-integer linear programs (MILPs), nonlinear programs (NLPs), and mixed-integer nonlinear programs (MINLPs). Furthermore, the linking constraints in \eqref{eq:optigraph_edges} are written as inequalities, but the abstraction is not limited to inequalities, and the software implementation in {\tt Plasmo.jl} supports both equality and inequality constraints on the edges. In this work, we will study how the graph abstraction facilitates the development of algorithms for the solution of LPs and MILPs; however, the abstraction can also potentially be used to implement algorithms for solving nonlinear problems.

\begin{figure*}[!htp]
    \centering
    \includegraphics[width = 0.9\textwidth]{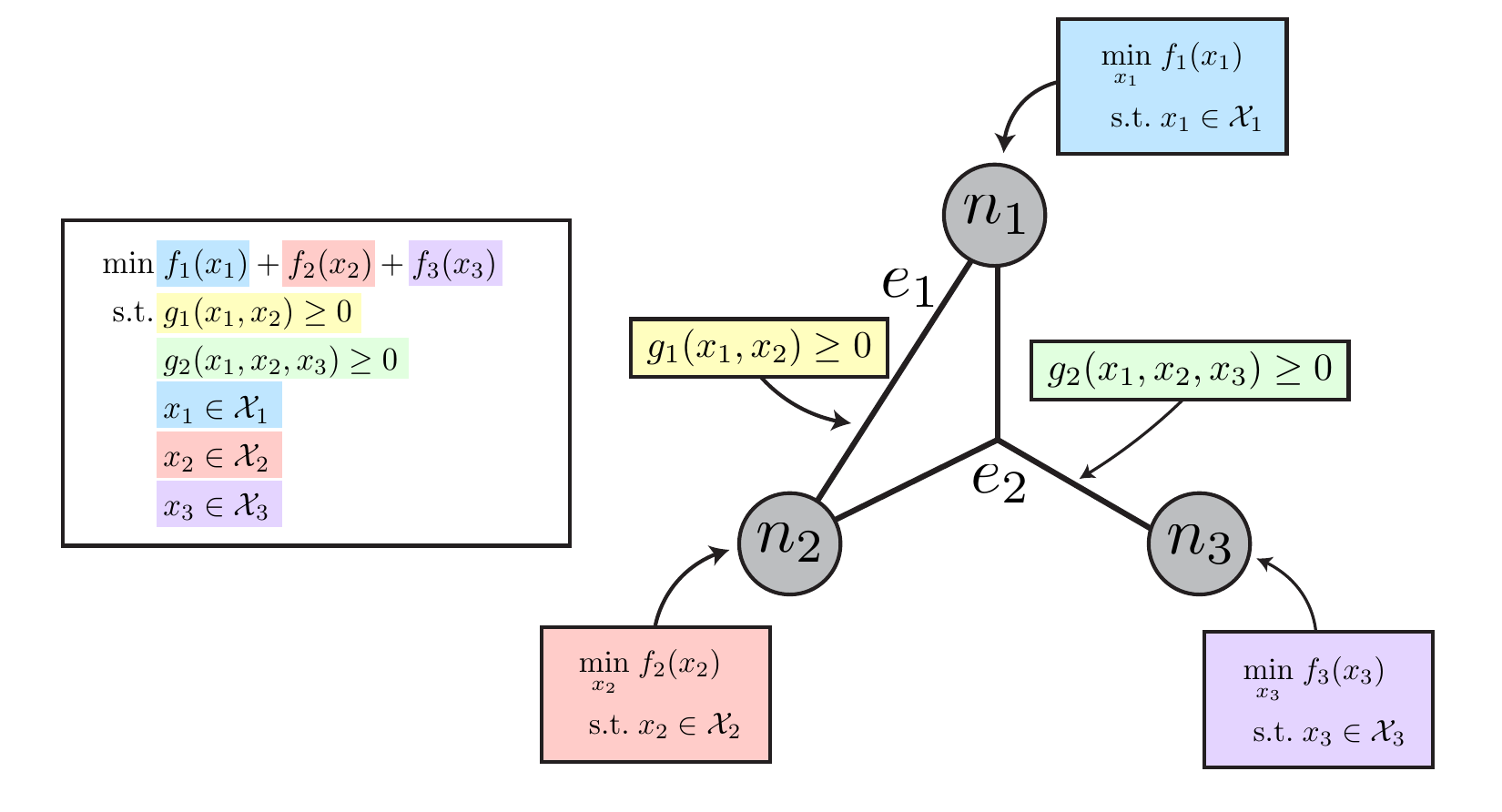}
        \vspace{-0.1in}
    \caption{Visualization of an OptiGraph where an optimization problem (left) is partitioned into the nodes and edges of a graph (right). Nodes contain their own objective, variables, constraints, and data, while edges contain constraints that link variables across nodes.}
    \label{fig:plasmo_overview}
\end{figure*}

Figure \ref{fig:plasmo_overview} provides an example of the OptiGraph structure, where the optimization problem shown on the left is modeled as an OptiGraph, with variables $x_{1}$, $x_{2}$, and $x_{3}$ placed on nodes $n_1$, $n_2$, and $n_3$, respectively. In this case, the objective functions are separable and can be placed in graph nodes. It is not always the case that the objective function is separable among the nodes, and the OptiGraph supports non-separable objective functions. This distinction between a node objective and a graph objective is primarily important for the software implementation (which will be discussed later), as the node objectives are not part of the mathematical definition of \eqref{eq:optigraph}.

The OptiGraph abstraction can be applied to optimization problems in a flexible manner. For instance, the coarsest graph representation of an optimization problem is obtained by placing all variables/constraints on a single node ($|\mathcal{N}(\mathcal{G}) | = 1$ and thus $\me(\mg) = \emptyset$), while the finest graph representation is obtained placing every variable on a separate node. A wide variety of intermediate graph structures can be obtained by partitioning/placing variables and constraints to different nodes. In building the graph, the number and placement of edges are determined automatically based on the placement of decision variables into the nodes of the graph, which determines whether a constraint becomes a nodal constraint, \eqref{eq:optigraph_nodes}, or a linking constraint \eqref{eq:optigraph_edges}. Often, there may be a natural graph structure that a user may want to use (e.g., arising from physical, spatial, or temporal connectivity), but in some cases the graph structure might need to be derived from advanced partitioning or aggregation procedures \cite{karypis1997metis,schlag2023high}.

The OptiGraph abstraction also supports hierarchical structures. For instance, any node in an OptiGraph can embed an OptiGraph, or one can encounter the situation in which an OptiGraph is composed of multiple OptiGraphs (subgraphs) that are linked together. In this way, nodes can belong to subgraphs that in turn belong to a higher-level OptiGraph. Here, a subgraph is defined as a subset of nodes and a subset of edges of the original OptiGraph, such that all nodes are a subset of the nodes of the original OptiGraph, and the subset of edges contains only edges connecting nodes in the subgraph. Each subgraph could likewise be further partitioned into additional subgraphs, and there are no limitations to the number of embedded subgraphs possible within an OptiGraph.  

We now formalize the notation and terminology that we will use for representing hierarchy and layering. The graph containing all nested subgraphs for a problem of interest will be given by $\mg$. We will use the term ``level'' to refer to how ``deep'' a subgraph is in a graph (i.e., how many times subgraphs are nested within other graphs). Subgraphs will be annotated by indices on $\mg$, and the number of indices will indicate the number of levels below the main OptiGraph $\mg$. For example, first-level graphs are denoted as $\mg_i, i \in \{1, ..., N\}$, where $N$ is the number of subgraphs contained in $\mg$. Since each subgraph $\mg_i, i \in \{ 1, ..., N\}$ can likewise contain nested subgraphs, these additional subgraphs are annotated by the addition of a second index on $\mg_i$, such as $\mg_{i, j}, i \in \{1, ..., N\}, j \in \{ 1, ..., N_i\}$, where $N_i$ is the number of subgraphs contained in $\mg_i$. Indices can be added for additional levels with no limit on the number of levels possible.

We will refer to the number of indices as the $``$depth" of subgraphs in the graph (e.g., $\mg_i$ is at a depth of one and $\mg_{i, j}$ is at a depth of two). For any subgraph on any level, we annotate the graph as 
\begin{equation}\label{eq:subgraph_notation}
    \mg_{\circ}(\{\mg_{\circ, i}\}_{i \in \{1, ..., N_{\circ} \}}, \mn_{\mg_{\circ}}, \me_{\mg_{\circ}}),     
\end{equation}
\noindent where $\circ$ represents any set of indices mapping the graph, $N_\circ$ is the number of subgraphs contained on $\mg_{\circ}$, $\{\mg_{\circ, i}\}_{i \in \{1, ..., N_* \}}$ represents the set of subgraphs on $\mg_{\circ}$, and $\mn_{\mg_{\circ}}$ and $\me_{\mg_{\circ}}$ are nodes and edges contained in $\mg_{\circ}$ but independent of its subgraphs $\{\mg_{\circ, i}\}_{i \in \{1, ..., N_* \}}$. We also define notation $\msg(\cdot)$ such that 
\begin{equation}\label{eq:subgraph_func}
 \msg(\mg_{\circ}) := \{\mg_{\circ, i}\}_{i \in \{1, ..., N_* \}}    
\end{equation}
\noindent is the set of subgraphs of a graph $\mg_{\circ}$ on the level directly below $\circ$ (i.e., this does not directly include subgraphs nested on $\mg_{\circ,i}, i \in \{1, ..., N_\circ \}$). Note also that subgraphs may have different numbers of nested subgraphs on them at any level; for instance, it is possible that $|\msg(\mg_{\circ, i})| \ne |\msg(\mg_{\circ,j})|$ for some $i \ne j$. For a graph  $\mg_{\circ}$ containing subgraphs, the notation $\mn(\mg_{\circ})$ and $\me(\mg_{\circ})$ represents the sets of all nodes and edges, respectively, on graph $\mg_{\circ}$ {\it including} the nodes and edges of the subgraphs. In other words, $\mn(\mg_{\circ}) = \mn_{\mg_{\circ}} \cup \{\mn(\mg_{\circ, i}) \}_{i \in \{1, ..., N_\circ \}}$ and $\me(\mg_{\circ}) = \me_{\mg_{\circ}} \cup \{\me(\mg_{\circ, i}) \}_{i \in \{1, ..., N_\circ \}}$.

\begin{figure*}[!htp]
    \centering
    \includegraphics[width = 0.5\textwidth]{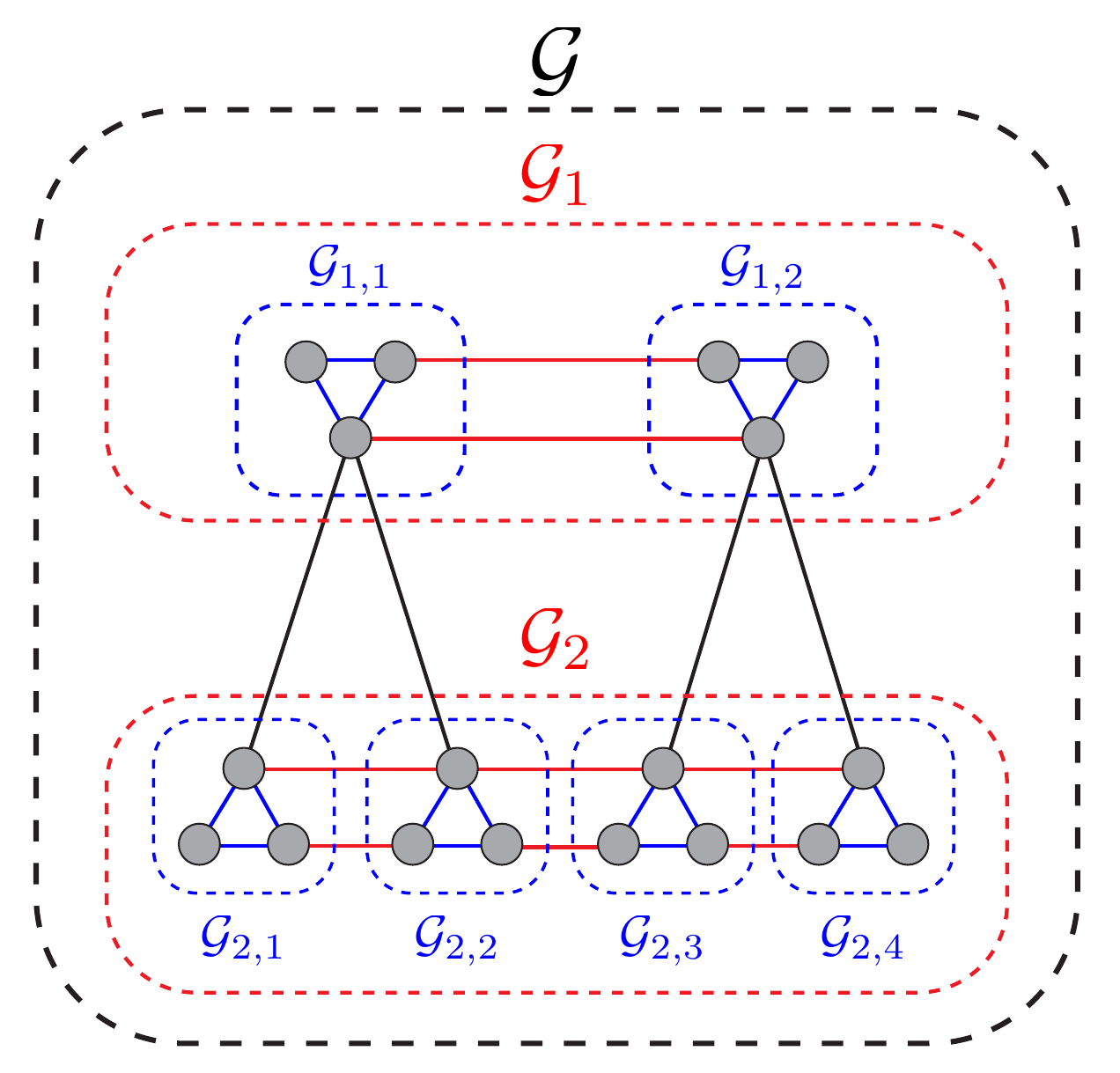}
        \vspace{-0.1in}
    \caption{Example of how OptiGraphs can capture hierarchical connectivity. Here, the red graphs $\mg_1$ and $\mg_2$ are subgraphs of the black $\mg$, highlighted by the dashed borders. Both $\mg_1$ and $\mg_2$ contain a set of lower-level (blue) subgraphs.}
    \label{fig:hier_example}
\end{figure*}

Figure \ref{fig:hier_example} provides an example of an OptiGraph with a hierarchical structure and illustrates the notation used. In this example, the primary graph $\mg(\{\mg_1, \mg_2\}, \emptyset, \mathcal{E}_\mg)$ (in black) contains subgraphs $\mg_1(\{\mg_{1,1}, \mg_{1,2}\}, \emptyset, \mathcal{E}_{\mg_1})$ and $\mg_2(\{ \mg_{2,1}, \mg_{2,2}, \mg_{2,3}, \mg_{2,4} \}, \emptyset, \mathcal{E}_{\mg_2})$ (in red). These subgraphs also contain nested subgraphs (each outlined in blue); $\mg_1$ contains two subgraphs ($\mg_{1,1}$ and $\mg_{1,2}$) while $\mg_2$ contains four subgraphs ($\mg_{2,1}$, $\mg_{2,2}$, $\mg_{2,3}$, and $\mg_{2,4}$). Here, the mid-level (red) subgraphs use only one index for identification, while the lowest-level (blue) subgraphs use two indices, where the first index corresponds to the mid-level subgraph to which they belong. Each subgraph also contains a specific set of edges, which are colored in the figure based on which subgraph the edges belong. For example, the set $\me_\mg$ contains the four black edges drawn from $\mg_{1,1}$ to $\mg_{2,1}$ and to $\mg_{2,2}$ and from $\mg_{1,2}$ to $\mg_{2,3}$ and to $\mg_{2,4}$. Likewise, $\mg_1$ and $\mg_2$ contain sets of edges highlighted in red (e.g., $\me_{\mg_1}$ includes the red edges drawn between $\mg_{1,1}$ and $\mg_{1,2}$). In this case, all nodes are contained within the lowest level (blue) subgraphs, but it would be possible for there to be nodes contained within a higher level subgraph.

Because OptiGraphs can have a nested structure, they are a flexible tool for modeling hierarchical problems. There are at least two ways in which this can be done. The OptiGraph can have an inherent hierarchy based on the nested layers in the graph. For instance, in Figure \ref{fig:hier_example}, one could consider the graphs $\mg_{2,1}$, $\mg_{2,2}$, $\mg_{2,3}$, and $\mg_{2,4}$ as separate systems in a chemical process that are connected together. Each individual graph may contain variables and objectives corresponding to the system it represents. The mid-level graph $\mg_2$ would then capture the higher-level connections between the separate systems, creating a hierarchy between $\mg_2$ and its subgraphs. Another approach to capturing hierarchy in the OptiGraph abstraction is by ordering subgraphs on the same level. In Figure \ref{fig:hier_example}, this can be thought of as $\mg_1$ and $\mg_2$ representing separate hierarchical layers (even though they are on the same hierarchical level of $\mg$). For instance, $\mg_1$ could represent a coarse-resolution planning layer while $\mg_2$ could represent a finer-resolution operations layer; such a structure arises frequently in applications featuring sequential strategic and tactical decisions (e.g., planning or scheduling followed by operations) \cite{jacobson2024computationally,long2016hierarchical,zhang2024decomposition}.

\begin{figure*}[!htp]
    \centering
    \includegraphics[width = 0.5\textwidth]{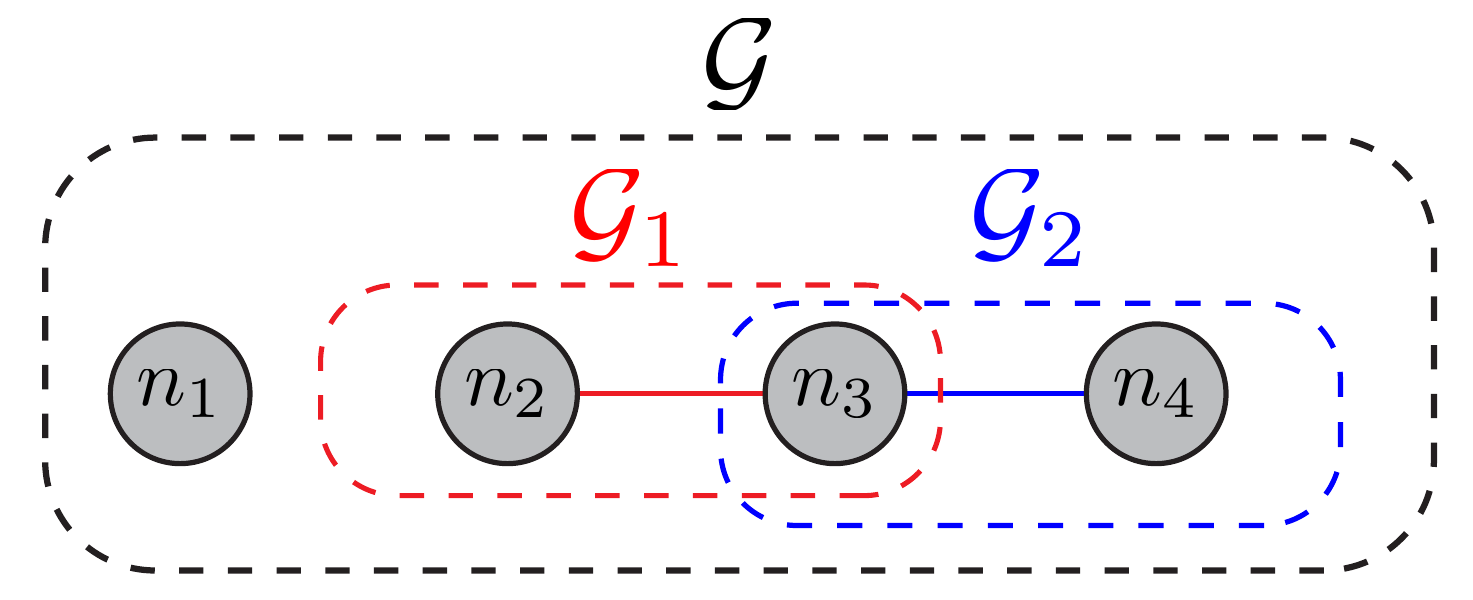}
        \vspace{-0.1in}
    \caption{Example of the hierarchical graph $\mg(\{\mg_1, \mg_2\}, \{n_1 \}, \emptyset)$, showing how a node can belong to multiple subgraphs, creating overlapping regions in the subgraphs.}
    \label{fig:overlapping}
\end{figure*}

There is no restriction that nodes/edges belong to a single graph. For instance, Figure \ref{fig:overlapping} shows the graph $\mg(\{\mg_1, \mg_2\}, \{n_1\}, \emptyset)$ where node $n_3$ belongs to {\it both} $\mg_1$ and $\mg_2$, creating overlapping subgraphs. Many problems will not have overlap, but the ability to express such overlap can have algorithmic implications (e.g., implementation of Schwarz decomposition algorithms \cite{jalving2022graph,shin2020,shin2021}).

\subsection{Software Implementation in {\tt Plasmo.jl}}

The OptiGraph abstraction is implemented in the open-source {\tt Julia} package {\tt Plasmo.jl}, and this provides a user-friendly interface for constructing problems. {\tt Plasmo.jl} leverages the modeling capabilities of {\tt JuMP.jl} \cite{lubin2023} for expressing node and graph subproblems. An OptiGraph is instantiated using the {\tt Plasmo.OptiGraph()} function. Nodes are added to an OptiGraph by calling the {\tt Plasmo.@optinode} macro, which takes as arguments the OptiGraph that the node will belong to and (optionally) a name for the node. Nodes operate similarly to {\tt JuMP.jl} {\tt model} objects and can have variables, constraints, expressions, and objectives added to them using the {\tt JuMP.jl} macros, {\tt @variable}, {\tt @constraint}, and {\tt @objective}, respectively. Edges are created in the graph by calling \texttt{\mbox{\hspace{0pt}} \allowbreak Plasmo.@linkconstraint}. Unlike the {\tt @constraint} macro, which requires a {\tt model} or {\tt OptiNode} object as an argument, the \texttt{ Plasmo.\allowbreak@linkconstraint} function requires an {\tt OptiGraph} object as an argument since the edge will be ``owned'' by the graph rather than by a single node.

To illustrate the basic hierarchical modeling capabilities of {\tt Plamso.jl}, we provide a small modeling example and present the code for assembling this problem. We consider a chemical process that converts a raw material $X$ into a product $Y$. The price for $Y$ fluctuates at different times, and only so much of $Y$ can be produced or sold at a given time. This problem will consider the sizing of a storage/inventory system that can hold $Y$ to be sold at the time periods of higher price. There is a cost associated with building the storage system based on its size, and there is a cost associated with buying raw material $X$. The problem will seek to maximize total profit and is written as: 
\begin{subequations}\label{eq:storage_problem}
    \begin{align}
        \min &\; \alpha\cdot s_{size} + \sum^T_{t = 1} \left( \beta_t\cdot x^{buy}_t - \gamma_t\cdot y^{sell}_t \right) \\
        &\; y^{stored}_{t+1} - y^{stored}_t = y^{save}_t, t = 1, ..., T-1 \label{eq:storage_mass_balance}\\
        &\; y^{save}_t + y^{sell}_t - \zeta x^{buy}_t = 0, t = 1, ..., T \label{eq:y_mass_balance} \\
        &\; 0 \le y^{stored}_t \le s_{size}, t = 1, ..., T \label{eq:storage_upper_bound} \\
        &\; 0 \le y^{sell}_t \le \overline{d}^{sell}, t = 1, ..., T \\
        &\; \underline{d}^{save} \le y^{save}_t \le \overline{d}^{save}, t = 1, ..., T \\
        &\; y^{stored}_1 = \bar{y}^{stored}
    \end{align}
\end{subequations}
\noindent where $s_{size}$ is a decision variable of the maximum storage size, $x^{buy}_t$ is the amount of raw material purchased, $y^{save}_t$ is the amount of product sent to storage (can be negative if it is being removed from storage), $y^{sell}_t$ is the amount of product sold, and $y^{stored}_t$ is the amount of product in storage. The parameters $\alpha$ and $\beta$ are costs associated with their respective variables, $\gamma$ is the price of product $Y$, and $\zeta$ is a conversion factor of $X$ to $Y$. The parameters $\overline{d}^{sell}$, $\underline{d}^{save}$, and $\overline{d}^{save}$ are upper or lower bounds on their respective variables, and $\bar{y}^{stored}$ is the initial amount in storage. Constraint \eqref{eq:storage_mass_balance} is a mass balance on the storage, constraint \eqref{eq:y_mass_balance} is a mass balance on product $Y$, and constraint \eqref{eq:storage_upper_bound} ensures that the storage cannot exceed $s_{size}$. This is a hierarchical problem in the sense that there is a planning/design decision in choosing storage system size and several operations decisions in choosing when to buy $X$ or store/sell $Y$.

\begin{figure*}[!htp]
    \centering
    \includegraphics[width = 0.6\textwidth]{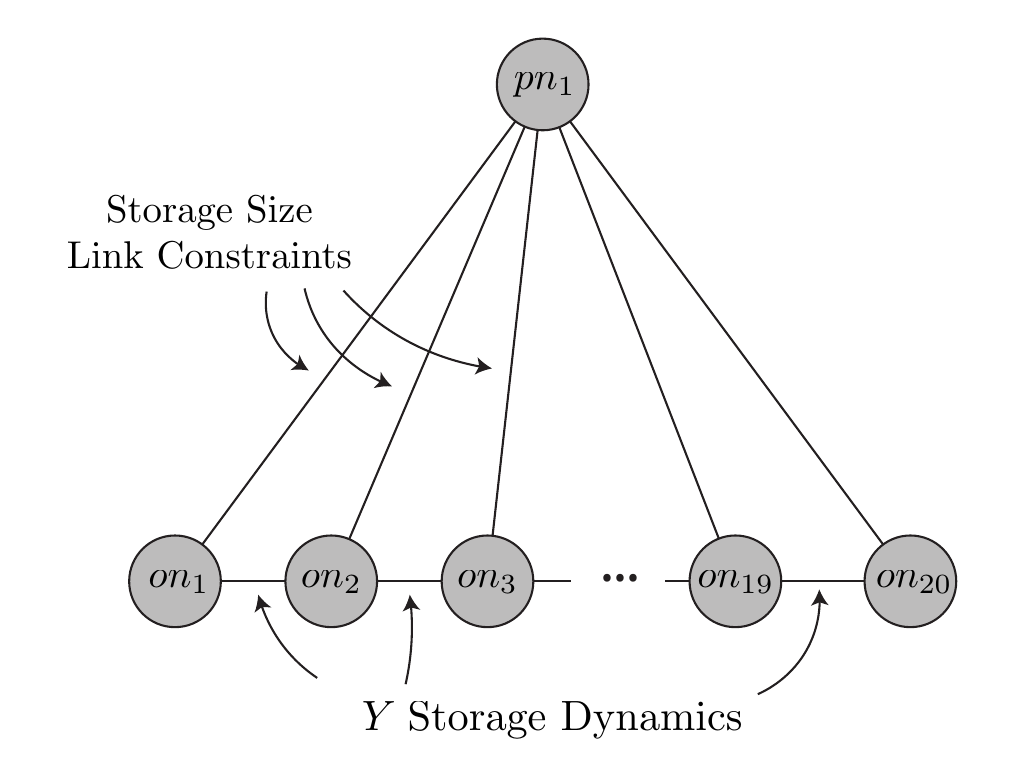}
    \caption{Hierarchical graph structure that can be used to represent the storage problem \eqref{eq:storage_problem} implemented in Code Snippet \ref{code:plasmo_example}. Planning variables are placed on the top/root node $pn_1$, while operations variables at each time point in $T=20$ are embedded within the children nodes $on_t$. The variables of the operation nodes are linked together because storage dynamics are coupled in time.}
    \label{fig:storage_example}
\end{figure*}

The problem in \eqref{eq:storage_problem} can be modeled in {\tt Plasmo.jl} using Code Snippet \ref{code:plasmo_example} for $T=20$ time points. To build this problem, we will place the planning/design variable ($s_{size}$) on its own node, and then link this node to a set of $T$ operations nodes containing the operations variables $x^{buy}_t$, $y^{save}_t$, $y^{sell}_t$, and $y^{stored}_t$. The constraints \eqref{eq:storage_upper_bound} and \eqref{eq:storage_mass_balance} will be captured by linking constraints. Problem data is set on Lines \ref{line:storage_data1} - \ref{line:storage_data2}. The OptiGraph is then instantiated (Line \ref{line:storage_optigraph}), a node is created for the planning variables (Line \ref{line:storage_pn}), and $T$ nodes are created for the operations variables (Line \ref{line:storage_on}). We then set a variable and objective on the planning node (Lines \ref{line:storage_pn_variable} and \ref{line:storage_pn_objective}).  We next loop through the operations nodes and add variables, a constraint, and an objective on each node (Lines \ref{line:storage_on_loop1} - \ref{line:storage_on_loop2}). The initial storage value can then be set (Line \ref{line:storage_initial}) and constraint \eqref{eq:storage_mass_balance} is enforced by adding linking constraints between each operations node (Lines \ref{line:storage_mass_balance1} - \ref{line:storage_mass_balance2}). The planning level decision is enforced by linking constraints between the planning node and each operation node (Line \ref{line:storage_links}). In this code snippet, we have set objectives directly on the nodes since our objective is separable between nodes; to set the graph objective, we will call the function {\tt set\_to\_node\_objectives} (Line \ref{line:storage_graph_obj}), which sums each of the node objectives together to get a scalar composite objective. Finally, we set the optimizer and call {\tt optimize!} on the graph (Lines \ref{line:storage_set_optimizer} and \ref{line:storage_optimize}). The final graph created by this modeling process is shown in Figure \ref{fig:storage_example}.

\begin{figure}[!htp]
    \begin{minipage}[t]{0.9\linewidth}
        \begin{scriptsize}
        \lstset{language=Julia, breaklines = true}
        \begin{lstlisting}[label = code:plasmo_example, caption = Code for generating the graph of storage problem \eqref{eq:storage_problem} in {\tt Plasmo.jl}.] 
# Load in packages
using Plasmo, HiGHS

# Set Problem data
T = 20 |\label{line:storage_data1}|
gamma = fill(5, T); beta = fill(20, T); alpha = 10; zeta = 2
gamma[8:10] .= 20; gamma[16:20] .= 50
d_sell = 50; d_save = 20; d_buy = 15; y_bar = 10 |\label{line:storage_data2}|
            
# Instantiate graph and nodes
graph = OptiGraph() |\label{line:storage_optigraph}|
@optinode(graph, planning_node) |\label{line:storage_pn}|
@optinode(graph, operation_nodes[1:T])|\label{line:storage_on}|
            
# Define planning node data which includes storage size
@variable(planning_node, storage_size >= 0) |\label{line:storage_pn_variable}|
@objective(planning_node, Min, storage_size * alpha)|\label{line:storage_pn_objective}|
            
# Loop through operations nodes and set variables, constraint, and objective
for (j, node) in enumerate(operation_nodes) |\label{line:storage_on_loop1}|
    @variable(node, 0 <= y_stored)
    @variable(node, 0 <= y_sell <= d_sell)
    @variable(node, -d_save <= y_save <= d_save)
    @variable(node, 0 <= x_buy <= d_buy)
            
    @constraint(node, y_save + y_sell - zeta * x_buy == 0)
    @objective(node, Min, x_buy * beta[j] - y_sell * gamma[j])
end |\label{line:storage_on_loop2}|

# Set initial storage level
@constraint(operation_nodes[1], operation_nodes[1][:y_stored] == y_bar) |\label{line:storage_initial}|
            
# Set mass balance on storage unit
@linkconstraint(graph, [i = 1:(T - 1)], operation_nodes[i + 1][:y_stored] - |\label{line:storage_mass_balance1}|
                        operation_nodes[i][:y_stored] == operation_nodes[i + 1][:y_save]) |\label{line:storage_mass_balance2}|
                        
# Link planning decision to operations decisiosn
@linkconstraint(graph, [i = 1:T], operation_nodes[i][:y_stored] <= planning_node[:storage_size]) |\label{line:storage_links}|
            
# Set graph objective
set_to_node_objectives(graph) |\label{line:storage_graph_obj}|
            
# Optimizer the graph
set_optimizer(graph, HiGHS.Optimizer) |\label{line:storage_set_optimizer}|
optimize!(graph)|\label{line:storage_optimize}|
        \end{lstlisting}
        \end{scriptsize}
    \end{minipage}
\end{figure}

With an understanding of basic {\tt Plasmo.jl} syntax, we will next show how hierarchical subgraphs can be constructed within {\tt Plasmo.jl}, as shown in Code Snippet \ref{code:hier_graphs}. This code creates the hierarchical graph structure shown in Figure \ref{fig:hier_example}. Since the example in Figure \ref{fig:hier_example} is not representing a specific optimization problem, the values in this code snippet in various constraints are arbitrary, and the code snippet is intended to create a problem with a structure like that of the figure, not to represent a specific case study. For simplicity, we first define a function called {\tt build\_graph} (Lines \ref{line:hier_graphs_func1} - \ref{line:hier_graphs_func2}) to create a basic graph structure with three nodes with edges between each node (this is the structure presented by each of the blue subgraphs shown in Figure \ref{fig:hier_example}). The subgraphs $\mg$, $\mg_1$, and $\mg_2$ from Figure \ref{fig:hier_example} are next defined as OptiGraphs named {\tt g}, {\tt g1}, and {\tt g2}, respectively (Line \ref{line:hier_graphs_high_graphs}). We then build $\mg_{1,1}$, $\mg_{1,2}$, $\mg_{2,1}$, $\mg_{2,2}$, $\mg_{2,3}$, and $\mg_{2,4}$ from Figure \ref{fig:hier_example} by calling the {\tt build\_graph} function (Lines \ref{line:hier_graphs_low_graphs1} - \ref{line:hier_graphs_low_graphs2}). Each subgraph is now created, but no hierarchical structure has been formed. Importantly, {\it no edges can yet be added between subgraphs since edges must be ``owned'' by a higher level graph}. We first add {\tt g11} and {\tt g12} to the higher-level OptiGraph {\tt g1} and then add {\tt g21}, {\tt g22}, {\tt g23}, and {\tt g24} to the higher-level OptiGraph {\tt g2} (Lines \ref{line:hier_graphs_add_subgraph1} - \ref{line:hier_graphs_add_subgraph2}). Edges can now be placed between the subgraphs {\tt g11} and {\tt g12} and between {\tt g21}, {\tt g22}, {\tt g23}, and {\tt g24} (Lines \ref{line:hier_graphs_links1} - \ref{line:hier_graphs_links2}), represented by the red colored edges shown in Figure \ref{fig:hier_example}. Finally, {\tt g1} and {\tt g2} are added to the main OptiGraph {\tt g} (Line \ref{line:hier_graphs_add_subgraph3}) and edges are placed between the nodes and subgraphs contained on {\tt g1} and {\tt g2} (Lines \ref{line:hier_graphs_links3} - \ref{line:hier_graphs_links4}), represented by the black colored edges shown in Figure \ref{fig:hier_example}. Since these edges link the mid-level subgraphs {\tt g1} and {\tt g2}, they must be owned by the highest-level graph {\tt g}.

\begin{figure}[!htp]
    \begin{minipage}[t]{0.8\linewidth}
        \begin{scriptsize}
        \lstset{language=Julia, breaklines = true}
        \begin{lstlisting}[label = code:hier_graphs, caption = {Code for generating the graph shown in Figure \ref{fig:hier_example} with one high-level subgraph, two mid-level subgraphs, and six low-level subgraphs.}] 
using Plasmo

# Define function for creating subgraphs with 3 nodes
function build_graph() |\label{line:hier_graphs_func1}|
    # Define OptiGraph
    g = OptiGraph()
            
    # Add nodes
    @optinode(g, nodes[1:3]) |\label{line:hier_graphs_nodes}|
            
    # Set variables on nodes
    @variable(nodes[1], x >= 0) |\label{line:hier_graphs_node_calls1}|
    @variable(nodes[2], x >= 0)
    @variable(nodes[3], x >= 0) |\label{line:hier_graphs_node_calls2}|
            
    # Set objective on nodes
    @objective(nodes[1], Min, 1 * nodes[1][:x])
    @objective(nodes[2], Min, 2 * nodes[2][:x])
    @objective(nodes[3], Min, 3 * nodes[3][:x])
            
    # Define link constraints between nodes
    @linkconstraint(g, nodes[1][:x] + nodes[2][:x] >= 1)
    @linkconstraint(g, nodes[1][:x] + nodes[3][:x] >= 1)
    @linkconstraint(g, nodes[2][:x] + nodes[3][:x] >= 1)

    set_to_node_objectives(g)
            
    # Return graph
    return g
end  |\label{line:hier_graphs_func2}|
            
# Define high-level and mid-level OptiGraphs
g = OptiGraph(); g1 = OptiGraph(); g2 = OptiGraph();  |\label{line:hier_graphs_high_graphs}|
            
# Define lower-level OptiGraphs that contain nodes and edges
g11 = build_graph(); g12 = build_graph();  |\label{line:hier_graphs_low_graphs1}|
g21 = build_graph(); g22 = build_graph();
g23 = build_graph(); g24 = build_graph(); |\label{line:hier_graphs_low_graphs2}|
            
# Add lower-level subgraphs to mid-level
add_subgraph!(g1, g11); add_subgraph!(g1, g12); |\label{line:hier_graphs_add_subgraph1}|
add_subgraph!(g2, g21); add_subgraph!(g1, g22);
add_subgraph!(g2, g23); add_subgraph!(g1, g24); |\label{line:hier_graphs_add_subgraph2}|
            
# Define link constraints between lower-level subgraphs;
# these constraints are "owned" by the mid-level graphs
@linkconstraint(g1, g11[:nodes][1][:x] + g12[:nodes][1][:x] >= 1) |\label{line:hier_graphs_links1}|
@linkconstraint(g1, g11[:nodes][3][:x] + g12[:nodes][2][:x] >= 1)
            
@linkconstraint(g2, g21[:nodes][1][:x] + g22[:nodes][1][:x] >= 1)
@linkconstraint(g2, g21[:nodes][3][:x] + g22[:nodes][1][:x] >= 1)
            
@linkconstraint(g2, g22[:nodes][1][:x] + g23[:nodes][1][:x] >= 1)
@linkconstraint(g2, g22[:nodes][3][:x] + g23[:nodes][1][:x] >= 1)
            
@linkconstraint(g2, g23[:nodes][1][:x] + g24[:nodes][1][:x] >= 1)
@linkconstraint(g2, g23[:nodes][3][:x] + g24[:nodes][1][:x] >= 1) |\label{line:hier_graphs_links2}|
            
# Add mid-level subgraphs to high-level OptiGraph
add_subgraph!(g, g1); add_subgraph!(g, g2); |\label{line:hier_graphs_add_subgraph3}|

@linkconstraint(g, g11[:nodes][1][:x] + g21[:nodes][1][:x] >= 1) |\label{line:hier_graphs_links3}|
@linkconstraint(g, g11[:nodes][1][:x] + g22[:nodes][1][:x] >= 1)
@linkconstraint(g, g12[:nodes][1][:x] + g23[:nodes][1][:x] >= 1)
@linkconstraint(g, g12[:nodes][1][:x] + g24[:nodes][1][:x] >= 1) |\label{line:hier_graphs_links4}|

set_to_node_objectives(g)
set_to_node_objectives(g1)
set_to_node_objectives(g2)
\end{lstlisting}
\end{scriptsize}
\end{minipage}
\end{figure}
\vspace{0.1in}

An important feature of the OptiGraphs is that they can be treated as independent optimization problems. OptiGraphs can also be passed to {\tt JuMP.jl}'s {\tt optimize!} function, independent of other OptiGraphs. For instance, {\tt g1} can be optimized independent of {\tt g} or {\tt g2}. Further, the optimization problem on {\tt g1} would consist of all variables and constraints contained on the six nodes of {\tt g11} and {\tt g12} as well as the constraints stored on the red edges connecting the two subgraphs (Figure \ref{fig:hier_example}). Similarly, {\tt g11} and {\tt g12} could likewise be optimized independently of {\tt g1} by calling {\tt optimize!(g11, Ipopt.Optimizer)} and {\tt optimize!(g12, Ipopt.Optimizer)}; however, the constraints on the red edges would not be enforced by these {\tt optimize!} calls since these constraints are only enforced on {\tt g1} (or {\tt g} if {\tt g} was optimized). This also means that any OptiGraph being optimized must have a solver set on them (e.g., calling {\tt set\_optimizer(g1, Ipopt.Optimizer)} does not set a solver on {\tt g11} or {\tt g12}). These capabilities highlight the high degree of flexibility that is available to optimize different graph components using different solvers.

It is important to note some syntax in Code Snippet \ref{code:hier_graphs}. First, nodes, like variables, can be indexed by an array. Line \ref{line:hier_graphs_nodes} shows how {\tt @optinode(g, nodes[1:3])} forms 3 OptiNodes, with indices 1, 2, and 3, which can be called directly (e.g., Lines \ref{line:hier_graphs_node_calls1} - \ref{line:hier_graphs_node_calls2}). Variables on nodes can be accessed by calling the node and then the symbol for the variable (e.g., {\tt nodes[1][:x]} for variable {\tt x}). Similar to the symbol notation used to access variables on nodes, nodes on a graph can be accessed in a similar way (e.g., Line \ref{line:hier_graphs_links1}). For instance, the variable {\tt x} on {\tt nodes[1]} of {\tt g11} can be accessed by calling {\tt g11[:nodes][1][:x]}. In addition, there may be times where subgraphs or nodes do not have a global identifier defined. The API functions {\tt local\_subgraphs} and {\tt local\_nodes} can be used to find the subgraphs or nodes contained directly on an OptiGraph, while {\tt all\_subgraphs} and {\tt all\_nodes} return {\it all} subgraphs or nodes contained on an Optigraph and all its subgraphs. For example {\tt local\_subgraphs(g)} will return a vector containing {\tt g1} and {\tt g2} and {\tt local\_nodes(g)} will return an empty vector (since no nodes have defined directly on {\tt g}). The function {\tt all\_subgraphs(g)} will instead return a vector of all subgraphs, including {\tt g1}, {\tt g2}, {\tt g11}, {\tt g12}, {\tt g21}, {\tt g22}, {\tt g23}, and {\tt g24}, while {\tt all\_nodes(g)} will return every node, including those in the lower-level subgraphs. 

\subsection{Graph Structure Manipulation}

The OptiGraph provides flexibility for manipulating graph structure, including apply partitioning and aggregation schemes. These can facilitate visualization (e.g., plotting an aggregated graph with fewer nodes may be significantly easier than the original graph) and can help apply algorithms that are suitable for specific graph structures. The ideas below were first introduced in \cite{jalving2022graph}, but we formalize these ideas for the hierarchical setting, and we show later how they can be used for exploring and decomposing hierarchical structures. We will introduce the concepts mathematically here, and details on the software implementation are contained in the Supporting Information. 
\\

Problem structure can be manipulated by creating induced subgraphs on an existing graph, where an $``$induced subgraph$"$ is one in which a subset of edges contains all edges of the original graph that connect nodes of the subset (an edge of the induced subgraph cannot connect nodes not contained in the subset) \cite{bosak1984induced}. Here, we assume that a subgraph can be created from an existing graph $\mg_{\circ}$ to obtain a graph with nested subgraphs that we denote as $\mg_{\circ}(\{\mg_{\circ, i}\}_{i \in 1, ..., N_*}, \mn_{\mg_{\circ}}, \me_{\mg_{\circ}})$.
\\

Partitioning is a common type of graph manipulation that involves placing graph nodes into subgraphs \cite{bichot2013graph}. For an initial graph $\mg_{\circ}(\emptyset, \mn, \me)$, partitioning involves splitting the node set $\mn$ into $N_{part}$ subsets; we define the node subsets as $\mn_{i}^{part} \subseteq \mn, i \in \{1, 2, ..., N_{part}\}$ and this satisfies $\mn^{part}_{i} \cap \mn^{part}_{j} = \emptyset, \; i, j \in \{ 1, 2, ..., N_{part}\}, i \ne j$ and $\bigcup_{i \in \{1, 2, ..., N_{part}\}} \mn^{part}_i = \mn$. Each node subset is used to induce a subgraph, which allow us to express the original graph as $\mg_{\circ}(\{\mg_{\circ, i}\}_{i \in 1, ..., N_{part}}, \emptyset, \me_{\mg_{\circ}})$; here, $\me_{\mg_{\circ}}$ is the set of all edges in $\me$ not contained in the subgraphs. 
\\

Aggregation is another common type of graph manipulation; we define the initial graph $\mg_{\circ}(\emptyset, \mn, \me)$ and a set of nodes to be aggregated (denoted as $\mn_{agg} \subseteq \mn$). The edges connecting the aggregated nodes are given by $\me_{agg} = \{ e : e \in \me, n \in \mn_{agg}, \; n \in \mn(e)\}$. Aggregation reduces nodes in the subset $\mn_{agg}$ into a single node, which we denote as $n_{agg}$. The graph after aggregation is given by $\mg_{\circ}^{agg}(\emptyset, \bar{\mn}, \bar{\me})$ where $\bar{\mn}=\{n_{agg} \cup \mn \} \setminus \mn_{agg}$ is the set of nodes in the aggregated graph and $\bar{\me}$ is the set of edges in the aggregated graph. The new node $n_{agg}$ will contain all variables and constraints of the nodes in $\mn_{agg}$, and the objective on $n_{agg}$ will be a summation of all objectives on nodes in $\mn_{agg}$. All constraints contained on $\me_{agg}$ will also now be stored as local constraints on $n_{agg}$. Note that this aggregation approach does not change the optimization problem as all variables/constraints still exist (and the graph can contain its same objective function); instead, only the structure of the graph is changed. 
\\

\begin{figure*}[!htp]
    \centering
    \includegraphics[width = 1.0\textwidth]{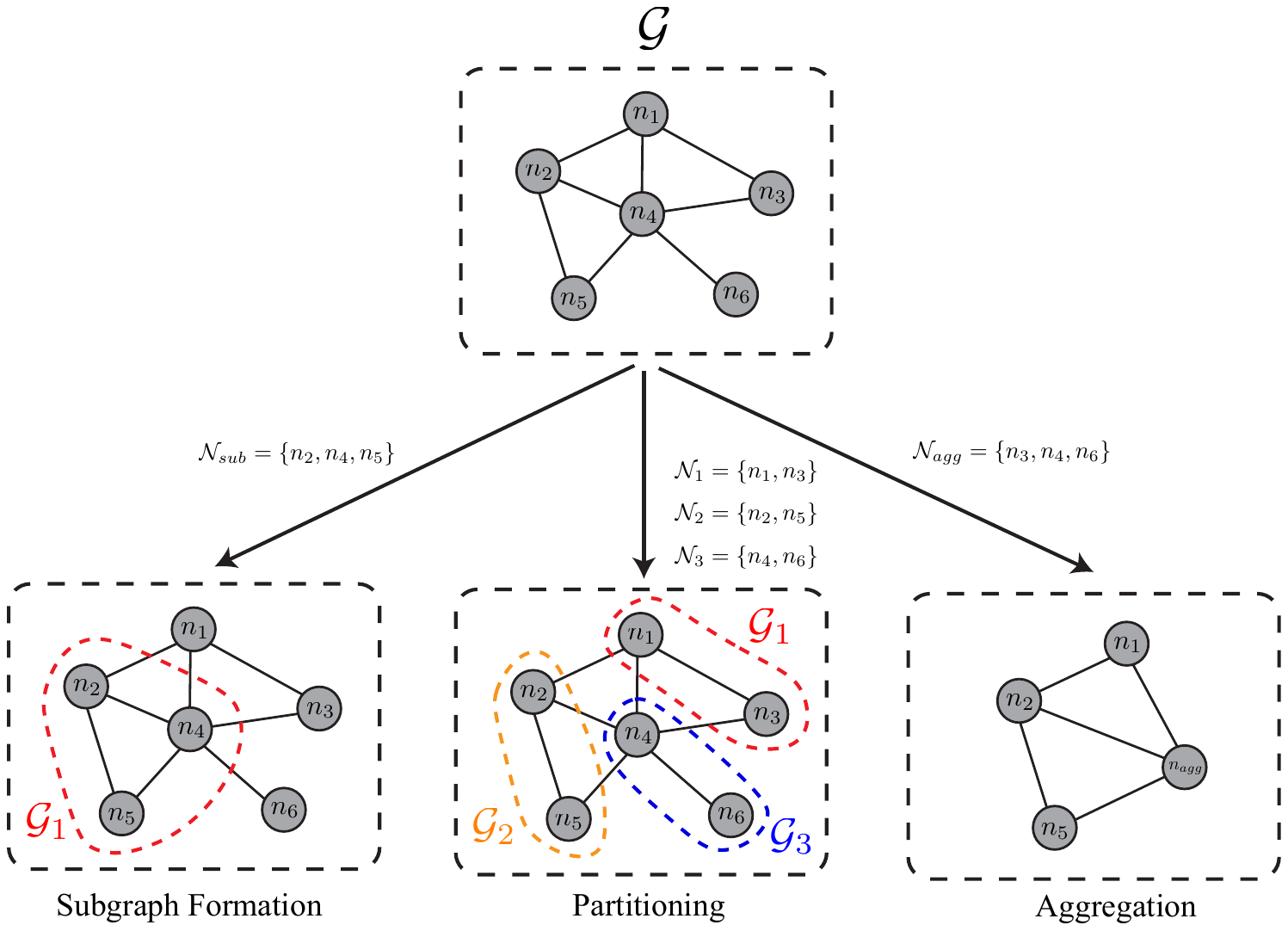}
        \vspace{-0.2in}
    \caption{Examples of how one can generate a subgraph from a graph $\mg$ using the node set $\mn_{sub}$, partitioning a graph $\mg$ with three partitions using the node sets $\mn_1$, $\mn_2$, and $\mn_3$, and aggregating nodes in $\mn_{agg}$ of a graph $\mg$ into a single new node. Each of these are manipulations of the graph or of its hierarchical structure.}
    \label{fig:structure_manipulation}
\end{figure*}

The graph manipulation approaches are visualized in Figure \ref{fig:structure_manipulation}; here, we show how the same graph can be manipulated by inducing subgraphs and partitions or by aggregating sets of nodes. By inducing subgraphs or adding partitions, we are creating a structure through the subgraphs without changing any of the nodes. This structure could have different characteristics than the original graph before inducing the subgraphs. Similarly, partitioning results in a new structure to the graph. These approaches can be useful for creating graphs with different characteristics; for instance, if a graph contains a single cycle, all nodes in the cycle could be aggregated into one single node, resulting in an acyclic graph.

The ability to manipulate structure has several applications. These are useful from a modeling perspective because they can simplify the formation of subproblems. For instance, in Sections \ref{sec:cs2} and \ref{sec:cs3}, a graph is modeled with no subgraphs and then partitioned into subproblem subgraphs for applying the graph-based Benders Decomposition. In \ref{sec:cs2}, this ability allows for exploring different formations of the master/root problem, while in \ref{sec:cs3} it is used to form subproblems of differing temporal sizes. Both of these approaches have important implications for run time and convergence of the algorithms. In addition, these structural manipulations can allow for visualizing different structures (e.g., visualizing temporal or spatial partitionings in a graph). 

\section{Graph-Based Decomposition}\label{sec:Decomposition}

Optimization problems with hierarchical structure can often be decomposed to enable scalable solutions or the derivation of approximate solutions. The OptiGraph abstraction can be used as a basis for defining such solution schemes and can simplify a practitioner's use of such solution algorithms. In the discussion that follows, we introduce a graph-based Benders Decomposition (gBD), which is a generalization of the nested Benders Decomposition (BD) \cite{kumar2021dual,lara2018deterministic,pereira1985stochastic,pereira1991multi,yagi2024nested,zou2019stochastic} to handle general hierarchical graph structures. We show how manipulations of the OptiGraph can be used to facilitate the derivation of gBD and introduce a {\tt Julia} package {\tt PlasmoBenders.jl} that implements gBD.

\subsection{Benders Decomposition}

\begin{figure*}[!ht]
    \centering
    \includegraphics[width = \textwidth]{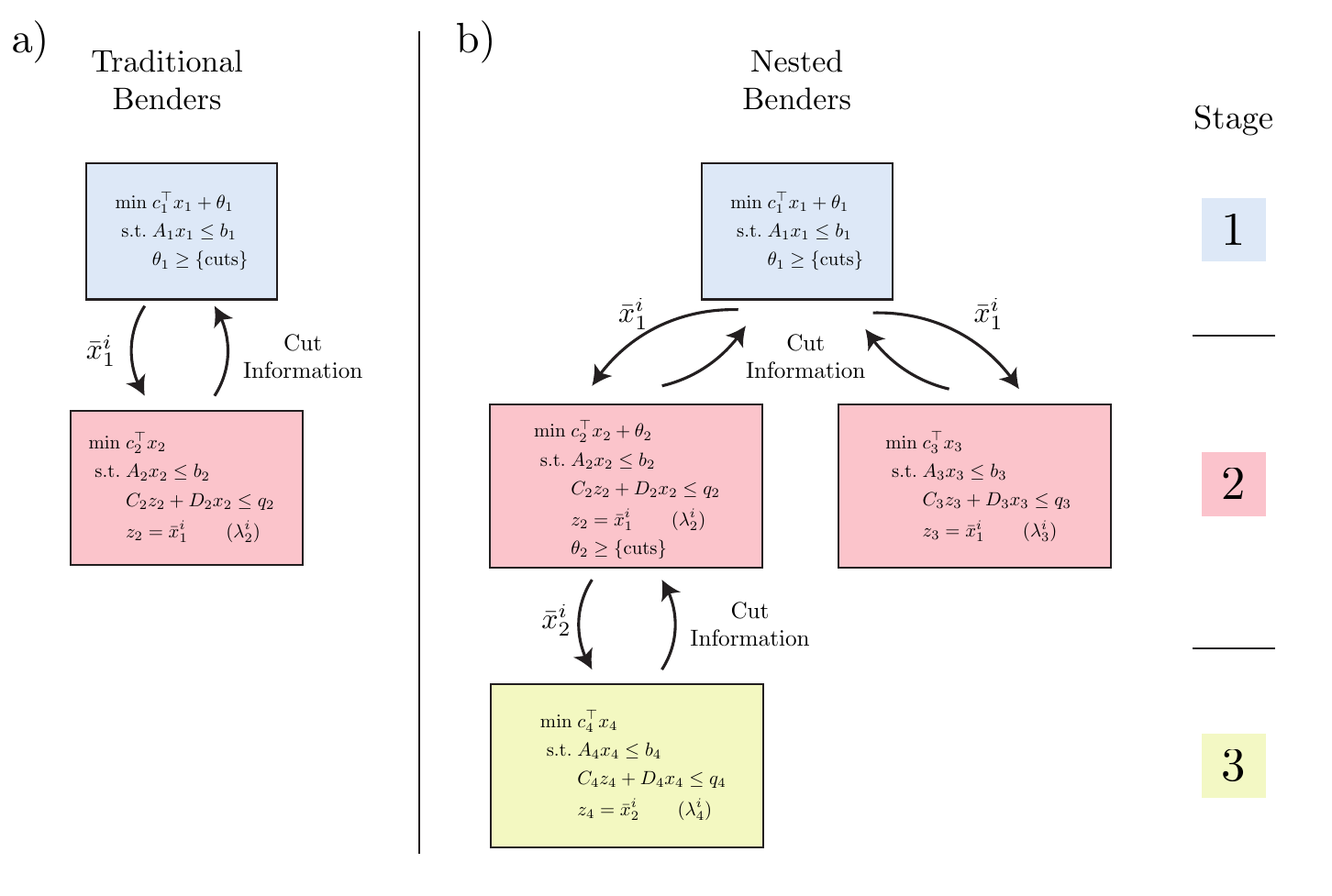}
        \vspace{-0.0in}
    \caption{Example of the Benders decomposition algorithm, with the traditional Benders case on the left and the nested Benders case on the right. Primal information is shared in a forward pass (the downward direction) and information for cuts is shared in a backward pass (the upward direction).}
    \label{fig:BD_subproblems}
\end{figure*}

BD is a powerful decomposition approach for solving linear and mixed-integer linear programs \cite{benders1962partitioning,rahmaniani2017benders}. The most general form of BD (also referred to as nested BD or dual dynamic programming) is applied to problems with tree graph structures. Here, the problem is partitioned into multiple, sequential ``stages'' (each stage has a corresponding optimization subproblem);  a visualization of this algorithm is given in Figure \ref{fig:BD_subproblems}. The general approach of BD is to solve each stage of a problem sequentially in a forward pass, passing the primal solution of the previous stage to the next stage, and then performing a backward pass, where information (often dual values) is passed from each stage to the previous stage to create cutting planes. Each pass of cutting planes helps to form a better lower bound on the problem. An iteration of BD consists of a forward pass and a backward pass, which produce an upper and lower bound, respectively, and the iterations continue until the upper and lower bounds converge or until another termination criterion is reached. The traditional BD algorithm consists of two stages (Figure \ref{fig:BD_subproblems}a), the first containing a master problem and the second containing a subproblem. However, BD can be applied to problems with multiple stages; see Figure \ref{fig:BD_subproblems}b. 

BD is applied to problems that can be separated into ``stages'' and ``subproblems''. By ``stages'' we are referring to a set of subproblems that depend on ``complicating variables'' (i.e., variables that are linked by constraints to the subproblem's variables) from the previous stage. These subproblems typically become easier to solve once those complicated variables are fixed. The form of these problems is given by
\begin{subequations}\label{eq:NBD}
    \begin{align}
        \min &\; \sum_{s=1}^S \left( \sum_{w \in W_s} c_w^\top x_w \right) \label{eq:NBD_obj}\\
        \textrm{s.t.} &\; A_w x_w \le b_w, \quad  w \in W_s, s \in \{1,\ldots,S\}, \label{eq:NBD_local_cons}\\
        &\; C_v x_w + D_v x_v \le q_v, \quad v \in \mc(w), w \in W_{s}, s \in \{1,\ldots,S-1\}\label{eq:NBD_links}
    \end{align}
\end{subequations}
\noindent where $\ms = \{1,\ldots,S\}$ is the set of stages (with $S$ the total number of stages); $W_s$ is the set of subproblem indices at stage $s\in \ms$; $x_w, w \in W_s$ are the decision variables (can be continuous or integer); $c_w, A_w, b_w, C_w, D_w$ and $q_w, w \in W_s$ are vector and matrix data; and $\mc(w), w \in W_s,  s\in\{1,\ldots,S-1\}$ is the index set of child subproblems for a subproblem $w$. By child subproblems, we are referring to the subproblems connected by the linking constraints \eqref{eq:NBD_links} that are in the next stage $s+1$ (note that a subproblem at stage $s+1$ cannot be a child of multiple subproblems in the previous stage). In other words, $\mc(w)$ is the set of all $v \in W_s, s \in \{2, ..., S\}$ where there are nonzero matrices $C_v$ and $D_v$ being multiplied by $x_w$ and $x_v$, respectively. Hence, for child subproblems of a subproblem $w$, variables $x_w$ can be considered the complicating variables (note that $S=2$ for the traditional BD case). Further, we note that \eqref{eq:NBD_links} is written as inequality constraints, but this formulation could easily be extended to equality constraints, since any equality constraint can be written with two inequality constraints. 

To solve \eqref{eq:NBD}, the BD algorithm partitions the problem into a root (or master) problem and subsequent subproblems across $S-1$ stages. The root (first stage) problem has the form
\begin{subequations}\label{eq:NBD_root}
    \begin{align}
        \underline{\phi}^i_r := \min &\; c_r^\top x_r + \theta_r\\
        \textrm{s.t.} &\; A_r x_r \le b_r \\
        &\; \theta_r \ge \{ \textrm{cuts} \}. \label{eq:NBD_root_cuts}
    \end{align}
\end{subequations}
\noindent where $r \in W_1$ and $W_1$ is a singleton, $\underline{\phi}_r^i$ is the optimal value of \eqref{eq:NBD_root}, $i$ is the iteration number, and  $\theta_r$ is the cost-to-go variable which underestimates the optimal value of the linked subproblems in the next stage. To approximate the optimal value of the subproblem, the root problem uses cutting planes (i.e., cuts) to constrain variable $\theta_r$. These cuts are based on information from the solution of the next stage subproblems and are said to be valid if they underestimate the optimal value of the subproblems. Note that rather than using a single cost-to-go variable, $\theta_r$ could be replaced in the objective by $\Sigma_{w \in W_2} \theta_{r, w}$, and constraint \eqref{eq:NBD_root_cuts} could be replaced by constraints on each individual $\theta_{r,w}$, where $\theta_{r,w}$ is a cost-to-go variable which underestimates the optimal value of individual subproblems at the subsequent stage (rather than an aggregation of subproblems at the subsequent stage). These two approaches are referred to as an aggregated cut when using $\theta_r$ and multi-cuts when using $\theta_{r,w}, w \in W_2$, and both can have their own benefits \cite{birge1997introduction,birge1988multicut,zverovich2012computational}. The multi-cut approach adds cuts after each iteration of the algorithm, and can be more efficient for some problems, but the aggregated approach results in fewer constraints being added to the root problem at each iteration, which can be important if there are already many constraints on the root problem or if $|W_2|$ is large. 

Each subsequent subproblem  $w \in W_s, s \in \{2, ..., S-1 \}$ with $\mc(w) \ne \emptyset$, has the form 
\begin{subequations}\label{eq:NBD_subroot}
    \begin{align}
        \overline{\phi}^i_w(\bar{x}^i_{p_w}) := \min &\; c_w^\top x_w + \theta_w\\
        \textrm{s.t.} &\; A_w x_w \le b_w \label{eq:NBD_subroot_cons} \\
        &\; C_w z_w + D_w x_w \le q_w \label{eq:NBD_link_constraint1}\\
        &\; z_w = \bar{x}_{p_w}^i \quad (\lambda^i_{w}) \label{eq:NBD_subroot_links}\\
        &\; \theta_w \ge \{ \textrm{cuts} \}. \label{eq:NBD_subroot_cuts}
    \end{align}
\end{subequations}
\noindent where $p_w$ is the index of the parent subproblem of $w$ (i.e., $w \in \mc(p_w)$), $\bar{x}_{p_w}$ is a feasible solution of \eqref{eq:NBD_subroot}, and where $\overline{\phi}_w$ is the optimal value of \eqref{eq:NBD_subroot}. Note that, unlike \eqref{eq:NBD}, which is written in terms of the parent subproblem, \eqref{eq:NBD_subroot} is written in terms of the child subproblem. As with the root problem, the cuts in \eqref{eq:NBD_subroot} could be of the form of aggregated cuts or multi-cuts. A feasible solution of the previous stage subproblem, $\bar{x}_{p_w}$, is also passed to \eqref{eq:NBD_subroot}, and the constraints linking the two subproblems, \eqref{eq:NBD_links}, are also enforced in \eqref{eq:NBD_link_constraint1}. The final subproblems (i.e., subproblems $w \in W_s, s \in \{2,\ldots,S\},$ with $\mc(w) = \emptyset$) have the form
\begin{subequations}\label{eq:NBD_final_sub}
    \begin{align}
    \overline{\phi}^i_w(\bar{x}^i_{p_w}) := \min &\; c_w^\top x_w\\
        \textrm{s.t.} &\; A_w x_w \le b_w \label{eq:NBD_final_sub_cons} \\
        &\; C_w z_w + D_w x_w \le q_w \label{eq:NBD_link_constraint2}\\
        &\; z_w = \bar{x}_{p_w}^i \quad (\lambda^i_{w}). \label{eq:NBD_final_sub_links}
    \end{align}
\end{subequations}
\noindent where \eqref{eq:NBD_final_sub} has a similar form to \eqref{eq:NBD_subroot} but where there is no cost-to-go variable since there are no children subproblems in future stages. We also note that at any given stage, not all subproblems of that stage may have children subproblems.

With the root and subproblems defined, we can now define the upper and lower bounds. Assuming relatively complete recourse of the subproblems (this can be ensured by adding slack variables to constraints and penalizing these slacks in the objectives), the upper bound at the $k$th iteration is given by 
\begin{equation}\label{eq:NBD_UB}
    UB^k_{BD} := \min_{i = 1,2, ..., k} \left( \sum_{s = 1}^S \sum_{w \in W_s} c_w^\top \bar{x}^i_w \right).
\end{equation}
\noindent The cost-to-go variables $\theta_w$ are an underestimate of the objective value of the subsequent subproblems, and if the cuts in \eqref{eq:NBD_root_cuts} or \eqref{eq:NBD_subroot_cuts} are valid, then at iteration $i$, $\sum_{w \in W_s} \theta_w \le \sum_{w \in W_{s+1}} \overline{\phi}^i_w$ (or for multicuts, $\theta_{w, v} \le \overline{\phi}^i_{v}, w \in W_s, v \in \mc(w)$). Thus, since cost-to-go variables in any stage are underestimating the next stage optimal value through every stage of the problem, the lower bound is given by the optimal solution of the root-problem, $LB^k_{BD} := \underline{\phi}^k_r$. We will omit the subscript $BD$ when it is obvious from context.

The BD algorithm is summarized in Algorithm \ref{alg:NBD}. The maximum number of iterations, $K_{max}$, the tolerance $\epsilon_{tol}$, and any cut data are set. An iteration of the algorithm consists of first solving the root subproblem, \eqref{eq:NBD_root}, which gives a lower bound on the solution. Then, each stage is sequentially solved (including all subproblems in that stage) with feasible solutions of the previous stage fixed in the current stage. The feasible solutions that are passed to the next stage are not required to be the optimal solutions of the previous stage; in some cases, other feasible solutions can be chosen, such as by using a regularization procedure \cite{lemarechal1995new,ruszczynski1986regularized,pecci2024regularized}, which can yield better cuts. After each subproblem is solved, the upper bound is computed. If the gap between the upper and lower bounds is below $\epsilon_{tol}$, the solutions at the iteration that yielded the upper bound are returned. If not, new cuts are computed and added to the respective root and subproblems, and the process is repeated. Note that the initial pass through the stages is often considered a forward pass because each stage is solved sequentially. The computation of new cuts is often considered a backward pass because it can likewise be a sequential process starting at the last stage and proceeding backward to the first, adding a cut to the previous stage before resolving it.
If the number of iterations reaches $K_{max}$ without reaching the solver tolerance, the algorithm terminates.

\begin{algorithm}
    \caption{Benders Decomposition}\label{alg:NBD}
    \begin{algorithmic}
    \State Set $K_{max}$, $\epsilon_{tol}$, $k = 1$, and cut data
    \While{$k \le K_{max}$}
    \State Solve \eqref{eq:NBD_root} to get $\underline{\phi}^k_r$ and set $LB^k = \underline{\phi}^k_r$
    \State Determine root-problem iterates $\bar{x}^k_r$
    \For{$s \in \{2, ..., S \}$}
        \For{$w \in W_s$}
            \State Solve \eqref{eq:NBD_subroot} and \eqref{eq:NBD_final_sub} for $w$ to get $\bar{\phi}^k_w$
            \State If $s<S$, determine decisions $\bar{x}^k_w$ to use in next stage subproblems.
        \EndFor
    \EndFor
    \State Compute $UB^k$ from \eqref{eq:NBD_UB} 
    \If{ $(UB^k - LB^k) / |LB^k| \le \epsilon_{tol}$}
        \State Set $j =\displaystyle \textrm{argmin}_{i = 1, ..., k} \left( \sum_{s =1}^S \sum_{w \in W_s} c_w^\top \bar{x}^i_w \right) $
        \State Return solutions $\bar{x}^j_w, w \in W_s, s=1,\ldots,S$ and $UB^j$
        \State Stop
    \Else
        \For{$s \in {S -1, ..., 1}$}
            \For {$w \in W_s$}
            \State Compute new cut(s) and add them to \eqref{eq:NBD_root_cuts} or  \eqref{eq:NBD_subroot_cuts}
            \State Solve \eqref{eq:NBD_root} or \eqref{eq:NBD_subroot} for $w$
            \EndFor
        \EndFor
        \State $k = k + 1$
    \EndIf
    \EndWhile
    \If{$(UB^k - LB^k) / |LB^k| > \epsilon_{tol}$}
        \State Set termination status to {\tt Maximum Number of Iterations Exceeded}.
    \EndIf
    \end{algorithmic}
\end{algorithm}
\vspace{0.1in}

BD is a useful decomposition approach, but it has some limitations. It is best suited for large and possibly intractable problems since BD results in smaller problems that are iteratively solved. However, BD presents a challenge in that it may have several stages that must be solved sequentially, which can limit the parallelization. One method for addressing this is to place multiple subproblems in a given stage, as each of these subproblems can be solved in parallel \cite{jacobson2024computationally}. The choice of which variables belong in the first stage (root) problem is an important modeling decision since placing linking variables in this stage can result in more subproblems in the second stage, each of which could be solved in parallel. 

There are some extensions of BD that endeavor to parallelize more of the computations for problems with many stages \cite{colonetti2022parallel,rahmaniani2024asynchronous,dos2023accelerated,santos2016new}, but the sequential nature of BD can be a shortcoming for some problem instances. Consequently, having a modeling framework that allows for exploring various subproblem structuring (e.g., which variables belong to which problem, how many stages or subproblems per stage to use) can be important for applying BD in an efficient manner. For instance, one can imagine having an intractable monolithic problem that they structure into a multi-stage BD problem with one subproblem per stage (e.g., as in the structure seen in \cite{lara2018deterministic,yagi2024nested}). One could then further structure the problem into a two-stage BD problem (e.g., as in the structure seen in \cite{pecci2024regularized}) by placing all linking variables between stages into a root or master problem, making each subproblem independent and parallelizable in the second stage.

\subsection{Graph-Based Benders Decomposition}

To define gBD, we identify basic structures and syntax that we will use to represent the decomposed subproblems (e.g., the root and subproblems of BD). So long as subgraphs do not overlap, subgraphs and nodes can be treated as interchangeable since each encodes an optimization problem. Thus, an algorithm generalized for an OptiGraph could use either subgraphs or nodes as subproblems. We will write the graph-based algorithms in terms of subgraphs rather than nodes because subgraphs can capture hierarchy and because {\tt Plasmo.jl} is designed to be more subgraph-centric. 

gBD requires a graph of the form $\mg(\{\mg_i\}_{i \in \{ 1, ..., N\}}, \emptyset, \me_\mg)$ where no subgraphs in $\mg$ overlap. The graph $\mg$ must also meet the following criteria, where we introduce the notation $\msg(e)$ as the set of subgraphs connected by an edge $e$:
\begin{enumerate}
    \item $|\msg(e)| = 2,  e \in \me_\mg$. i.e., any edge connecting the subgraphs $\msg(\mg)$ only connects two subgraphs. This condition is required because $\mg$ is a hypergraph and therefore a single edge has the capability of connecting more than two subgraphs.
    \item $\bigcup_{e \in \me_\mg} \msg(e) = \msg(\mg)$. i.e., each subgraph in $\mg$ must be connected to at least one other subgraph in $\mg$ so that there are no disconnected subgraphs.
    \item The subgraphs form a tree structure. That is, for any subgraph  $g_i \in \msg(\mg)$, it is not possible to move along edges to another subgraph $g_j \in \msg(\mg)$ and return to $g_i$ without traversing across the same set of subgraphs again. In other words, there are no cycles among any subgraphs $\msg(\mg)$. 
\end{enumerate}
\noindent Some examples of the required structure are visualized in Figure \ref{eq:NBD_graph}, where each subgraph (i.e., subproblem) is represented by a node and where each depicted edge could contain many edges between OptiNodes in the different subgraphs. The subgraph structure for each graph is an acyclic (tree) graph. The first graph shows the traditional BD case where there is a root problem ($\mg_1$) and a subproblem ($\mg_2$). The second graph is a BD case where the subproblem subgraphs ($\mg_2$ - $\mg_n$) can be solved in parallel since there are no links between the second stage subgraphs. The third graph is a linear graph; these graphs commonly arise in dynamic optimization, where each time point is constrained by the previous and next time points. The final graph is a more complex tree graph, which includes multiple subgraphs in each stage. In all three cases, we treat $\mg_1$ as the ``root'' or first-stage subgraph. However, any subgraph could be chosen as the root subgraph (this could change the number of stages), and the choice of root subgraph can have important algorithmic implications. We will notate the root subgraph as $\mg_r$.

\begin{figure*}[!htp]
    \centering
    \includegraphics[width = 1.0\textwidth]{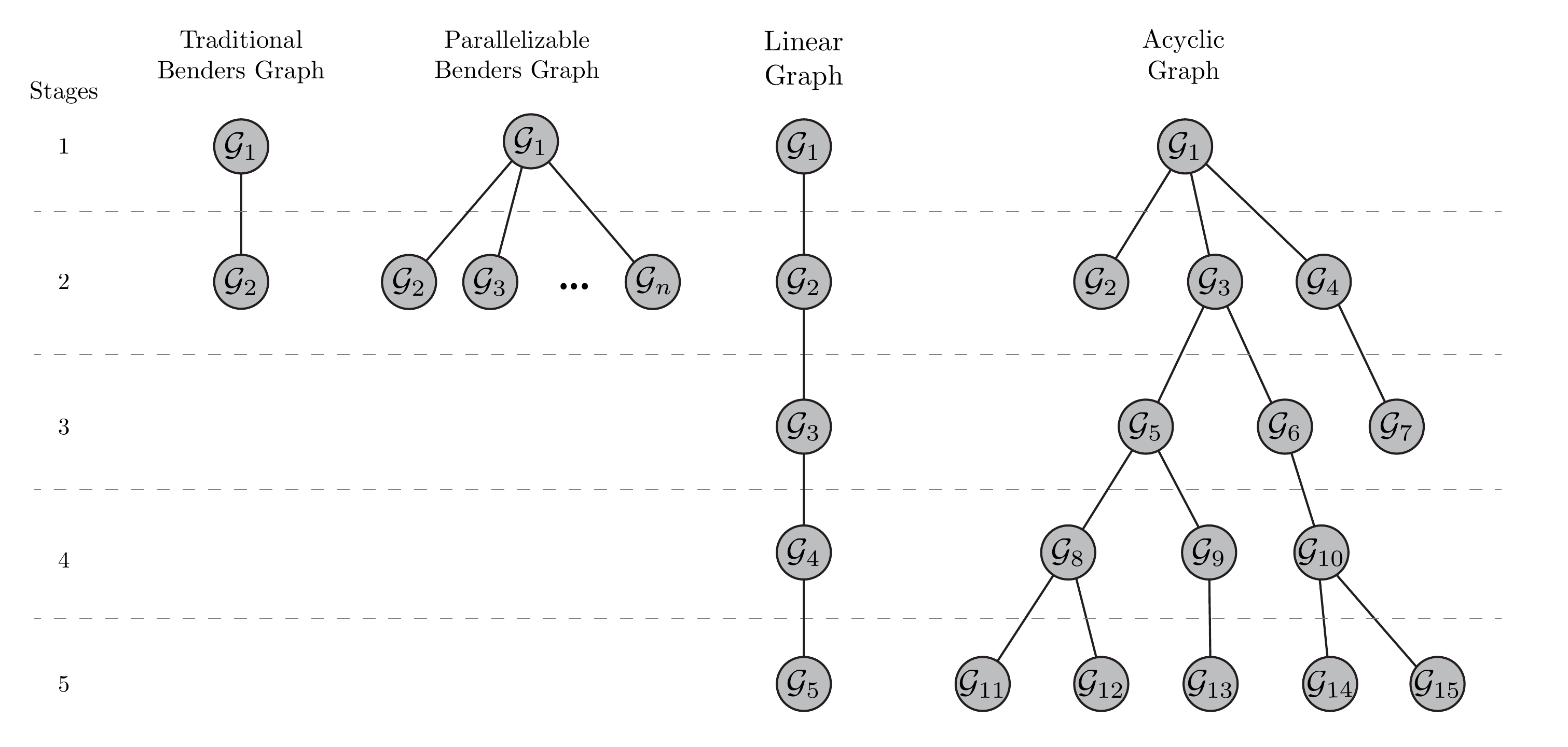}
        \vspace{-0.0in}
    \caption{Examples of the graph structure required for gBD, including a traditional BD case, a BD case with parallelizable subproblems, a linear graph, and a more complex acyclic graph. Here, subgraphs are represented by circles (nodes) but can embed more structure, nodes, and complexity. The root subgraph, $\mg_r$ is set to be $\mg_1$. Once the root subgraph is set, the remainder of the subgraphs are placed into stages.}
    \label{fig:NBD_graph}
\end{figure*}

We will also define the mathematical form of the problems defined on the subgraphs and edges. Each subgraph $g \in \msg(\mg)$ corresponds to the optimization problem: 
\begin{subequations}\label{eq:NBD_graph}
    \begin{align}
        \min &\; c_{g}^\top x_{g} \label{eq:GNBD_obj}\\
        \textrm{s.t.} &\; A_{g} x_{g} \le b_{g}. \label{eq:GNBD_local_cons}
    \end{align}
\end{subequations}
Each edge $e \in \me_\mg$ contains constraints of the form
\begin{align}\label{eq:benders_edge}
    \sum_{g \in \msg(e)} F_{g, e} x_g \le q_e
\end{align}
\noindent where $|\msg(e)| = 2$ by the definition of the gBD structure and where $F_{g,e}$ is a matrix. Note that the optimization problem on the overall OptiGraph $\mg$ has the form of \eqref{eq:NBD}, where \eqref{eq:NBD_obj} is represented by the summation of the subgraphs objectives, \eqref{eq:GNBD_obj}, the constraints \eqref{eq:NBD_local_cons} are represented by the constraints \eqref{eq:GNBD_local_cons}, and the linking constraints \eqref{eq:NBD_links} are represented by \eqref{eq:benders_edge}.

The optimization problem on each subgraph in \eqref{eq:NBD_graph} will be updated to create the subproblem structures given by \eqref{eq:NBD_root}, \eqref{eq:NBD_subroot}, and \eqref{eq:NBD_final_sub}. An important task in gBD is choosing the root subgraph, and any $g \in \msg(\mg)$ can be chosen. The root subgraph is updated to have the form
\begin{subequations}\label{eq:NBD_rootgraph}
    \begin{align}
        \underline{\phi}_{\mg_r} := \min &\; c_{\mg_r}^\top x_{\mg_r} + \theta_{\mg_r} \\
        \textrm{s.t.} &\; A_{\mg_r} x_{\mg_r} \le b_{\mg_r} \\
        &\; \theta_{\mg_r} \ge \{\textrm{cuts} \} \label{eq:GNBD_rootcuts}
    \end{align}
\end{subequations}
\noindent which matches the form of \eqref{eq:NBD_root}. The subsequent subgraphs $g \in \msg(\mg)\setminus \mg_r$ with $\mc(g) \ne \emptyset$ (where $\mc(g)$ are the child subgraphs or the subgraphs connected by an edge to $g$ in the next stage) take the form
\begin{subequations}\label{eq:GNBD_sub}
    \begin{align}
        \overline{\phi}^i_g := \min &\; c_g^\top x_g + \theta_g\\
        \textrm{s.t.} &\; A_g x_g \le b_g \label{eq:GNBD_sub_con1} \\
        &\; C_{p_g, e} z_g + C_{g, e} x_g \le q_e, \quad e \in \{e:e\in \me_\mg, g, p_g \in \msg(e)\} \label{eq:GNBD_sub_con2} \\
        &\; z_g = \bar{x}_{p_g}^i \qquad (\lambda_g^i) \label{eq:GNBD_sub_links} \\
        &\; \theta_g \ge \{ \textrm{cuts} \} \label{eq:GNBD_subroot_cuts}
    \end{align}
\end{subequations}
\noindent which matches the form of \eqref{eq:NBD_subroot} and where $p_g$ is the parent graph of $g$ (the graph connected by an edge to $g$ and which is in the previous stage). Finally, the subgraphs with a parent subgraph but no children subgraphs take the form 
\begin{subequations}\label{eq:GNBD_final_sub}
    \begin{align}
        \overline{\phi}^i_g := \min &\; c_g^\top x_g\\
        \textrm{s.t.} &\; A_g x_g \le b_g \label{eq:GNBD_final_sub_con1} \\
        &\; C_{p_g, e} z_g + C_{g, e} x_g \le q_e, \quad e \in \{e:e\in \me_\mg, g,p_g \in \msg(e)\} \label{eq:GNBD_final_sub_con2} \\
        &\; z_g = \bar{x}_{p_g}^i \qquad (\lambda_g^i) \label{eq:GNBD_final_links}
    \end{align}
\end{subequations}
\noindent which matches the form of equation \eqref{eq:NBD_final_sub}. In both \eqref{eq:GNBD_sub} and \eqref{eq:GNBD_final_sub}, the edge constraints \eqref{eq:benders_edge} are enforced as \eqref{eq:GNBD_sub_con2} and \eqref{eq:GNBD_final_sub_con2} in the subgraph that is downstream from the root subgraph. The upper bound at the $k$th iteration of gBD is given by 
\begin{align}\label{eq:GNBD_UB}
    UB^k_{gBD} := \min_{i = 1, 2, ..., k} \sum_{g \in \msg \left(\mg \right)} c_g^\top \bar{x}^i_w
\end{align}
\noindent where we will omit the subscript gBD when it is obvious from context. As with BD, the lower bound is given by $LB^k_{gBD} := \underline{\phi}^k_{\mg_r}$.

The gBD algorithm is summarized in Algorithm \ref{alg:graph_NBD}. A root subgraph, $\mg_r$ is defined for the graph $\mg$. The root subgraph is updated to take the form of \eqref{eq:NBD_rootgraph}. The graph then must be traveled through (traversed) to update each subgraph and note its stage. This is done by moving from the root subgraph to each stage of the new graph using two new sets: a queue set, $\msg_{queue}$ that records the subgraphs that need to be explored next, and a visited set, $\msg_{visited}$, which includes a set of all subgraphs that have been visited (i.e., subgraphs that have been updated to take the form of \eqref{eq:NBD_rootgraph}, \eqref{eq:GNBD_sub}, or \eqref{eq:GNBD_final_sub}). This is effectively a breadth-first search \cite{bundy1984breadth,lee1961algorithm}. We also define a function $f_{stage}$ that maps a subgraph to its stage. After the subgraphs have been updated, the maximum number of iterations, $K_{max}$, the tolerance, $\epsilon_{tol}$, and cut data are set, as well as the total number of stages. The graph is then sequentially solved, starting at the root subgraph (which also yields the lower bound), and then progressing through each stage, where the solutions of the previous stage are passed to the next stage and those solutions fixed as in \eqref{eq:GNBD_sub_links} and \eqref{eq:GNBD_final_links}. Note that each set of subgraphs at a given stage can be solved in parallel. After each subgraph is solved, the upper bound is computed. If the relative gap between the upper and lower bounds is below the tolerance, the algorithm terminates, and the solutions of the iteration that yielded the best lower bound are returned. If the relative gap is not below the tolerance, new cuts are added to their respective subproblems, and another iteration is run. 

\begin{algorithm}
    \caption{Graph Benders Decomposition}\label{alg:graph_NBD}
    \footnotesize
    \begin{algorithmic}
    \State Set OptiGraph $\mg(\{\mg_i\}_{i \in \{ 1, ..., N\}}, \emptyset, \me_\mg)$ and root subgraph $\mg_r \in \msg(\mg)$ 
    \Ensure $\mg$ meets the three criteria for graph structure for gBD
    \State Define root subgraph $\mg_r$
    \State Define a new set for subgraphs in queue, $\msg_{queue} := \mc(\mg_r)$
    \State Define a new set for subgraphs that have been visited, $\msg_{visited} := \{\mg_r\}$
    \State Define a function, $f_{stage}$, that maps a sugraph to its stage.
    \State Set $f_{stage}(\mg_r) = 1$
    \While{$|\msg_{visited}| \ne |\msg(\mg)|$}
        \For{$g \in \msg_{queue}$}
            \If{$\mc(g) \ne \emptyset$}
                \State Formulate problem~\eqref{eq:GNBD_sub} for $g$
                \State Set $\msg_{queue} \leftarrow \msg_{queue} \cup \mc(g)$
            \Else
                \State Formulate problem~\eqref{eq:GNBD_final_sub} for $g$
            \EndIf
            \State Set $f_{stage}(g) = f_{stage}(p_g) + 1$
            \State Set $\msg_{queue} \leftarrow \msg_{queue} \setminus g$
            \State Set $\msg_{visited} \leftarrow \msg_{visited} \cup \{ g\}$
        \EndFor
    \EndWhile
    \State Set $K_{max}$, $\epsilon_{tol}$, $k = 1$, and cut data, and set $S = \max_{g \in \msg(\mg)} f_{stage}(g)$
    \While{$k \le K_{max}$}
    \State Solve root subgraph to compute $\underline{\phi}_{\mg_r}^k$ and set $LB^k = \underline{\phi}^k_{\mg_r}$
    \State Determine root-problem decisions $\bar{x}^k_{\mg_r}$
    \For{$s \in \{2, ..., S\}$}
        \For{$g \in \{g: f_{stage}(g) = s, g \in \msg(\mg)\}$}
            \State Solve the corresponding subproblem and compute $\overline{\phi}_g^k$
            \State Determine iterates to use in next stage subproblems, $\bar{x}^k_g$
        \EndFor
    \EndFor
    \State Compute $UB^k$ from \eqref{eq:GNBD_UB} 
    \If{ $(UB^k - LB^k) / |LB^k| \le \epsilon_{tol}$}
        \State Set $j = \textrm{argmin}_{i = 1,2, ..., k} \sum_{g \in \msg(\mg)} c^\top_g \bar{x}^i_g $
        \State Return solutions $\bar{x}^j_g, g \in \msg(\mg)$ and $UB^k$
        \State Stop.
    \Else
        \For $s \in \{S-1, ..., 1\}$
            \For {$g \in \{g: f_{stage}(g) = s, g \in \msg(\mg)\}$}
            \State Compute new cut(s) and add them to \eqref{eq:GNBD_rootcuts} or  \eqref{eq:GNBD_subroot_cuts}
            \State Solve \eqref{eq:NBD_rootgraph} or \eqref{eq:GNBD_sub} for $g$
            \EndFor
        \EndFor
        \State $k = k + 1$
    \EndIf
    \EndWhile
    \If{$(UB^k - LB^k) / |LB^k| > \epsilon_{tol}$}
        \State Set termination status to {\tt Maximum Number of Iterations Exceeded}.
    \EndIf
    \end{algorithmic}
\end{algorithm}
\vspace{0.1in}

Several variants of nested BD exist, and Algorithm \ref{alg:graph_NBD} can be easily extended in several ways, many of which are implemented in our software implementation. For instance, Algorithm \ref{alg:graph_NBD} assumes an aggregated cut approach, but this can be easily extended to multi-cuts. In addition, the most common nested BD approach for the backwards pass is outlined in this algorithm, where the backwards pass sequentially solves stages in reverse order, starting at stage $S-1$ and proceeding backwards to stage $1$, passing a cut each time. In doing so, it allows information to propagate backwards each iteration. However, other variants of this approach exist, and the variant used in the third case study below is one in which the ``backwards pass'' is done in parallel, meaning cuts are not added from stage $i+1$ to $i$ before solving stage $i$ at each iteration. This variant has an added benefit for LPs because the backwards pass has no computational cost since dual information is available from the forward pass. This approach has been taken in literature (see \cite{kumar2021dual} or \cite{avila2022parallel}), and both implementations of the backwards pass are supported in {\tt PlasmoBenders.jl} beginning in version 0.2.0.

\begin{figure*}[!htp]
    \centering
    \includegraphics[width = \textwidth]{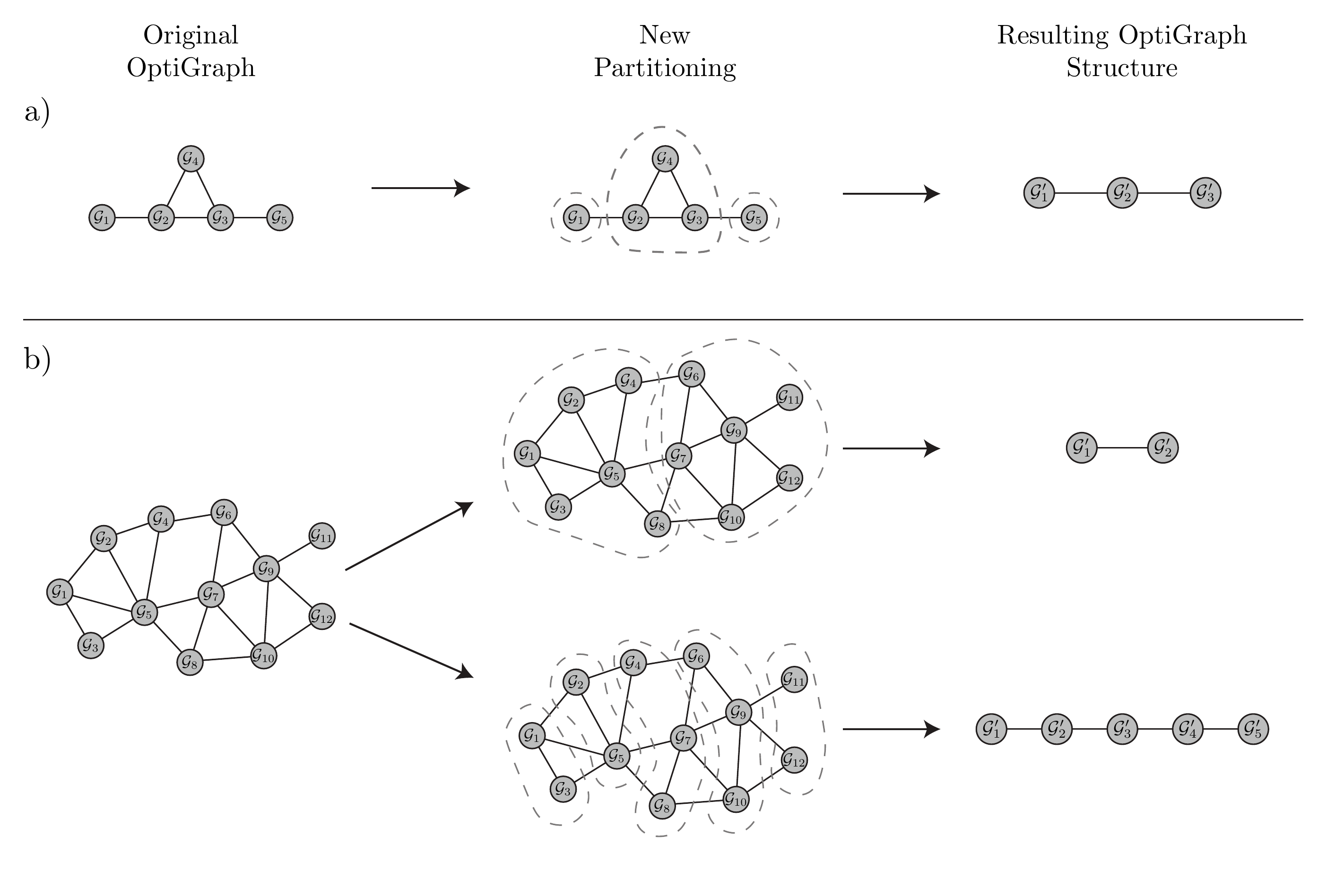}
        \vspace{-0.1in}
    \caption{Examples of cycle removal by aggregation. a) Shows a simple structure with a cycle where the cycle is destroyed by placing the cycle on its own partition. b) Shows a more complex structure that can be partitioned into subgraphs for the traditional BD approach or partitioning into a longer linear graph.}
    \label{fig:NBD_cycles}
\end{figure*}

The requirement that there be no cycles among subgraphs for gBD may seem restrictive. However, the flexibility of the OptiGraph abstraction can allow the user to obtain the required structure from any graph via aggregation/partitioning  procedures. For example, Figure \ref{fig:NBD_cycles}a shows a subgraph structure containing a cycle (among $\mg_2$, $\mg_3$, and $\mg_4$). However, if the graph is partitioned with the cycle being placed on a single subgraph, it results in a subgraph structure without cycles. This can be thought of as aggregating the cycle into a single node or subgraph. More complex structures can similarly be partitioned into the required structure, as shown in Figure \ref{fig:NBD_cycles}b, where we give two examples. First, we can arbitrarily partition any set of subgraphs into two separate subgraphs, forming the traditional BD case. In fact, any optimization problem with linking constraints can be partitioned into this traditional BD structure, since we can arbitrarily choose to place variables into one of two subgraphs (in this sense, {\em the gBD algorithm is broadly applicable to any graph structure}). Second, we can set partitions in such a way that we create an acyclic structure. In this case, we set five partitions and create a linear graph structure (such as that seen in dynamic programming). If we consider the original graph to be a spatial structure, we can see that we can recover a linear graph structure that is commonly seen in temporal problems (this is the structure required by dual dynamic programming), and to which we can apply gBD.

\vspace{0.1in}

\begin{figure*}[!htp]
    \centering
    \includegraphics[width = \textwidth]{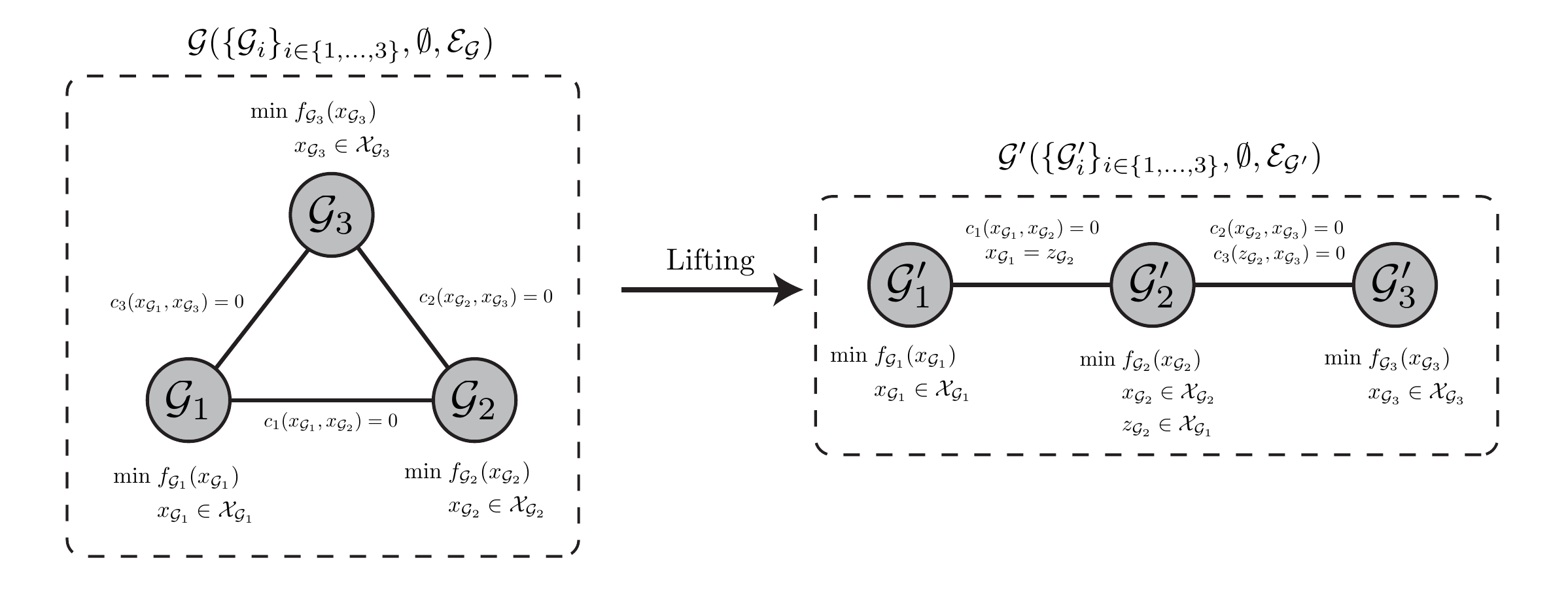}
        \vspace{-0.1in}
    \caption{Visualization showing how lifting of an optimization problem changes the graph structure. The left graph, $\mg$, contains cycles among its subgraphs. However, dummy variables $z_{\mg_2}$ can be added to subgraph $\mg_2$ and set equal to the variables $x_{\mg_1}$ and then the constraint $c_3(x_{\mg_1}, x_{\mg_3})$ can be replaced by $c_3(z_{\mg_2}, x_{\mg_3})$ which eliminates the connection between $\mg_1$ and $\mg_3$ and results in the graph $\mg'$ which has no cycles. }
    \label{fig:lifting}
\end{figure*}

While partitioning and aggregation functionalities allow for removing cycles, structures can also be specified directly by the modeler. Figure \ref{fig:lifting} shows an example of a graph with three subgraphs that form a cycle. The cycle can be removed by creating dummy variables on one of the subgraphs such that a linking constraint can be removed (this is a lifting procedure). In this case, a variable $z_{\mg_2}$ is added to $\mg_2$ that is constrained to be equal to the variables on $\mg_1$ (i.e., $z_{\mg_2} = x_{\mg_1}$). Then, the constraint $c_3(x_{\mg_1}, x_{\mg_3}) = 0$ can be placed instead between $\mg_2$ and $\mg_3$ by the constraint $c_3(z_{\mg_2}, x_{\mg_3})$ which forms a linear graph. This changes the structure of the problem but results in an optimization problem that will have the same solution as before the lifting procedure. This same procedure could be applied to more complex graphs with many cycles.

\begin{figure*}[!t]
    \centering
    \includegraphics[width = 0.9\textwidth]{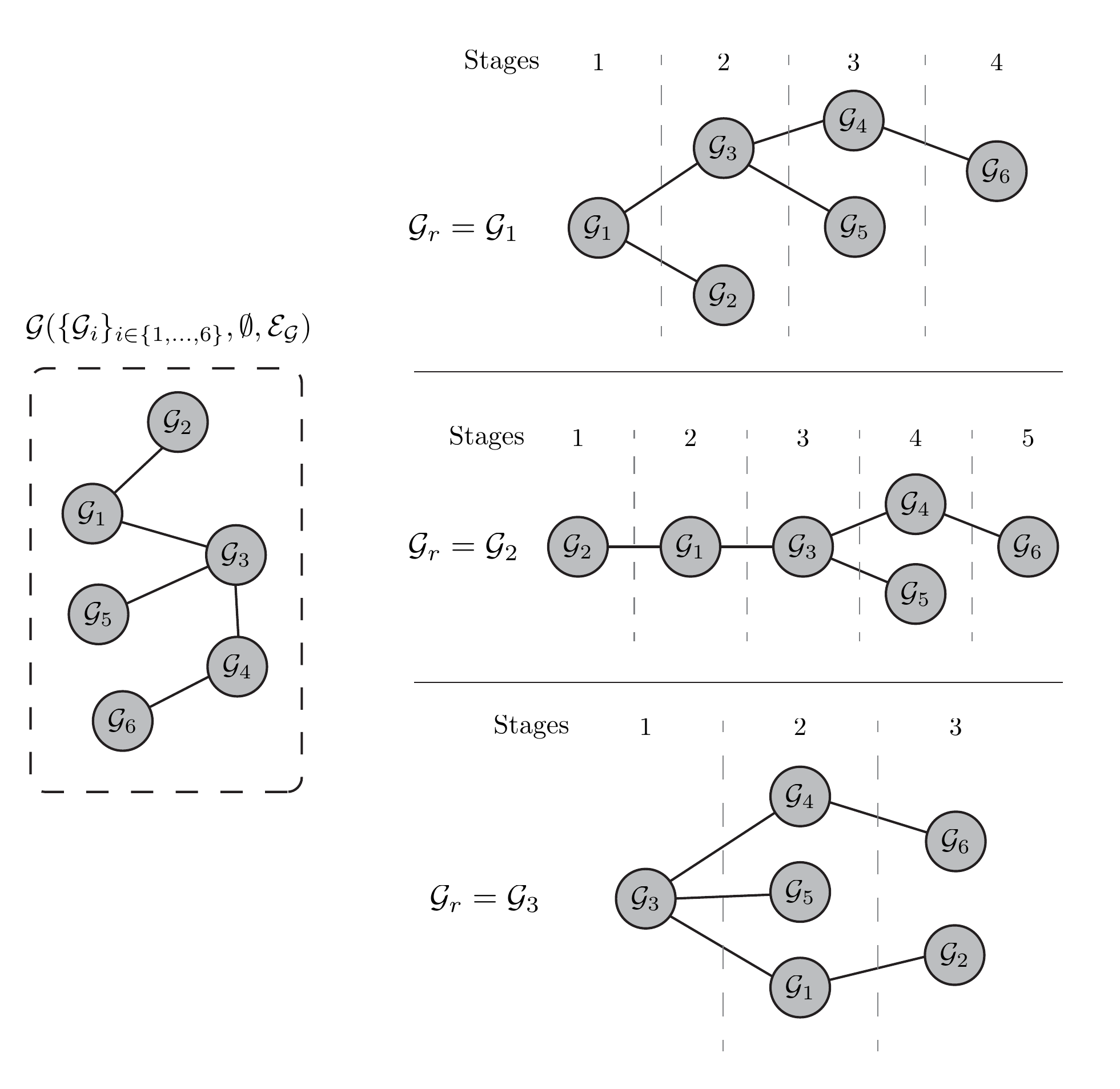}
        \vspace{-0.2in}
    \caption{Visualization of how the choice of root subgraph influences how the gBD algorithm is applied. The graph, $\mg$, on the left contains six subgraphs but has the required structure for gBD. On the right are the resulting stages for setting $\mg_1$, $\mg_2$, and $\mg_3$ as the root subgraph. Each results in a different number of stages and could result in different cuts being added to the subproblems. For instance, when using multi-cuts, using $\mg_1$ or $\mg_2$ as the root subgraph results in two sets of cuts to $\mg_3$ at each iteration, but using $\mg_3$ as the root subgraph results in three sets of cuts being added to $\mg_3$ at each iteration.}
    \label{fig:NBD_root_graphs}
\end{figure*}

The gBD algorithm (Algorithm \ref{alg:graph_NBD}) also provides flexibility in how this is applied based on the choice of the root subgraph. Any subgraph in the graph can be set as the root subgraph, so long as the subgraphs form the required structure noted above. The choice of root subgraph could have algorithmic and convergence implications. For instance, Figure \ref{fig:NBD_root_graphs} shows the stages for three different root subgraphs for a given OptiGraph. Each decision results in a different number of stages. Further, the number of connections to the next stage can also be different and could result in different cuts (or number of cuts). With $\mg_r = \mg_2$, the subgraph $\mg_3$ has two child subgraphs, but $\mg_r = \mg_3$ results in $\mg_3$ having three child subgraphs. As a further example, linear graphs can arise from temporal problems and can be solved using a nested BD approach (i.e., dual dynamic programming) with the root problem being the subgraph with the first time point. However, we could instead apply this same algorithm starting at the middle time point, resulting in half as many stages. With the flexibility provided by the partitioning and aggregation capabilities of the OptiGraph, many different root subgraphs could be explored for the same problem, such as exploring different partitionings of integer and continuous variables into subgraphs or choosing a root subgraph with a larger number of integer variables. This ability to define any subgraph as the root subgraph in gBD is an important and flexible ability of the algorithm.

\vspace{0.1in}

Cutting planes are added to \eqref{eq:NBD_root_cuts}, \eqref{eq:NBD_subroot_cuts}, \eqref{eq:GNBD_rootcuts}, or \eqref{eq:GNBD_subroot_cuts} in part based on problem structure. During the gBD algorithm, the graph is spanned in the forward direction, starting at the root subgraph and moving across stages. For instance, in Figure \ref{fig:NBD_graph}, the forward pass will start in Stage one and sequentially move to the next stages. After the forward pass, cuts can be added in a backward pass. Often, this will include starting at the final stage and spanning the graph backwards. This is often the case when subgraphs contain integer variables, since the forward pass requires solving an MILP, the backward pass can require solving an LP relaxation, depending on the types of cuts chosen. 

The cuts added to \eqref{eq:NBD_root_cuts}, \eqref{eq:NBD_subroot_cuts}, \eqref{eq:GNBD_rootcuts}, or \eqref{eq:GNBD_subroot_cuts} can take several forms \cite{lara2018deterministic,zou2019stochastic}. Below, we will detail Benders cuts, Lagrangian cuts, and strengthened Benders cuts based on the work of \cite{kumar2021dual,lara2018deterministic,zou2019stochastic}. We note that \cite{zou2019stochastic} applied their work to the case of binary state variables (i.e., binary complicating variables), but \cite{lara2018deterministic} extended their work to allow for both binary and continuous variables, which comes at the expense of losing the finite convergence property due to the potential of a duality gap. For this reason, we did not specify any variables being continuous or integer in the above algorithm definitions. The cuts we discuss below will have the general form:
\begin{align}\label{eq:general_cuts}
    \Psi_w^i(\Phi^i_w, \pi^i_w, \bar{x}_{p_w}) := \Phi_w^i + \pi_w^{i \top} (x_{p_w} - \bar{x}_{p_w}^i)
\end{align}
\noindent where $\Psi_w^i(\Phi_w^i, \mu_w^i, \bar{x}_{p_w})$ is the lower bound on the subproblem $w$, $\Phi_w^i$ and $\pi_w^i$ are information coming the child subproblems, and $x_{p_w}$ are the complicating variables that come from the parent subproblem, $p_w$ (with $\bar{x}_{p_w}$ being the fixed feasible solution; see \eqref{eq:NBD_subroot_links}, \eqref{eq:NBD_final_sub_links}, \eqref{eq:GNBD_sub_links}, and \eqref{eq:GNBD_final_links}). Here, the set of aggregated cuts at iteration $k$ is given by $\theta_{p_w} \ge \sum_{w \in \mc({p_w})} \Psi^i(\bar{\Phi}^i_w, \pi_w^i, \bar{x}_{p_w}), \; i = 0, 1, ..., k-1$ and the multi-cut is given by $\theta_{p_w, w} \ge \Psi^i(\bar{\Phi}^i_w, \pi_w^i, \bar{x}_{p_w}), \; w \in \mc(p_w), \; i= 0, 1, ..., k-1$. Each algorithm given above allows the user to set initial cut data (i.e., $i = 0$), where additional cuts are added after each iteration. 

Benders cuts come from LP duality properties \cite{kumar2021dual}. If \eqref{eq:NBD_subroot} or \eqref{eq:NBD_final_sub} are LPs, the Benders cut for decision $x_{p_w}$ at iteration $i + 1$ for the $w$th subproblem can be given by $\Psi^i(\bar{\phi}^i_w, \bar{\lambda}_w^i, \bar{x}_{p_w})$. This cut comes from the fact that \eqref{eq:NBD_subroot} and \eqref{eq:NBD_final_sub} can be written in dual form with the objective $\overline{\phi}^{i + 1}_w(\bar{x}^{i+1}_{p_w}) = \lambda_w^{i+1 \top} x_{p_w} + \nu^{i + 1}$, where $\nu^{i + 1}$ denotes the cost component of constraints \eqref{eq:NBD_subroot_cons}, \eqref{eq:NBD_link_constraint1}, \eqref{eq:NBD_final_sub_cons}, and \eqref{eq:NBD_link_constraint2}. The cost to go at iteration $i + 1$ can then be approximated by using information from the previous iteration, where we replace $\lambda_w^{i + 1}$ with $\bar{\lambda}_w^{i}$ and replace $\nu^{i+1}$ with $\bar{\phi}_w^{i}(\bar{x}^{i}_{p_w}) - \bar{\lambda}_w^{i \top } \bar{x}_{p_w}^{i}$, where $\bar{\phi}_w^{i}$ is the objective value of \eqref{eq:NBD_subroot} or \eqref{eq:NBD_final_sub} at iteration $i$. When \eqref{eq:NBD_subroot} or \eqref{eq:NBD_final_sub} are a MILP, rather than using $\bar{\phi}^i_w$, the Benders cut uses the solution of their LP relaxation. The above cuts are also easily extended to the graph-based case.

Lagrangian cuts come from a relaxation of the constraint \eqref{eq:NBD_subroot_links} or \eqref{eq:NBD_final_sub_links}. This relaxation for iteration $i$ is given by:
\begin{subequations}\label{eq:lagrange_relax}
    \begin{align}
        \phi^{LR,i}_w(\bar{x}_{p_w}^i, \mu_w) := \min &\; c_w^{\top} y_w - \mu^{\top}_w(z_w - \bar{x}_{p_w}^i)\\
        \textrm{s.t.} &\; B_w y_w \le b_w \\
        &\; C_w z_w + D_w y_w \le q_w
    \end{align}
\end{subequations}
\noindent for some Lagrange multiplier, $\mu_w$. Optimal Lagrange multipliers $\bar{\mu}_w^i$ are found by solving 
\begin{align}\label{eq:lagrange_multiplier}
    \Phi_w^{LR,i}(\bar{x}_{p_w}^i) = \max_{\mu_w} \phi^{LR,i}_w(\bar{x}_{p_w}^i, \mu_w).
\end{align}
Once the optimal lagrange multipliers, $\bar{\mu}^i_w$, have been determined, the lagrange cut takes the form $\Psi_w^i(\Phi^{LR, i}_w(\bar{x}_{p_w}^i), \bar{\mu}^i_w, \bar{x}_{p_w})$.

Benders cuts can be quick to compute but tend to be weak, while Lagrangian cuts can be expensive to compute because the solution of \eqref{eq:lagrange_multiplier} can be expensive. Strengthened Benders cuts were introduced by \cite{zou2019stochastic} and seek to mitigate performance issues. Rather than solving \eqref{eq:lagrange_multiplier}, they solve $\phi^{LR,i}_w(\bar{x}_{p_w}^i, \bar{\lambda}_w^i)$, where $\bar{\lambda}_w^i$ comes from \eqref{eq:NBD_subroot_links} or \eqref{eq:NBD_final_sub_links} at iteration $i$ (or the LP relaxation if it is a MILP). Thus, the strengthened cut has the form $\Psi^i_w (\bar{\phi}^{LR,i}_w (\bar{x}_{p_w}^i, \bar{\lambda}_w^i), \bar{\lambda}_w^i, \bar{x}_{p_w}^i)$. 

\subsection{Software Implementation of gBD -- {\tt PlasmoBenders.jl}}

To highlight how the graph-based abstraction provides a framework for implementing decomposition approaches, we now introduce the Julia package {\tt PlasmoBenders.jl} which implements gBD as outlined above. Using the OptiGraph abstraction, we can generalize the decomposition scheme which enables the user to build a problem as a graph and then pass it to the generalized, graph-based solver. Additional solution options and algorithmic approaches (e.g., different types of cuts, regularization schemes) can likewise be generalized to the graph-based approach and implemented in the solver as we will show below. In this case, we use {\tt PlasmoBenders.jl} to highlight these capabilities, but many of the principles described in this manuscript could be applied to other decomposition approaches. 

The primary object in {\tt PlasmoBenders.jl} is the {\tt BendersAlgorithm} object, which is motivated by the gBD algorithm (Algorithm \ref{alg:graph_NBD}). This algorithm is applied to a graph containing subgraphs with edges between subgraphs, and its implementation is dependent on setting a root subgraph. Similarly, after a user has defined an OptiGraph object in {\tt Plasmo.jl}, they can create a {\tt BendersAlgorithm} object by passing that OptiGraph object, $\mg$, as well as a root subgraph, $\mg_r \in \msg(\mg)$. After setting the graph and root subgraph, the {\tt BendersAlgorithm} constructor will make the changes to the subgraphs noted in Algorithm \ref{alg:graph_NBD} (e.g., adding cost-to-go variables or enforcing the linking constraints). The {\tt BendersAlgorithm} object also includes additional data structures for implementing gBD (e.g., best upper bounds, lower bounds, solver tolerance). 

{\tt PlasmoBenders.jl} gives the user access to several additional options and functionalities. More details on the package operation can be found in the documentation at \url{https://github.com/plasmo-dev/PlasmoAlgorithms.jl/tree/main/lib/PlasmoBenders} and in the Supporting Information, but we will highlight some of these here. Note that some of these topics (like regularization) are able to be generalized to the graph structure, simplifying their implementation. Available options include:
\begin{itemize}
    \item Solver information such as $K_{max}$ or $\epsilon_{tol}$
    \item Whether to use strengthened Benders cuts rather than regular Benders cuts (if the problem is mixed-integer)
    \item Whether to use aggregated or multi-cuts
    \item Whether to parallelize certain computations, such as solutions of Benders subproblems when $S=2$ and  $|W_2|>1$
    \item Whether to add slack variables to constraints \eqref{eq:NBD_link_constraint1} or \eqref{eq:NBD_link_constraint2} to help ensure recourse. 
    \item Whether to warm-start the next iteration with the current best solution
    \item Whether to use a regularization scheme to get the variables $\bar{x}^i_g$ that are fixed in each subproblem at iteration $i$ (see \cite{pecci2024regularized} and Supporting Information).
    \item Beginning in version 0.1.0, whether to use feasibility cuts if subproblems are infeasible (only supported for two stage problems)
    \item Beginning in version 0.2.0, whether to use a sequential backwards pass, where the stages are resolved in reverse order after cuts are passed to the previous stages.
\end{itemize}

As {\tt PlasmoBenders.jl} is used in some of the case studies below, we give an example in Code Snippet \ref{code:plasmobenders_example} to show how it can be applied to an OptiGraph defined in {\tt Plasmo.jl}. First, we define a graph with three nodes  and two edges (Lines \ref{line:PB_graph_start} - \ref{line:PB_graph_end}), where one edge is between node one and node two and one edge is between node two and node three (i.e., this is a linear graph). The graph is then partitioned so that each node is on its own subgraph (Lines \ref{line:PB_partition_start} - \ref{line:PB_partition_end}). We then define a solver for the subgraphs (Line \ref{line:PB_subproblem_solver}) which we will pass to the constructor to set on each of the subgraphs internally. Next, we can build the {\tt BendersAlgorithm} object by calling the function {\tt BendersAlgorithm} (Lines \ref{line:PBO_start} - \ref{line:PBO_end}). This function takes the overall OptiGraph object {\tt graph} as well as the root subgraph, which we use here as the first subgraph (containing node one). We can also set other information, including the maximum number of iterations, the solver tolerance, to use strengthened cuts, not to use regularization, not to add slack variables, not to use multicuts, and what to use for the subproblem solver. We can then call {\tt run\_algorithm!} on the {\tt BendersAlgorithm} object (line \ref{line:PB_optimize}) which applies gBD to solve the problem. Solution information can be queried from the {\tt BendersAlgorithm} object using {\tt PlasmoBenders.jl} functions and extended {\tt JuMP.jl} functions (Lines \ref{line:PB_query1} - \ref{line:PB_query2}). 

\begin{figure}[!htp]
    \begin{minipage}[t]{0.9\linewidth}
        \begin{scriptsize}
        \lstset{language=Julia, breaklines = true}
        \begin{lstlisting}[label = code:plasmobenders_example, caption = {Code for applying graph-based Benders and Nested Benders decompositions to an OptiGraph defined in {\tt Plasmo.jl}. The package {\tt PlasmoBenders.jl} is able to apply these decomposition schemes based on the user-defined graph and root subgraph.     \vspace{0.1in} }] 
using Plasmo, PlasmoBenders, HiGHS

# Define a graph with three nodes and two edges
graph = OptiGraph() |\label{line:PB_graph_start}|
@optinode(graph, nodes[1:3])

# Set local variable, constraint, and objective information
for node in nodes
    @variable(node, x, Bin)
    @variable(node, y >= 0)
    @constraint(node, x + y >= 1.3)
    @objective(node, Min, x + 2 * y)
end

# Add linking constriants
for i in 1:2
    @linkconstraint(graph, nodes[i][:x] + nodes[i + 1][:y] >= i)
end |\label{line:PB_graph_end}|

# Create subgraphs with one node on each subgraph
part_vector = [1, 2, 3] |\label{line:PB_partition_start}|
partition = Partition(graph, part_vector)
apply_partition!(graph, partition) |\label{line:PB_partition_end}|

# Define a solver to use for the subproblems
solver = optimizer_with_attributes(HiGHS.Optimizer, "output_flag" => false) |\label{line:PB_subproblem_solver}|

# Construct the BendersAlgorithm object
benders_alg = BendersAlgorithm( |\label{line:PBO_start}|
    graph, # Set the overall graph
    local_subgraphs(graph)[2]; # Set the root subgraph
    max_iters = 20, # Maximum iterations
    tol = 1e-7, # Solver tolerance
    parallelize_benders = true, # Whether to parallelize the second stage solutions
    strengthened = true, # Whether to use strengthened cuts
    regularize = false, # Whether to use a regularization scheme
    add_slacks = false, # Whether to add slack variables for recourse
    multicut = true, # Whether to use multi-cuts rather than aggregated
    solver = solver # Set the solver on the subgraphs
) |\label{line:PBO_end}|

# Solve the problem using gBD
run_algorithm!(benders_alg) |\label{line:PB_optimize}|

# Query Solutions
println("Objective value is ", JuMP.objective_value(benders_alg)) |\label{line:PB_query1}|
println("Relative gap is ", PlasmoBenders.relative_gap(benders_alg))
println("Value of x on node 1 is ", [JuMP.value(benders_alg, nodes[1][:x])]) |\label{line:PB_query2}|
\end{lstlisting}
\end{scriptsize}
\end{minipage}
\end{figure}

The generalization of the BD algorithm to graphs and its implementation in \texttt{\mbox{\hspace{30pt}} \allowbreak PlasmoBenders.jl} allows for a simple user-interface for applying what otherwise might be complex algorithms. In Code Snippet \ref{code:plasmobenders_example}, it requires less than fifteen lines of code and only two function calls to solve a problem using gBD. Because the algorithm is generalized to a graph, the problem can be modeled as a {\tt Plasmo.jl} OptiGraph, and the structure of that OptiGraph, along with the definition of the root subgraph, determines the complicating variables and subproblem sizes. It also gives the user flexibility in how they will apply these algorithms such as allowing them to use strengthened cuts or regularization. Further, the graph above is a linear graph with three subgraphs. In Code Snippet \ref{code:plasmobenders_example}, we solve using gBD and set the root subgraph to be {\tt local\_subgraphs(graph)[2]} (the second subgraph) which forms two stages, with two subgraphs in the second stage. We could instead use almost the same code but set the root subgraph to be {\tt local\_subgraphs(graph)[1]} (the first subgraph) which would have three stages. If we use the first subgraph as the root subgraph, the {\tt multicut} option would have no influence because there is only one subpgraph in each stage. The {\tt parallelize\_benders} option also would throw an error since this parallelization is only implemented for two-stage problems. Moving between these algorithm options is simple to do in practice because it can be done by changing a few arguments in the {\tt BendersAlgorithm} constructor function.

\subsection{Other Decomposition Approaches}

While this section has focused on gBD as a solution approach, other decomposition approaches for graphs are possible and these are not restricted to LP/MILP problems. Other approaches have been applied to graph-based optimization problems, either on the modeling level (using graphs to form subproblems) or on the linear algebra level (exploiting the structure formed by the graph for different computations). Jalving and co-workers applied an overlapping Schwarz decomposition scheme for LP and NLP graph structured problems \cite{jalving2022graph,shin2020decentralized,shin2021}. A Schur decomposition has also been applied to graphs in two-stage stochastic NLPs \cite{cole2022julia}, and an ADMM approach was applied to the graph structure of a NLP power system problem \cite{shin2019hierarchical}. In addition, a Hansknecht et al. \cite{hansknecht2024framework} showed how the linear algebra resulting from tree graph structures can be exploited or decomposed for applications in non-linear programming or sequential quadratic programming. These or other decomposition schemes could be generalized to graph structure as well, and there are instances of approaches such as ADMM being applied to problems that we believe could be readily modeled by the OptiGraph abstraction; for example, multiple studies introduce different distributed control strategies where ADMM is used as a solution approach for coordinating distributed subsystems, forming a hierarchical structure \cite{behrunani2024distributed,braun2018hierarchical,wang2021admm}. 

\begin{figure}
    \centering
    \includegraphics[width=0.9\linewidth]{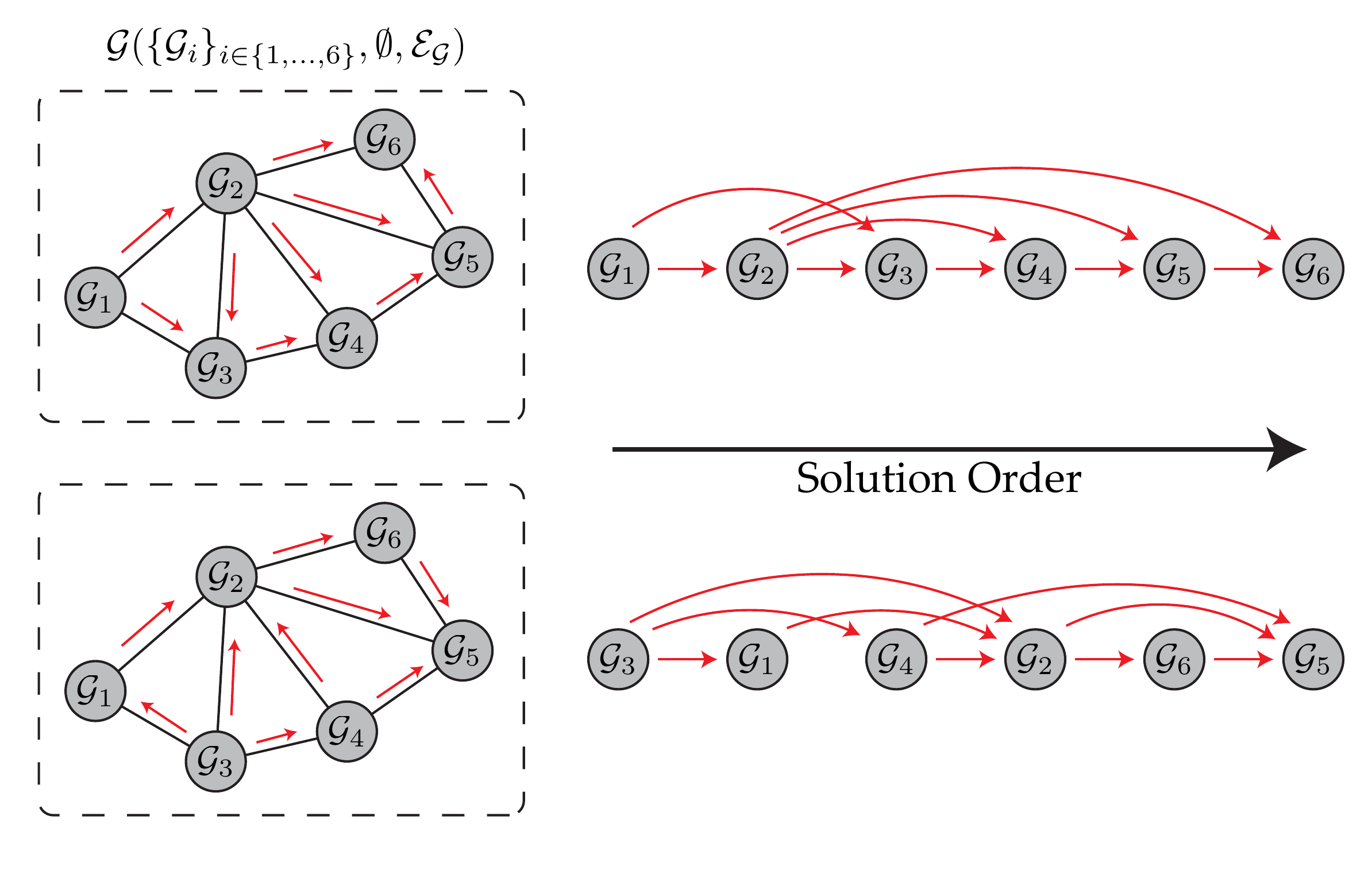}
    \caption{Visualization showing how an OptiGraph can be solved sequentially to obtain a feasible solution, with two different solution orders shown. Subgraphs are solved in a specific order and solutions are passed to the next subgraph in order. Red arrows represent the direction in which solutions are shared. Linking constraints on the edges are enforced on the ``downstream'' subgraph. Solution order can be arbitrary, and a tree structure is not necessarily required. }
    \label{fig:approximation_schemes}
\end{figure}

Hierarchical structures are also often leveraged for obtaining {\em approximate solutions}, such as partitioning problems into subproblems and solving those subproblems sequentially; this is often done when the problem is computationally intractable or when it is impractical to solve the full problem to optimality. For instance, doing a single iteration of gBD results in a feasible solution (assuming complete recourse) from the forward pass through the graph and could be used as a solution to the problem. This sequential solution approach is known as receding-horizon decomposition and is widely used in dynamic programming and model predictive control (for linear or tree graphs). Sequential decomposition approaches are used in many applications, including to hierarchical optimization problems (where decisions flow from one layer to the next), multi-scale problems (where a coarse problem is solved to create targets in a higher resolution problem), and to clustered problems, where clusters are iteratively solved and solutions shared between clusters; further, this approach is not restricted to LPs, MILPs, NLPs, or MINLPs. A visualization of how a sequential solution approach can be applied to a graph is shown in Figure \ref{fig:approximation_schemes}. Note that the choice of solution order, like the choice of root subgraph in gBD, is a choice made by the user. In addition, the sequential solution approach is not necessarily restricted to tree structures. In this case, each subgraph is solved in a set order and solutions (represented by the red arrows) are shared to the next subproblem. The linking constraints contained on the edges would be enforced on the downstream subgraph. Assuming recourse, this sequential approach yields an upper bound on the optimal solution. In addition, a lower bound can be obtained by solving each subgraph in parallel, ignoring the linking constraints between subgraphs.

There are many reported applications of sequential approximation approaches. Power systems studies frequently solve a series of separate optimization problems, such as first solving an MILP unit commitment problem, then solving a LP economic dispatch problem and/or LP or NLP system dynamics problems \cite{conejo2018power,lara2024powersimulations,wood2013power}. This sequential solution approach is often required to ensure tractability or because certain information, such as future demand, is not yet realized. Supply chain problems often include multiple echelons or decision layers \cite{liu2022optimization,yue2014biomass} which can be solved sequentially, such as is done by \cite{sharma2011decision} for biorefineries, where multiple MILPs for flexibility and structural evaluation are solved sequentially. Supply chains often result in complex structures in addition to hierarchical layers, and in \cite{guerra2019integrated}, they solve coupled shale gas and water supply chains by clustering the problem spatially (rather than by echelons) and sequentially solving and sharing solutions between clusters. Similarly, problems in health care (often routing and scheduling problems) have been addressed using sequential strategies, such as by solving two problems in series (sharing solutions of the first problem with the second \cite{rais2011operations}) or by solving initial clustered problems and using the solutions to inform a secondary problem \cite{bertels2006hybrid,fikar2017home}. In heating, ventilation, and air conditioning systems, \cite{long2016hierarchical} proposed a hierarchical control approach for coordinating temperature control of multiple rooms, where a coarse representation of the full system is solved, and the solution is then passed to lower-level controllers which operate independently. These different sequential approaches can be implemented within the OptiGraph abstraction by using subgraphs to capture the different subproblems. In addition, the freedom provided by the OptiGraph structuring capabilities (e.g., partitioning and aggregation) allows a user to explore additional representations of the same problem, which could lead to other applications of these sequential solution approaches. For instance, exploring different partitioning approaches and solving the resulting subgraphs sequentially could potentially yield multiple upper bounds on the optimal solution (one for each partitioning), while solving these same subgraphs in parallel (ignoring linking constraints between them) would provide different lower bounds.

\section{Case Studies}\label{sec:case_studies}
In this section, we use challenging case studies that show how the graph approach can be used for modeling hierarchical structures and obtain approximate solutions (case study I) and how gBD can be applied to solve hierarchical  problems (case studies II and III). Scripts for reproducing all case studies, as well as a script for creating some of the visuals in case study 3, are available at \url{https://github.com/zavalab/PlasmoBendersExamples}.

\subsection{Case Study 1:Hierarchical Electricity Markets}

\begin{figure}[!ht]
    \centering
    \vspace{-0.1in}
    \includegraphics[width=1.1\textwidth]{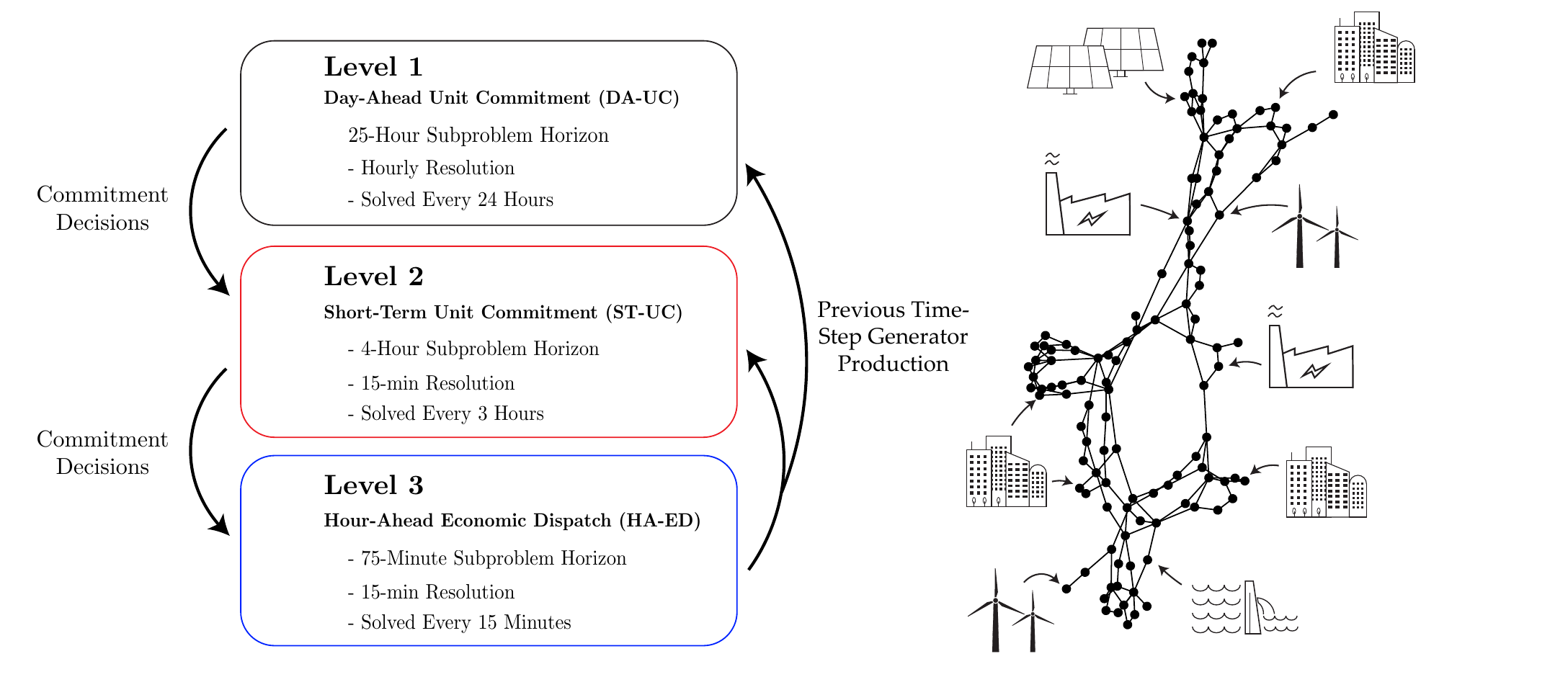}
    \caption{The tri-level hierarchical architecture of Atakan and co-workers \cite{atakan2022towards} of an electricity market.} 
    \label{fig:HP_framework}    
\end{figure}

This problem is based on the work of Atakan and co-workers \cite{atakan2022towards} who presented a tri-level hierarchy of market operations (Figure \ref{fig:HP_framework}). This included three hierarchical layers that could each have different time resolutions and time horizons: a day-ahead unit commitment layer (DA-UC), a short-term unit commitment layer (ST-UC), and an hour-ahead economic dispatch layer (HA-ED). Atakan and co-workers \cite{atakan2022towards} applied this framework to a 118-bus system \cite{pena2017extended}, which we will likewise follow below. They treated this as a stochastic problem, but we will only use the deterministic form of the problem for simplicity. We will give a general overview of this problem and show how we model and solve it using OptiGraphs. Additional mathematical details can be found in the Supporting Information.

An overview of the three hierarchical layers are described below and in Table \ref{tab:problem_sizes}. Each hierarchical layer consists of different numbers of subproblems which each have different horizons and degrees of overlap between subproblems. The DA-UC subproblem is an MILP that spans a 25-hour period, is solved once per day, and runs with a coarse resolution. This layer commits an initial subset of generators while considering operation of these generators to meet a demand profile (e.g., it considers ramp up/down constraints and minimum up/down time constraints). The ST-UC subproblems are MILPs that span 4-hour periods at a 15 minute resolution and which are solved every 3 hours; they use the commitment decisions of the DA-UC layer and also commit an additional subset of generators. The ST-UC layer considers operation of {\it both} the DA-UC generators and ST-UC generators, but maintains the commitment decisions of the DA-UC generators from the previous layer. The commitment decisions of the DA-UC and ST-UC layers are passed to the HA-ED layer, which include LP subproblems over a 75-minute time span and which are solved every 15 minutes. Each layer also includes power provided by renewable generators as well as demand data. The demand data can vary between layers (e.g., the realized demand of the HA-ED layer may differ from the projected demand of the DA-UC layer; see Supporting Information). Operation of the power systems in each layer was also subject to transmission constraints using DC approximations of the power flow. For this case study, we used the 118-bus system \cite{pena2017extended} used by Atakan and co-workers \cite{atakan2022towards}. 

\begin{table}
\caption{Details for the DA-UC, ST-UC, and HA-ED layers of the tri-level hierarchical optimization problem of Atakan et al. \cite{atakan2022towards}.}
\label{tab:problem_sizes}
\begin{tabular}{llll}
\hline\noalign{\smallskip}
    & DA-UC & ST-UC & HA-ED\\
\noalign{\smallskip}\hline\noalign{\smallskip}
    Subproblem Horizon (hr) & 25 & 4 & 1.25 \\
    Discretization Size (min) & 60 & 15 & 15 \\
    Subproblems per Day & 1 & 8 & 96 \\
    Variables per Subproblem & 31,725 & 19,008 & 4,770 \\
    Binary Variables per Subproblem & 11,775 & 3,744 & 0 \\
    Constraints per Subproblem & 87,525 & 63,440 & 18,265  \\
    Ren. Generators  & 31 & 31 & 31 \\
    Conv. Generators & 157 & 235 & 235\\
    Conv. Gens. Committed & 157 & 78 & 0 \\
    Total Generators & 188 & 266 & 266  \\
\noalign{\smallskip}\hline
\end{tabular}
\end{table}

\begin{figure}
    \centering
    \includegraphics[width=\linewidth]{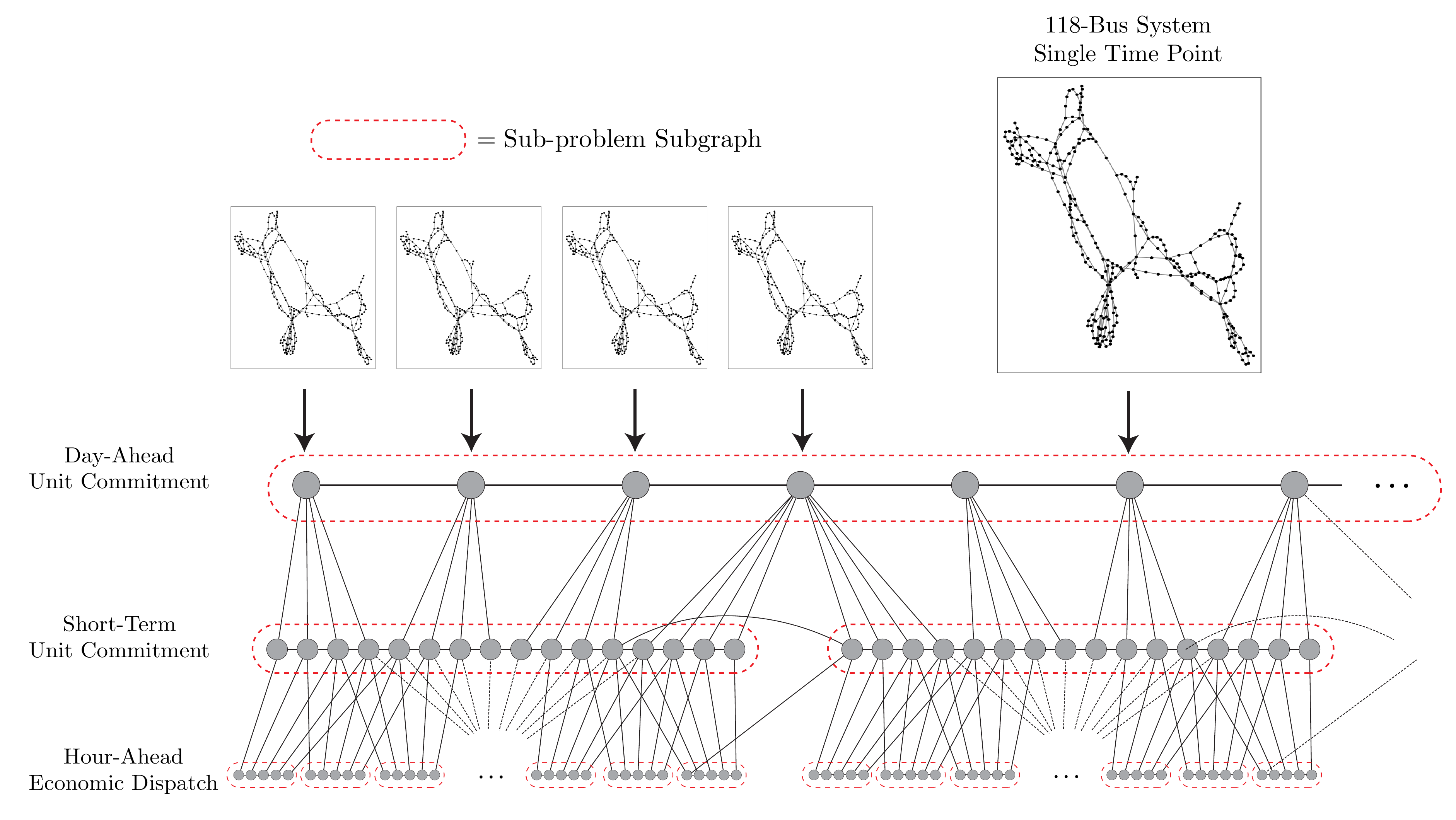}
    \caption{Visualization of the hierarchical graph structure that is formed by the tri-level markets problem presented by Atakan and co-workers \cite{atakan2022towards}. Here, subproblem subgraphs are formed for each subproblem in the DA-UC, ST-UC, and HA-ED layers and nested on an overall monolithic graph for each day of operation. Each ``node'' in this figure is a nested subgraph representing a single time point; this subgraph encodes the spatial structure of the 118-bus system as an OptiGraph, which buses and transmission lines represented by nodes. Edges are placed between the different subgraphs and across hierarhical layers and subproblems.}
    \label{fig:CS1_graph}
\end{figure}

This problem is constructed as a graph in {\tt Plasmo.jl}. An overall, monolithic graph, representing a single day, contains a subgraph for each subproblem (i.e., the overall graph contains 1 DA-UC subgraph, 8 ST-UC subgraphs, and 96 HA-ED subgraphs). Each subproblem subgraph contains nested sugraphs for any given time point. Each time point subgraph is made up of nodes for each bus in the system and nodes for each transmission line in the system. Each of these nodes contains variables corresponding to their objects (e.g., bus nodes contain variables for generators and line nodes contain variables for power flow). Edges are placed between the time point subgraphs in the same subproblem and between bus and line nodes (e.g., for up/down time constraints on generators or transmission constraints on lines). Additional edges are also placed between subproblem subgraphs which link solutions of different layers. A visualization of the resulting graph structure is shown in Figure \ref{fig:CS1_graph}. Here, the subproblem subgraphs are highlighted by red bubbles, and each time point is represented by a node (which encodes the 118-bus system OptiGraph). Edges are placed that link solutions of each layer. There are also additional edges between the DA-UC and HA-ED layers that are not shown which enforce the DA-UC commitment decisions in the HA-ED subgraphs. In addition, each of the subproblem subgraphs on this figure exists on an overall, monolithic graph. The system we consider here has 32 days of data, so there are 32 different monolithic graphs that will be created, with the solutions of the previous day problems being fixed in the subsequent day graph (note that we only keep 3 graphs in memory at a time to reduce overall memory requirements). Each monolithic graph contained 641,709 total variables (41,727 binary), highlighting the inherent high dimensionality of the problem. Further details on the graph representation can be found in the Supporting Information.

\begin{figure}
    \centering
    \includegraphics[width=0.9\linewidth]{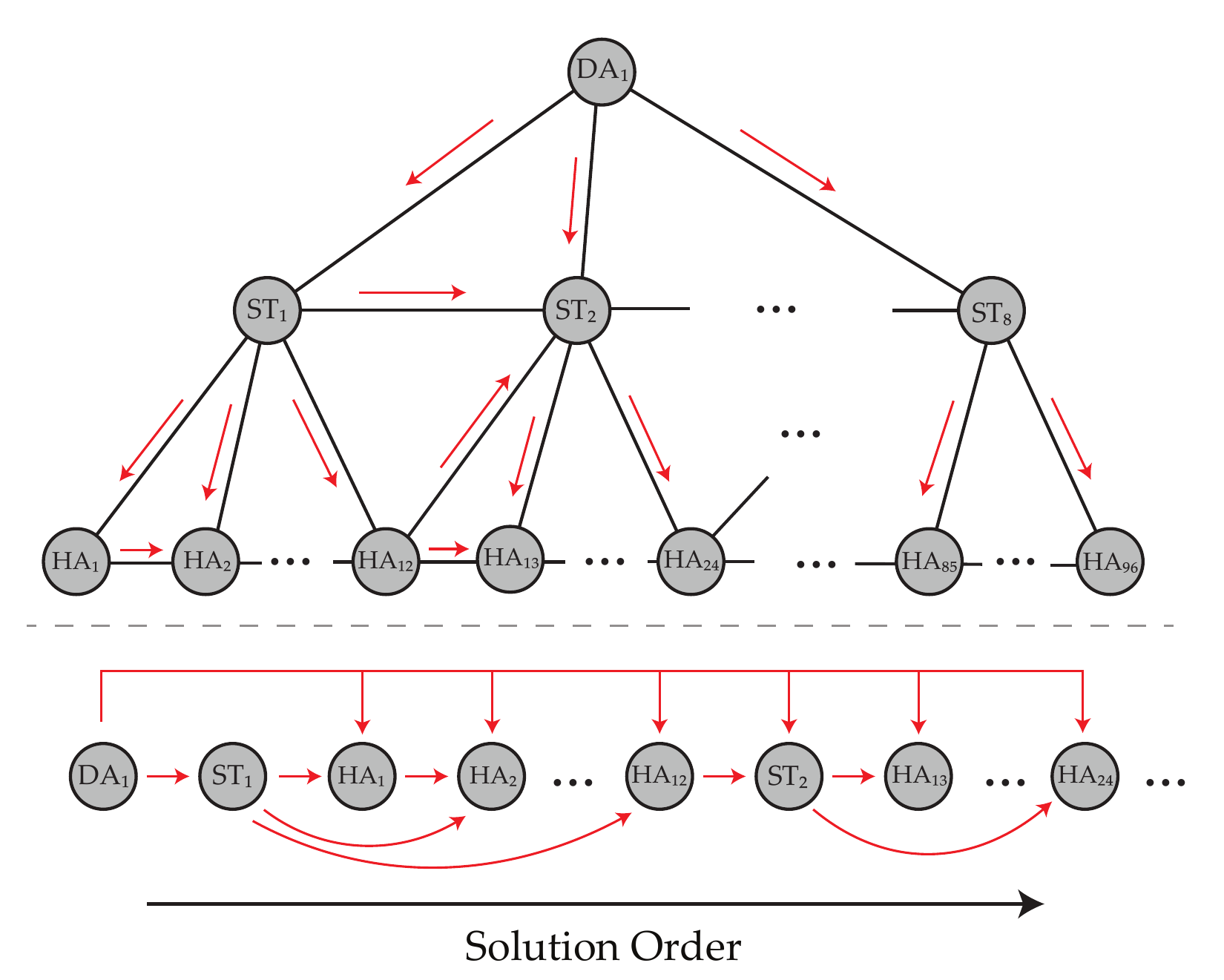}
    \caption{Visualization of the receding horizon decomposition approach for the hierarchical energy markets case study. Each subproblem subgraph of Figure \ref{fig:CS1_graph} are represented by a node, with DA$_i$, ST$_i$, and HA$_i$ representing the $i$ subproblem of the DA-UC, ST-UC, and HA-ED layers, respectively. Red arrows show the direction in which solutions are shared.}
    \label{fig:receding_horizon}
\end{figure}

We solve this problem using two different approaches. The first approach is to solve the overall ``monolithic'' graph as a single problem for each day. The second is a ``receding horizon'' approximation approach, where each subproblem subgraph is solved sequentially. In the receding horizon approach, we first solve the DA-UC subgraph, then the first ST-UC subgraph, then the 12 HA-ED subgraphs linked to the first ST-UC subgraph, then the second ST-UC subgraph, and so forth. This approach is represented in Figure \ref{fig:receding_horizon}, where each subgraph is represented by a node, and red arrows highlight information being shared to the next subproblem. These two approaches were repeated over all 32 days for which we have data, and both approaches used the same demand profiles in each layer for comparison. This highlights how the graph representation allows us to explore different ways of computing exact (when possible) or approximate solutions. We also note that the implementation of these approaches required slight modeling differences since, in the receding horizon approach, the linking constraints like those between layers had to be enforced on the subgraphs as well.

These approaches result in different overall solutions and solution times. We solved these problems with Gurobi with access to 20 threads on a 40-core Intel(R) Xeon(R) CPU E5-2698 v4 processor running @2.2 Ghz. The overall cost of economic dispatch (based on the first time point of each HA-ED subproblem; this is the realized performance) was \$ 199.6 million and  \$ 222.1 million for the monolithic and receding horizon approaches, respectively (a difference of 10.8\%). However, the monolithic approach took 48.7 hours to solve, compared to 7.2 hours for the receding horizon approximation. In addition, we used a tighter gap for the receding horizon DA-UC and ST-UC layers (0.5\% or less, with a maximum time limitation) compared to a 1\% gap on the monolithic problem (though using an even larger gap could reduce runtime). Further details on the solutions are included in the Supporting Information.

This case study highlights how the graphs facilitate hierarchical model building and enable different solution approaches. Often, the receding-horizon approach above is used because the monolithic problem is too large to solve in a reasonable time, because some information is not yet know (e.g., realized demand 12 hours from now), or a combination of both. The OptiGraph abstraction provides a framework in which to test and apply different approaches that operate within these types of constraints; for instance, a receding-horizon approach could be applied at an even higher resolution, where each time point subgraph is solved in series. Alternatively, problems could be partitioned in space (rather than their current time partitioning), or different numbers of hierarchical layers or different time scales could be applied, and the OptiGraph abstraction could help facilitate these analyses. 

\subsection{Case Study 2: Integration of Planning and Operations in Power Grid Expansion} \label{sec:cs2}

This problem highlights how gBD can be applied to a capacity expansion model (CEM) with both a planning level (long-term decisions) and an operations level (short-term decisions). CEMs are models of power systems that seek to minimize total system cost by considering investment and retirement (planning-level) decisions for power systems. CEMs also often model higher resolution operation, such as regular (e.g., hourly) operation of generators (unit commitment and dispatch decisions), energy storage, power transmission, and flexible demands. This allows the CEM to capture the impact of planning level decisions on the day-to-day operation of the system and the resultant approximation of likely total system cost, permitting co-optimization of both fixed and variable costs. The problem we consider is implemented using the open-source software package, {\tt GenX.jl} \cite{jenkins2017enhanced}. This implementation includes variables at the planning-level (including policy decisions and investment decisions) and variables at the operations-level, where the operations-level is split into multiple sub-periods, each constituting a series of linked operational steps (e.g. a \textit{n}-day long periods of hourly sequential operations).

A compact form of the CEM (adpated from \cite{jacobson2024computationally} and \cite{pecci2024regularized}) is given by
\begin{subequations}\label{eq:CEM}
\begin{align}
    \min &\; c_p^\top x_p + \sum_{w \in W_o} c_w^\top x_w \label{eq:CEM_obj}\\
    \textrm{s.t.} &\; A_w x_w + B_w x_p \le b_w, \quad w \in W_o \label{eq:CEM_link_cons}\\
    &\; \sum_{w \in W_o} Q_w x_w \le d \label{eq:CEM_policy_cons} \\
    &\; x_w \ge 0, \quad w \in W_o \\
    &\; x_p \in \mathcal{X}_p
\end{align}
\end{subequations}
\noindent where $W_o$ is the set of operations-level sub-periods, $x_p$ is the set of variables on the planning level, $x_w$ is the set of variables on the operations-level for the $w$th sub-period, $\mathcal{X}_p$ is a feasible region for the planning variables, and $c_p$, $c_w$, $A_w$, $B_w$, $b_w$, $Q_w$, and $d$ are problem data (vectors and matrices). The objective \eqref{eq:CEM_obj} defines total system cost, where $c_p$ represent fixed costs (fixed operation costs and investment costs) and $c_w$ represent variable costs (e.g., fuel prices, variable operation costs, penalties for non-served energy or constraint violations). The constraints \eqref{eq:CEM_link_cons} include operational constraints that do not link sub-periods, such as constraints on power generation resources, transmission and ramping limits, and start-up/shut-down constraints. The constraints \eqref{eq:CEM_policy_cons} include policy constraints, such as CO$_2$ emission limits or renewable portfolio standards which couple operational subproblems. 

We apply this formulation to an application in {\tt GenX.jl} example systems. The specific application includes 11 operational sub-periods, each 168 hours long with hourly resolution. The problem models three zones in the northeastern United States with transmission allowed between zones. Each zone includes thermal generators, variable renewable energy generators (e.g., wind, solar PV), battery storage, and demands. The problem also includes policy constraints on CO$_2$ production, capacity reserve requirements, energy share requirements, and minimum and maximum capacity requirements. The final problem was an LP and included 138,972 variables and 166,612 constraints. This problem is not especially large but serves to highlight how gBD can be applied. Further, this model can scale rapidly, such as by considering  additional operations sub-periods (e.g., weather years) or increasing the spatial resolution, and this same approach could be readily applied to a larger, scaled-up version of the problem where conventional monolithic solution methods struggle.

\begin{figure}
    \centering
    \includegraphics[width=0.9\linewidth]{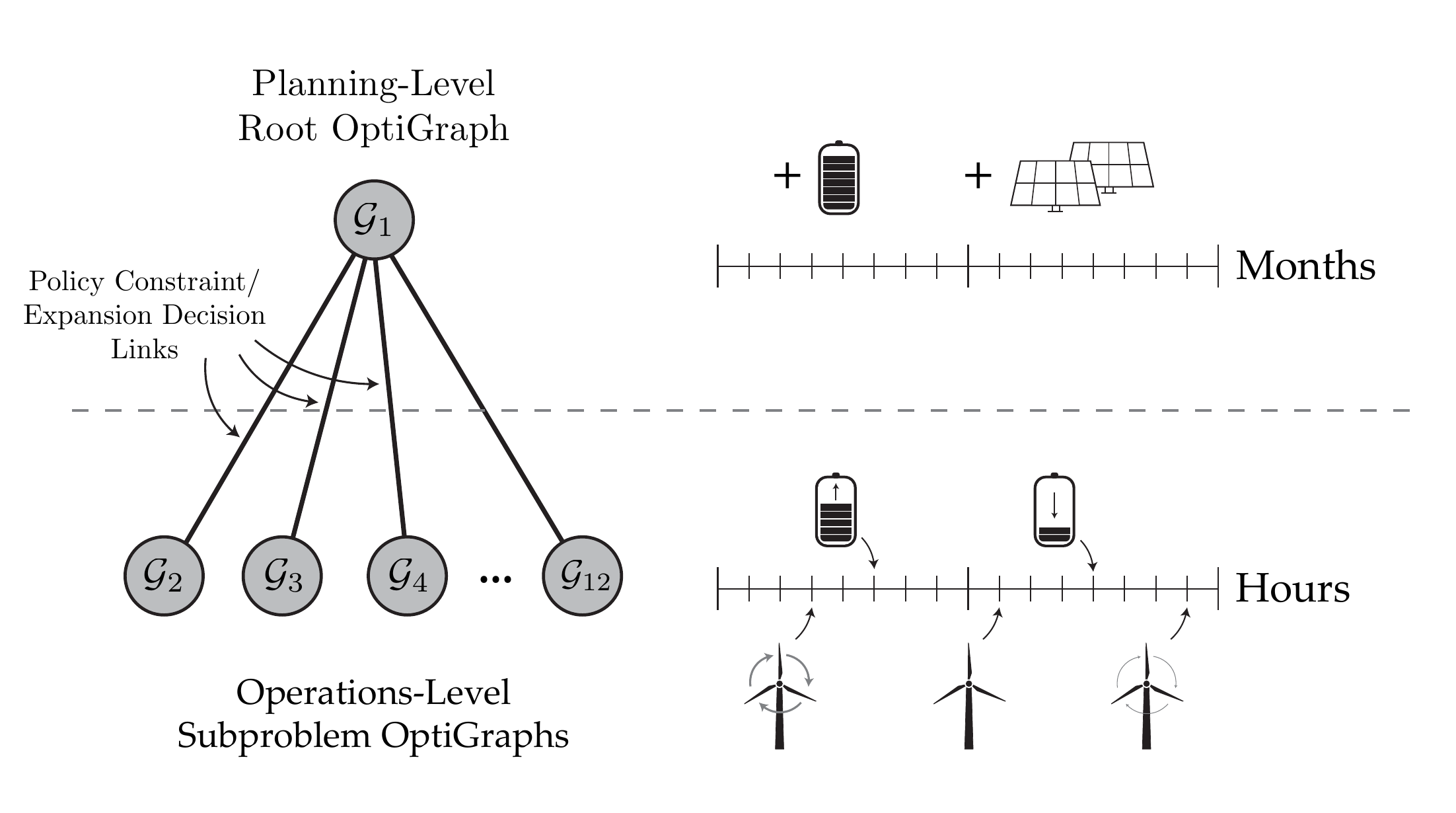}
    \caption{Graph representation of CEM obtained from {\tt GenX.jl}. The problem includes a planning-level containing expansion and retirement decisions and policy constraints which are contained in subgraph $\mg_1$. The operations-level includes different subperiods containing operations variables and constraints. The planning level decisions can be considered as operating over a longer time scale (order of months) than the operations level (order of hours).}
    \label{fig:CS2_graph}
\end{figure}

This problem can be represented as a graph by using a subgraph for the planning-level containing the variables $x_p$ and $|W_o|$ subgraphs for the operations-level (one for each operations subperiod), each with variables $x_w, w \in W_o$. This approach results in \eqref{eq:CEM_link_cons} being stored on edges connecting each operations-level subgraph and the planning-level subgraph. Importantly, the policy constraints \eqref{eq:CEM_policy_cons} form a hyper edge between all the operations-level subgraphs, meaning we do not form the tree structure required for gBD. To address this, budgeting variables $q_{p, w}, w \in W_o$ are added to the planning-level to capture the policy constraints for each operations sub-period; then, \eqref{eq:CEM_policy_cons} is replaced by $\sum_{w \in W_o} q_{p,w} \le d$ which is stored on the planning-level subgraph, and the constraint $Q_w x_w \le q_{p,w}, w \in W_o$ are added to the edges between each operations-level subgraph and the planning-level subgraph. This approach is similar to the lifting procedure shown in Figure \ref{fig:lifting} and was implemented in previous works with {\tt GenX.jl} \cite{jacobson2024computationally,pecci2024regularized}. The resulting graph structure can be seen in Figure \ref{fig:CS2_graph}, where the planning-level subgraph, $\mg_1$, contained 86 variables (the complicating variables) and each operations subgraph contained 12,626 variables. Here, we will set the planning-level subgraph $\mg_1$ as the root subgraph $\mg_r$, as this results in a two-level tree where the second stage is fully parallelizable (see \cite{jacobson2024computationally}. Note however that, for the gBD algorithm, we could set any of the subgraphs as the root subgraph. 

\vspace{0.1in}

\begin{figure}
    \centering
    \includegraphics[width=0.6\linewidth]{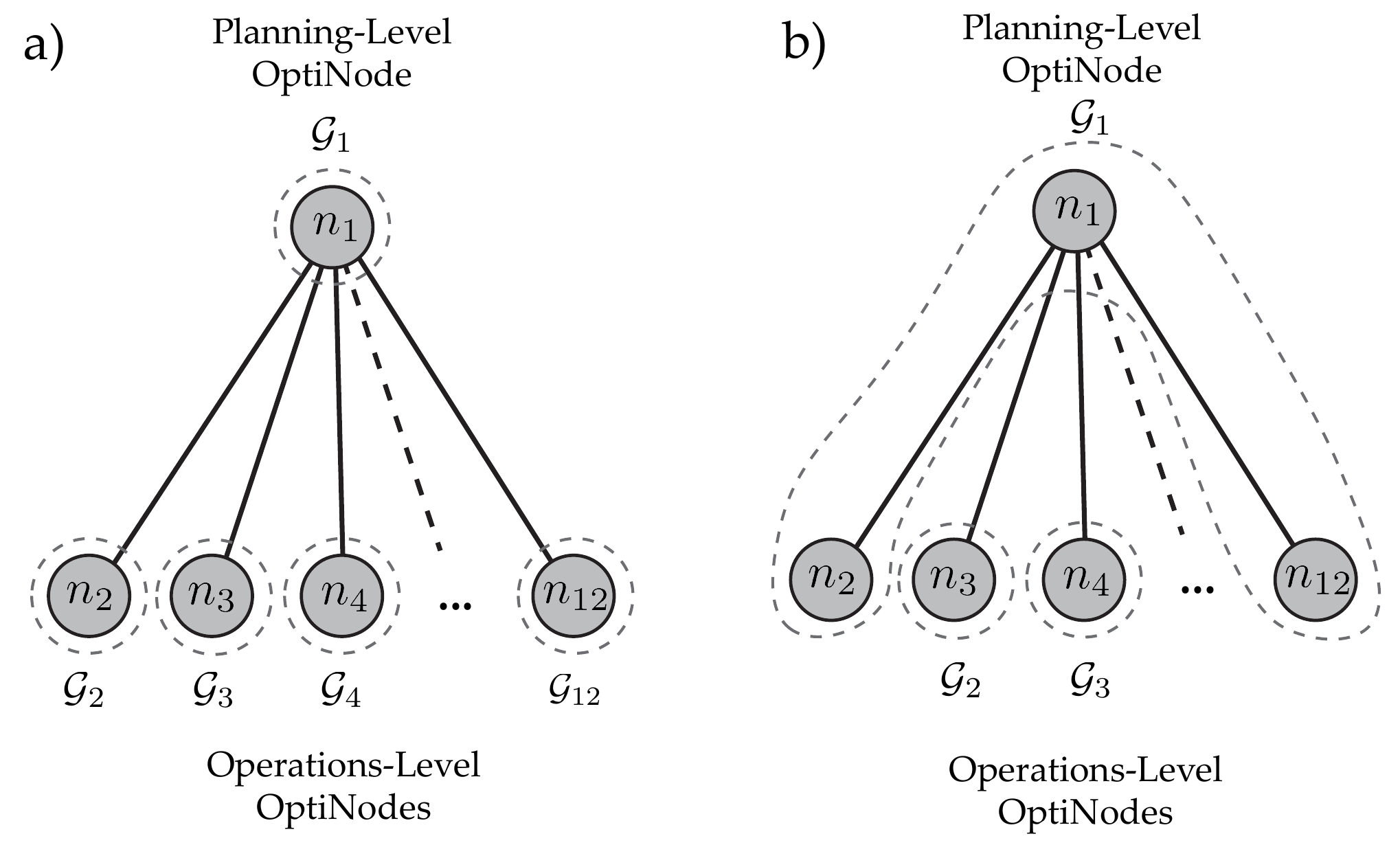}
    \caption{The two different partitioning approaches used for the capacity expansion model. a) shows the original partitioning with a single OptiNode on each subgraph, where there is one planning level optinode and 11 operations level OptiNodes, each representing a different subperiod. b) shows an alternative partitioning where two of the operational OptiNodes are placed with the planning-level OptiNode on the root subgraph.}
    \label{fig:CS2_partitioned}
\end{figure}

\begin{figure}
    \centering
    \includegraphics[width=0.6\linewidth]{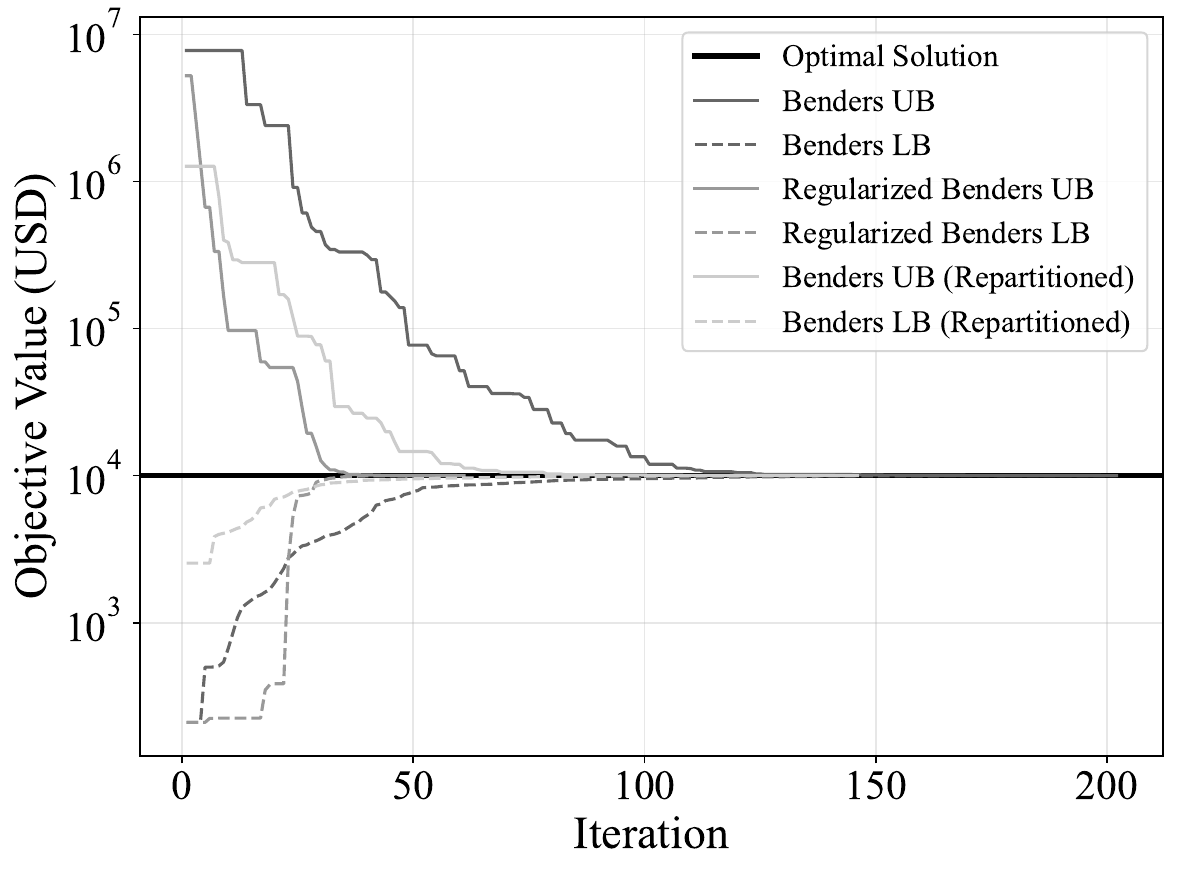}
    \caption{Solution of \eqref{eq:CEM} with gBD using {\tt PlasmoBenders.jl}. The problem is solved without a regularization scheme and each subperiod on its own subproblem subgraph (``Benders''), with a regularization scheme and each subperiod on its own subproblem subgraph (``Regularized Benders''), and with two subperiods included in the root/master problem and without a regularization scheme (``Benders (repartitioned)'').}
    \label{fig:CS2_solutions}
\end{figure}

The problem modeled as a graph was solved in {\tt PlasmoBenders.jl} using the gBD algorithm outlined in Algorithm \ref{alg:graph_NBD}. The planning level and each operational subproblem were each represented as a node and initially placed on separate subgraphs. We set $\mg_1$ (containing the planning problem node) as the root subgraph (i.e., $\mg_r = \mg_1$) which results in two stages with multiple subproblems in the second stage, allowing for the use of the multi-cut approach. The convergence tolerance was set to $1 \times 10^{-4}$. We also solved using a regularization scheme outlined by \cite{pecci2024regularized} and which is implemented generally for the graph-based problem within {\tt PlasmoBenders.jl} (see Supporting Information). In addition, {\tt PlasmoBenders.jl} supports parallelizing the solution of each of the operations sub-periods to speed up overall solve time, and this functionality was used for this problem. Finally, as a third option, we solved the problem by placing two (the first and the last) of the operational subproblem nodes into a subgraph with the planning level node as shown in Figure \ref{fig:CS2_partitioned}. This results in a larger root problem, but one which is ``better informed'' by knowledge of some subproblems \cite{constante2025relaxation}. The results of these three approaches can be seen in Figure \ref{fig:CS2_solutions}. 

\vspace{0.1in}

\vspace{0.1in}

We note as well a challenge that arose in this problem in that the exact number of iterations was not always consistent between different runs. This inconsistency is likely due to degenerate solutions in the root problem due to policy constraints (e.g., the CO$_2$ budgeting variables for the maximum allowable emissions in a subproblem do not show up in the root problem objective function, meaning different CO$_2$ budgets can be chosen with the same root problem objective but different subproblem performance at an early iteration). Consequently, we ran each problem 10 times and recorded the number of iterations used, with the results in Figure \ref{fig:CS2_solutions} being the run with the fifth fewest iterations. These results are all included in the repository, but the three different versions (without regularization, with regularization, and the repartitioned problem with two subperiods included in the root subgraph without regularization) took on average 202, 61, and 139 iterations, respectively, with overall ranges of 149-271, 57-65, and 109-158 iterations.  

Overall, {\tt PlasmoBenders.jl} was able to solve \eqref{eq:CEM} using the gBD algorithm (Algorithm \ref{alg:graph_NBD}) and reached the required tolerance between the upper and lower bounds, serving as a proof-of-concept for the graph-based approach. In addition, {\tt PlasmoBenders.jl} made it trivial to toggle on/off the use of a regularization strategy (by changing a single settings parameter), which substantially reduced the required number of iterations for convergence. This regularization strategy is implemented in a general manner in {\tt PlasmoBenders.jl}, such that it can easily be applied to other BD problems. {\tt Plasmo.jl} also readily supports altering the graph structure, such as moving operational subproblems into the root problem. Doing so for this case study resulted in overall fewer iterations, though the root problem solve time did increase. Under a bespoke implementation of BD, making an alteration to the root problem like this would be non-trivial, but it can be done with just a few lines of code because of {\tt Plasmo.jl}'s flexibility and the gBD algorithm being based solely off of graph structrue. This highlights how the graph approach provides a framework for developing and implementing general algorithms and simplifies modeling and implementation for users.

While this case study has focused on a CEM, many problems can have a similar structure to the problem above. BD has been used in many applications, and the hierarchical or multi-scale nature of this problem like that seen in Figure \ref{fig:CS2_graph} could arise in a variety of areas (see for instance \cite{jemai2013home,kong2017hierarchical,liu2021hierarchical,vegetti2022ontology,zhang2024decomposition}). Thus, the approach used herein could therefore be applied to other problems or disciplines as well. 

\subsection{Case Study 3: Production Cost Modeling}\label{sec:cs3}

This study illustrates how gBD can be used to solve an electricity system production cost model (PCM). PCMs are routinely used in power systems to determine the cost of operating a power grid over a given timeframe (typically the next 24 hours) \cite{kahn1995regulation,lara2024powersimulations} by considering the dispatch and operation of a set of generators (e.g., fossil fuel, wind, solar) to meet demand in a region. These models also frequently include a degree of ``look-ahead'' which allows the model to consider expected future demand and  avoid myopic decisions (e.g., the model may choose to leave some generators turned on to meet expected demand arising after the next 24 hour period). These dispatch and operation decisions are subject to several constraints, such as up/down time and ramping constraints. 

\begin{figure}
    \centering
    \includegraphics[width=1.0\textwidth]{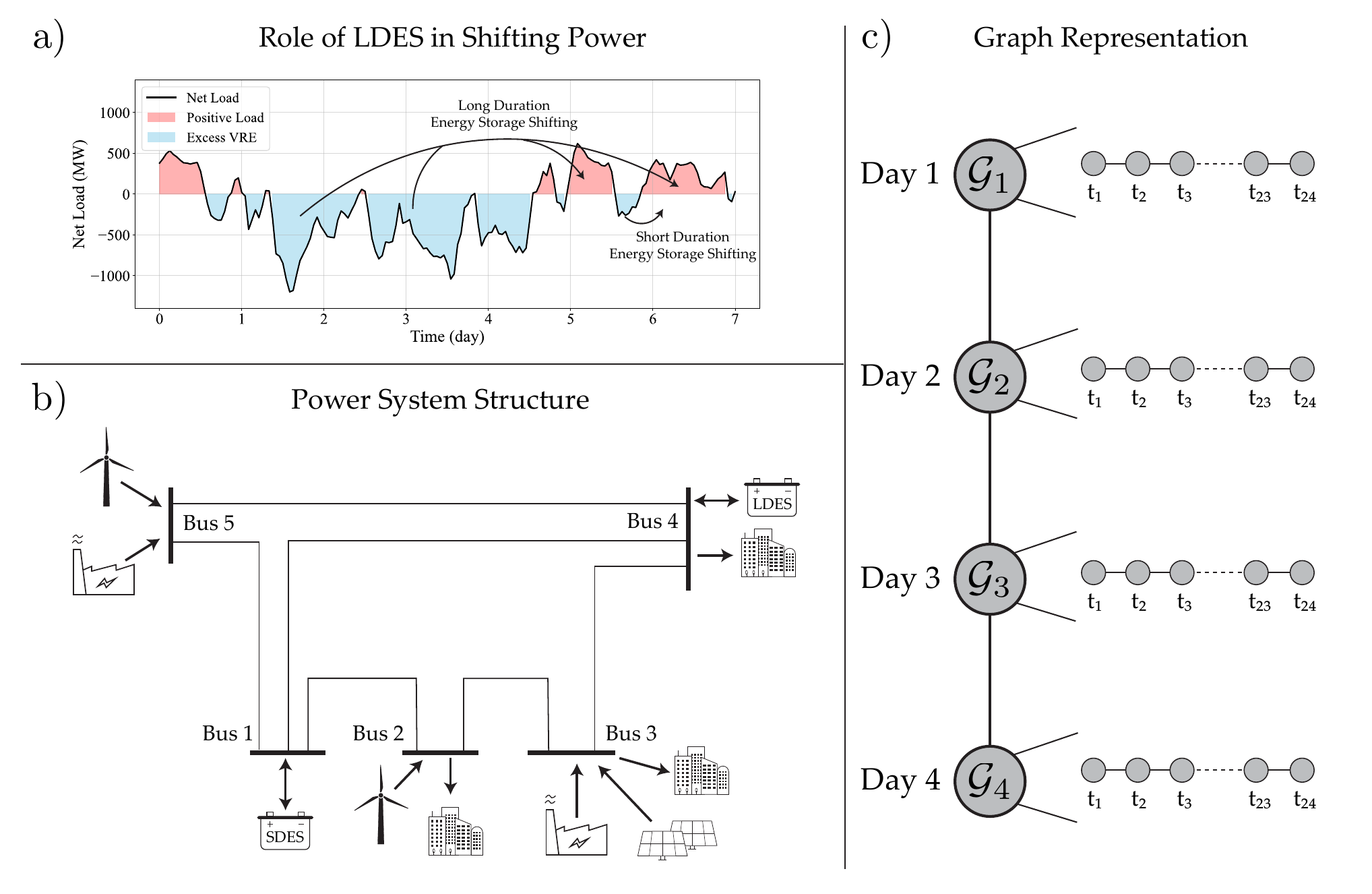}
    \caption{A visualization highlighting PCMs and how they can be represented as a graph. LDES devices copmlicate PCM operation because they play the role of shifting power across longer time scales (a) for a power system like the one shown in (b). We can represent these problems as a graph (c), where each subgraph represents a period of time (forming a coarse time resolution among the subgraphs). Each time period is effectively a hierarchical decision that impacts the subsequent time period, resulting in the linear graph structure shown. Adapted from \cite{cole2025toward}.}
    \label{fig:PCM_visual}
\end{figure}

The hierarchy in this problem arises due to the inclusion of a long-duration energy storage (LDES) device in the system. LDES devices are important for high-renewable grids \cite{blair2022storage,denholm2023moving,guerra2021beyond}, but they complicate the operation of the grid (i.e., the solution of PCMs) because they play the role of shifting power across larger time scales so capturing their true value can require looking farther into the future (see Figure \ref{fig:PCM_visual}a). Unfortunately, adding longer look-ahead horizons to PCMs can increase required computational time and resources. For instance, many PCMs operate with a total horizon of 48 hours (a 24-hour day ahead period with an additional 24 hours of ``look-ahead''), but PCMs with LDES can require using look-ahead horizons of 3 or more days which can lead to an exponential increase in required solution time \cite{guerra2025towards}. In this way, the presence of the LDES device creates a hierarchy in decision time scales: the PCM operation on the order of hours and the LDES operation on the order of days. We capture this hierarchy by representing a series of time points as subgraphs (e.g., representing a sequence of hours or days as a single subproblem) that are connected to the next series of time points, with the battery state of charge being linked across the subgraphs; in other words, one could consider that the subgraphs are representing the time scale on which the LDES operates, and each subgraph contains a finer resolution problem on which the PCM operates. This creates the linear graph hierarchical structure seen in Figure \ref{fig:NBD_graph}. An example of this hierarchical structure is presented in Figure \ref{fig:PCM_visual}c, which is applied to a system such as that of Figure \ref{fig:PCM_visual}b.

The problems we consider are implemented in the {\tt Julia} package \texttt{\mbox{\hspace{60pt}} \allowbreak PowerSimulations.jl} \cite{lara2024powersimulations} and converted into a {\tt Plasmo.jl} OptiGraph. \texttt{\mbox{\hspace{50pt}} \allowbreak PowerSimulations.jl}  constructs an optimization model for a power system given by the user. Based on the power system data and components (e.g., generators, loads, transmission lines), {\tt PowerSimulations.jl} defines the needed objective function, decision variables, and constraints to represent the PCM of the system. \texttt{\mbox{\hspace{40pt}} \allowbreak PowerSimulations.jl}  also provides various options for constraints and variables used for modeling components (for instance, a thermal generator can be modeled in several ways). Further details of {\tt PowerSimulations.jl} can be found in the package documentation and in \cite{lara2024powersimulations}. 

The general formulation for the PCMs we used in this case study is described below. As the purpose of this work is to highlight how graphs can be used for representing and decomposing hierarchical problems, we will not introduce every variable and constraint and instead only give a high-level overview. The PCMs we consider are MILPs that seek to minimize cost, including production and start-up/shut-down costs, while also considering penalties on slack variables (which ensure feasibility, such as by allowing load shedding). Here, production costs apply to thermal generators (e.g., coal, gas, nuclear) and consider the cost of producing energy with any given generator. Discrete startup-shutdown decisions apply only to thermal generators and incur costs associated with turning on or shutting down a generator. The PCMs are subject to constraints including generator ramp up/down, generator minimum up/down time constraints after start-up/shut-down, storage operating constraints, reserve requirements, transmission constraints, and overall power balance constraints. Transmission constraints included constraints on power flow throughout the system using a linearized DC approximation of power flow. There are also constraints on the battery operation, with the state of charge of the battery linked across every time point. The full sets of constraints are dependent in part on the system data; for instance, {\tt PowerSimulations.jl} is built on the {\tt PowerSystems.jl} data structure, which supports different types of thermal generators, including standard generators and multi-start generators, so the exact problem formulation is dependent on the system defined.

The examples below include two different systems: a 5-bus system and the 73-bus Reliability Test System (RTS; adapted from \cite{li2011testing} and \cite{gonzalez2013composite}, respectively). Details of these systems can be found at \url{https://github.com/NREL/LDES-Sizing-and-siting-model}. The 5-bus system included five buses, six transmission lines, two thermal generators, one solar generator, two wind generators, and three loads. The RTS included 73 buses, 120 transmission lines and tap transformers, 54 thermal generators, 58 solar generators, five wind generators, and 51 loads. Each system included a short-duration (3.6 hr) and long-duration (13.5 hr) energy storage device. The presence of LDES in the system requires using longer look-ahead horizons in the PCM simulation.

The {\tt PowerSimulations.jl} model can be represented as a graph by using nodes to represent buses and arcs at each time point (arcs are {\tt PowerSystems.jl} objects that can contain lines or tap transformers; for our purposes, they are containers to hold transmission lines). Each generator, storage device, and load are contained on a specific bus, and each transmission line is contained on an arc. Thus, each node corresponding to a bus will contain variables for any generator, storage device, or load that are contained on that bus (e.g., active power variables, on variables). Each node corresponding to an arc will contain variables for any transmission line on that arc. As with Case Study 1, each node corresponds to a single point in time, such that for a system with $N_B$ buses, $N_A$ arcs, and $N_T$ time points, there will be $(N_B + N_A) \times N_T$ total nodes. Edges are placed between nodes based on constraints in the {\tt PowerSimulations.jl} model. For instance, edges are placed between nodes for an arc and two buses if that arc has a line that transports power between those two nodes. Edges are also placed across time, such that nodes at different time points are connected; for instance, a node representing a bus containing a generator with up-time limits will be linked to next and previous time points (and if the up-time limit is greater than 1, it will be linked to multiple nodes in both forward and backward directions). 

\begin{figure}
    \centering
    \includegraphics[width=\textwidth]{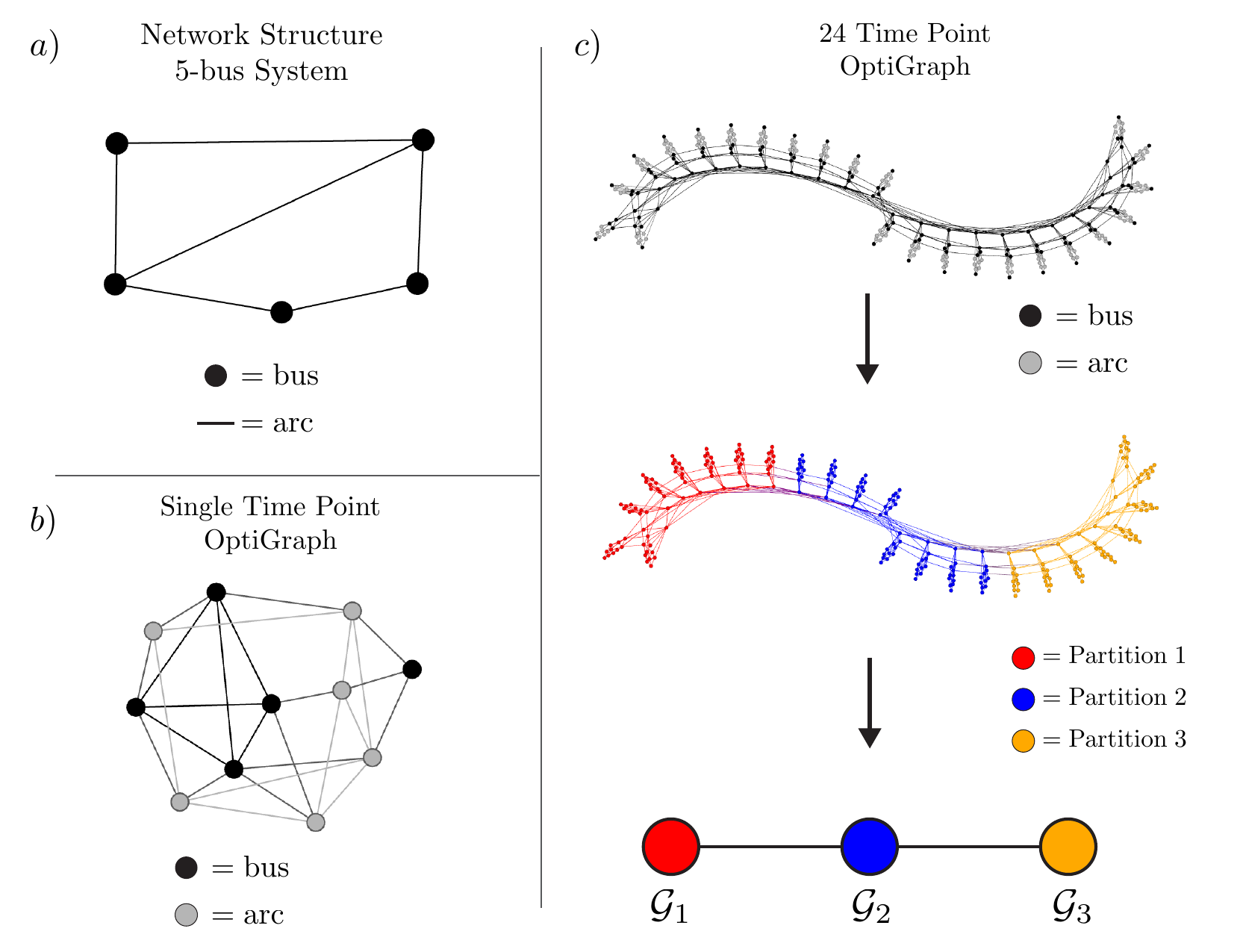}
    \caption{Example of the 5-bus system used in Case Study 3, represented as a graph, including a) the basic network structure with buses represented by a node and system arcs (transmission lines) represented by edges; b) a single time point OptiGraph with buses and arcs represented by nodes (each node can contain variables and constraints); c) a 24 time point OptiGraph where each time point is also linked across time (e.g., via up/down time constraints, storage device state of charge constraints). This results in a graph with many cycles, but the graph can be partitioned to form the subgraph structure shown on the bottom, to which gBD can be applied. Adapted from \cite{cole2025toward}}.
    \label{fig:5bus_system}
\end{figure}

\begin{figure}
    \centering
    \includegraphics[width=1.0\textwidth]{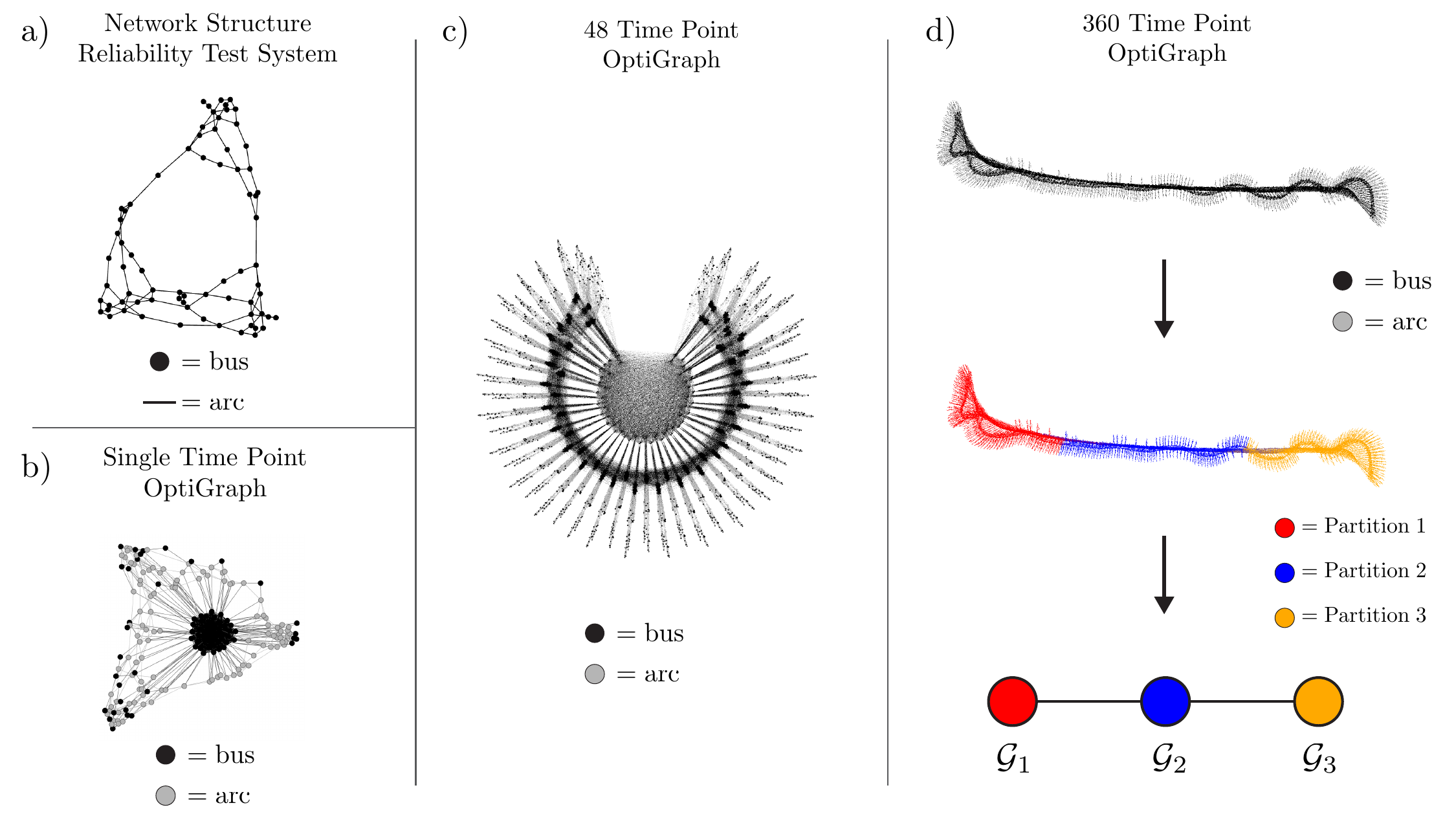}
    \caption{Example of the RTS used in Case Study 3, represented as a graph, including a) the basic network structure with buses represented by a node and system arcs (transmission lines) represented by edges; b) a single time point OptiGraph with buses and arcs represented by nodes (each node can contain variables and constraints); c) a 48 time point OptiGraph where each time point is also linked across time (e.g., via up/down time constraints, storage device state of charge constraints); d) a 360 time point OptiGraph with more than 65,000 nodes, highlighting the complexity that can be present in these types of problems. This graph contains many cycles, but can be partitioned to form the subgraph structure shown on the bottom, to which gBD can be applied.}
    \label{fig:RTS_system}
\end{figure}

The 5-bus system and RTS used in this case study can both be represented as a graph as explained above, and visualizations of these systems are given in Figures \ref{fig:5bus_system} and \ref{fig:RTS_system}. Figures \ref{fig:5bus_system}a and \ref{fig:RTS_system}a show the basic system structure with each bus represented by a node and arcs (containing transmission lines) represented by edges (this is the most intuitive graph representation and is not the OptiGraph). Figures \ref{fig:5bus_system}b and \ref{fig:RTS_system}b show the OptiGraph representation for a single time point using the methods outlined above. For the single time point graph, the edges are linking constraints which capture flow constraints as well as overall system energy balance constraints (which includes all nodes containing generator variables; each bus also contains a local balance for power flowing in and out of the node). Figures \ref{fig:5bus_system}c, \ref{fig:RTS_system}c, and \ref{fig:RTS_system}d show the system expanded across time, where the linking constraints now also include the state of charge of the battery, ramping constraints, and up- and down-time constraints. Note here that Figure \ref{fig:RTS_system}c shows the RTS with 48 time points and the down-time constraint can best be seen in this figure: two generators in this system have down times of 48-hours, and the central ``cloud'' in the middle of this graph are the down-time constraints, where each of the nodes containing these generator variables are linked to each other since the down-time limit is equal to the time span the entire graph.  Importantly, none of these expansions results in the required structures for gBD. Figures \ref{fig:5bus_system}c and \ref{fig:RTS_system}d show an OptiGraph with a longer horizon (24-hours and 15-days for the 5-bus system and RTS, respectively), which can then be partitioned across time to form the subgraph structure shown in the bottom of the subfigure. This resulting subgraph structure can be solved using gBD. 

\begin{figure}
    \centering
    \includegraphics[width=\linewidth]{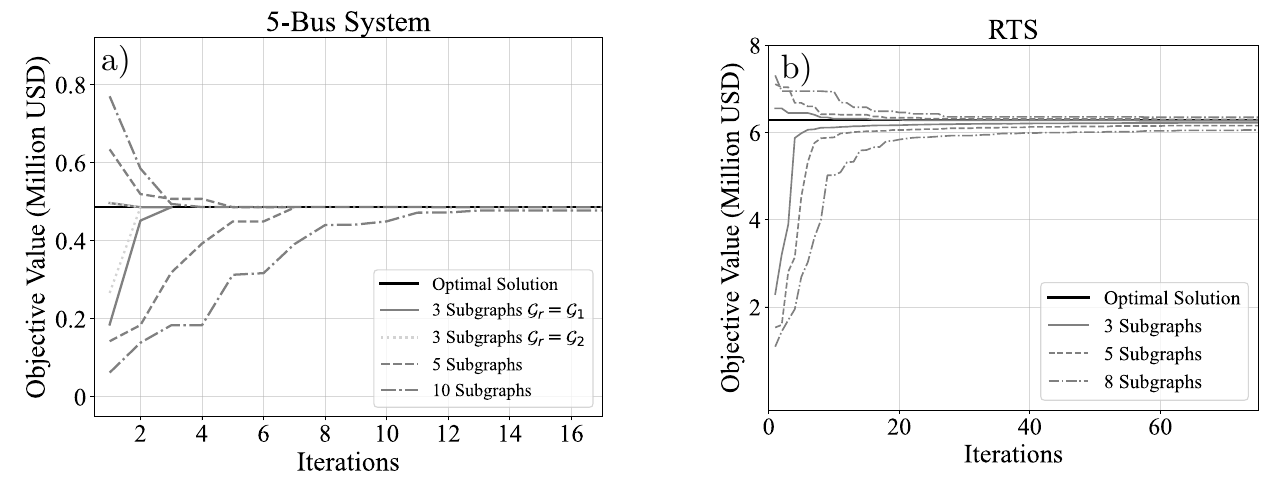}
    \caption{Results for the gBD approach for a) the 5-bus system with a 90-day time horizon using 3, 5, and 10 partitions, and b) the RTS with a 15-day time horizon using 3, 5, and 8 partitions. gBD was applied by setting the first subgraph (containing the first time point of the model) to be the root subgraph. The lines above the black line represent the upper bound at each iteration and the lines below the black line represent the lower bound at each iteration. For the 5-bus system, we also tested gBD using both $\mg_1$ and $\mg_2$ as the root subgraphs.}
    \label{fig:CS3_solutions}
\end{figure}

To highlight how gBD can be applied to graph-based problems, we solved an instance of the 5-bus system and the RTS using {\tt PlasmoBenders}' implementation of gBD (see Algorithm \ref{alg:graph_NBD}) with the backwards pass parallelized (i.e., when solving stage $i$ of iteration $k$ in the backwards pass, we do not wait for the cut from stage $i+1$ and add cuts to all problems at once when all problems have finished solving). Each system included 1-hour time steps. The 5-bus system used a 90-day total horizon (2160 time points) and the RTS used a 15-day total horizon (360 time points). The 5-bus system included 127,450 variables (17,280 binary) and 157,680 constraints while the RTS included 406,655 variables (59,040 binary) and 360,360 constraints (constraints do not include variable bounds). Each system was structured as a graph, and the convergence tolerance was set to be 0.5\%. To create the structure required for gBD, we must partition the nodes into subgraphs; the partitions used are arbitrary, so long as they result in a tree structure in the subgraphs (see for instance the bottom graph in Figures \ref{fig:5bus_system}c and \ref{fig:RTS_system}d). We chose to test different partition sizes and report the results for the gBD solution algorithm in Figure \ref{fig:CS3_solutions}. For the 5-bus system, we used 3, 5, and 10 partitions (subgraphs), where each subgraph of the partition represented consecutive time points (i.e., partition sizes in the three cases were 30 days, 15 days, and 9 days, respectively). For the RTS system, we used 3, 5, and 8 partitions (i.e., partition sizes were 5, 3, and 2 days, respectively; note that the final partition in the 8-partition case only contained 1 day). Each of these partitioning decisions resulted in a linear graph, and we set the subgraph containing the first time point as the root subgraph in each case (e.g., for the subgraph structure in Figure \ref{fig:RTS_system}d, we set $\mg_1$ as the root subgraph). To highlight how the decision of root subgraph can impact the algorithm, we also tested two different root subgraphs for the 3-partition case of the 5-bus system by solving with a root subgraph containing the initial time point of the PCM ($\mg_r = \mg_1$) and the ``middle'' subgraph (e.g., $\mg_2$ in Figure \ref{fig:5bus_system}c; additional slack variables had to be added to this problem to ensure recourse, but the final solution was equal to that of the problems without these slack variables). 

For both the 5-bus system and RTS, the gBD algorithm in {\tt PlasmoBenders} was able to get to a high quality solution, though partitioning decisions impacted convergence (Figure \ref{fig:CS3_solutions}; see Tables 1 and 2 in the Supporting Information for CPU time and memory usage). For the 5-bus system, the 3 subgraph problem converged in only 3 iterations, while the problem with 5 subgraphs required 7 iterations. With 10 subgraphs, the problem was not able to close the gap beyond a gap of 1.7\% (which it reached after 17 iterations) even if it ran for more than 50 iterations. Interestingly, using the middle subgraph as the root subgraph in the 3 subgraph case was able to converge in only 2 iterations. For the RTS, the problem did not converge to the 0.5\% tolerance for any of the partitionings within the set $K_{max}$ iterations (for each partitioning, the gap was still closing very slowly at termination). The 3, 5, and 8 subgraph problems only reached a gap of 0.7\%, 2.1\%, and 4.7\% (for $K_{max}$ of 75), respectively. This may be due in part to the presence of a duality gap in the mixed integer problem. 

Overall, the graph-based approach was effective at constructing problems that can be decomposed by gBD. Representing these problems under the OptiGraph abstraction can give flexibility in applying different graph-based decomposition schemes (such as enabling different numbers of partitions or selection of alternative root subgraphs as shown above) and could lead to non-intuitive applications of algorithms. For instance, changing the root subgraph of the 3 subgraph case of the 5-bus system was able to converge effectively. Temporal decompositions like those above are not typically solved in this way, but by representing these problems as graphs, we can identify alternative solution approaches.

\section{Conclusions and Future Work}\label{sec:future_work}

We have presented a graph abstraction for modeling and solving optimization problems that contain hierarchical optimization problems. The OptiGraph abstraction, implemented in {\tt Plasmo.jl}, allows for constructing hierarchical subgraphs to capture hierarchical relationships in optimization problems. The OptiGraph abstraction also gives significant flexibility by enabling partitioning and aggregation of graph structures. We have shown how the OptiGraph can can be used as a foundational abstraction for developing generalized solution algorithms and for implementing approximation schemes. We have given an example of generalizing algorithms to a graph with graph-based Benders decomposition, gBD, which is implemented in the package {\tt PlasmoBenders}. By using this graphical approach, it is possible to identify alternative non-intuitive applications of algorithms or flexibly reconfigure problem structure to improve convergence properties or better leverage available computing resources (e.g., parallel computing hardware). We have provided different case studies showing how the OptiGraph abstraction can be used for modeling and solving hierarchical problems.

There are several areas we would like to explore in future work. There are possible extensions of gBD that can be explored. For instance, there have been different works presenting parallel implementations of nested BD \cite{colonetti2022parallel,dos2023accelerated,rahmaniani2024asynchronous,santos2016new}, such as by solving subproblems asynchronously. Similar extensions could be applied to graphs, or new algorithms could be developed. For example, if the linear graph in Figure \ref{fig:NBD_graph} could be partitioned such that all odd numbered subgraphs are in one partition and all even numbered subgraphs in another partition (this would form a bipartite graph); this would be the structure required for the traditional BD application, but within each subgraph, there would be two or three additional subgraphs with no original connections to other subgraphs in their partition. It is possible that, since there are no connections of subgraphs within the two partitions, a parallel algorithm could be more readily developed. We would also like to improve the gBD implementation in {\tt PlasmoBenders.jl} to parallelize the solution within stages (currently, this is only supported for problems with two stages), and we would like to implement feasibility cuts for multi-stage problems. We are also interested in considering how the structuring decisions impact the strength of cutting planes. As an example, for the linear graph in Figure \ref{fig:NBD_graph}, by changing the root subgraph from $\mg_1$ to $\mg_3$, we also could change the number of cutting planes added to $\mg_3$. When $\mg_r = \mg_1$, $\mg_3$ would only have cutting planes from the solution of $\mg_4$, but when $\mg_r = \mg_3$, then $\mg_3$ would also receive cutting planes from the solution of $\mg_2$. Understanding how to choose the root subgraph and how these cutting planes impact the solution could be informative. Fourth, we would like to extend gBD to nonlinear problems, such as by following generalized BD frameworks \cite{geoffrion1972generalized}. We would also like to explore how regularization schemes could be applied within gBD for problems with many stages as {\tt PlasmoBenders.jl} currently only applies regularization to two stage problems. Finally, we are interested in adding functionality for different types of cutting planes or other advancements (such as the adaptive oracles of \cite{zhang2024stabilised}) within {\tt PlasmoBenders.jl}.

 At present, {\tt Plasmo.jl} does not readily facilitate distributing memory with the OptiGraphs, but this can be important for large-scale decomposition problems on high-performance compute clusters or for cloud-based computing. We would like to develop a distributed OptiGraph abstraction for distributing memory and solution of large-scale problems. Second, we would like to consider different applications of the hierarchical graphs; one such application would be stochastic programming, and in particular, multi-horizon stochastic programming (MHSP) \cite{kaut2014multi,zhang2024decomposition,zhang2024stabilised}. In MHSP, these problems are frequently visualized as graphs and decomposition approaches like BD can be required for solving them. Third, we are interested in implementing additional decomposition schemes for OptiGraphs, such as developing a graph-based  progressive hedging, ADMM. or other Lagrangian decompositions. Finally, these hierarchical graphs could facilitate including surrogate models into an optimization problem (such as replacing a single subgraph or a hierarchical layer with a neural network). Surrogate models have been embedded in optimization problems in other works \cite{ceccon2022omlt,henao2011surrogate,tjeng2017evaluating}, but the hierarchical subgraphs could provide a powerful tool for including these kinds of models. 

\section*{Acknowledgements}
This material is based on work supported by the U.S. Department of Energy under grant DE-0002722. We also acknowledge support from the National Science Foundation under award CBET-2328160. F.P also received support from the Eric \& Wendy Schmidt Fund for Strategic Innovation (grant G-23-64740) and the Princeton Zero-carbon Technology Consortium, supported by unrestricted gifts from Google and Breakthrough Energy. We also thank Dr. Jordan Jalving for his technical support and continued maintenance of {\tt Plasmo.jl}.

\bibliography{PlasmoDecompositions}

\end{document}


\title{{\bf Supporting Information}\\ {Graph-Based Modeling and Decomposition of Hierarchical Optimization Problems}}

\author{David L. Cole$^{1}$, Filippo Pecci$^{2}$, Omar J. Guerra$^{3}$, Harsha Gangammanavar$^{4}$, Jesse D. Jenkins$^{2, 5}$, and Victor M. Zavala$^{1, 6}$\\
 {\small $^{1}$ Department of Chemical and Biological Engineering}\\
 {\small University of Wisconsin-Madison, Madison, WI 53706, USA}\\
 {\small $^{2}$ Andlinger Center for Energy and Environment}\\
 {\small Princeton University, Princeton, NJ 08540, USA}\\
 {\small $^{3}$ Grid Planning and Analysis Center}\\
 {\small National Renewable Energy Laboratory, Golden, CO 80401, USA}\\
 {\small $^{4}$ Operations Research and Engineering Management}\\
 {\small Southern Methodist University, Dallas, TX 75205, USA}\\
 {\small $^{5}$ Department of Mechanical and Aerospace Engineering}\\
 {\small Princeton University, Princeton, NJ 08540, USA}\\
 {\small $^{6}$ Mathematics and Computer Science Division}\\
 {\small Argonne National Laboratory, Lemont, IL 60439, USA}\\
}


\maketitle
\tableofcontents

\section{Software Implementation for OptiGraph Structure Manipulation}

{\tt Plasmo.jl} includes functionality for automatically performing partitioning and aggregation for OptiGraphs. Partitioning is done by creating a {\tt Plasmo.Partition} object and then either applying that partition to an existing graph or creating a new graph from that partition. In {\tt Plasmo.jl}, the {\tt Partition} object contains lists of OptiNodes, OptiEdges, and subpartitions (i.e., additional lists of nodes, edges, and subpartitions). This is analagous to the structure of an OptiGraph given in \eqref{eq:optigraph}, where the lists of nodes and edges of the partition object are similar to $\mn_{\mg}$ and $\me_{\mg}$ of \eqref{eq:optigraph}. Consequently, The partition object can be used for either creating new subgraphs on an existing graph or for creating a new graph with the partitioned structure. The {\tt Partition} object requires identifying which nodes should belong to which partition. This is done by passing either i) a vector of integers with a length equal to the number of nodes in the graph that maps each node to a partition number, or ii) a vector of OptiNode vectors, where each OptiNode vector contains the nodes that belong to the same partition. In each case, the number of partitions is determined by the highest integer value or by the number of OptiNode vectors. 
\\

\begin{figure*}[!ht]
    \centering
    \includegraphics[width = 0.6\textwidth]{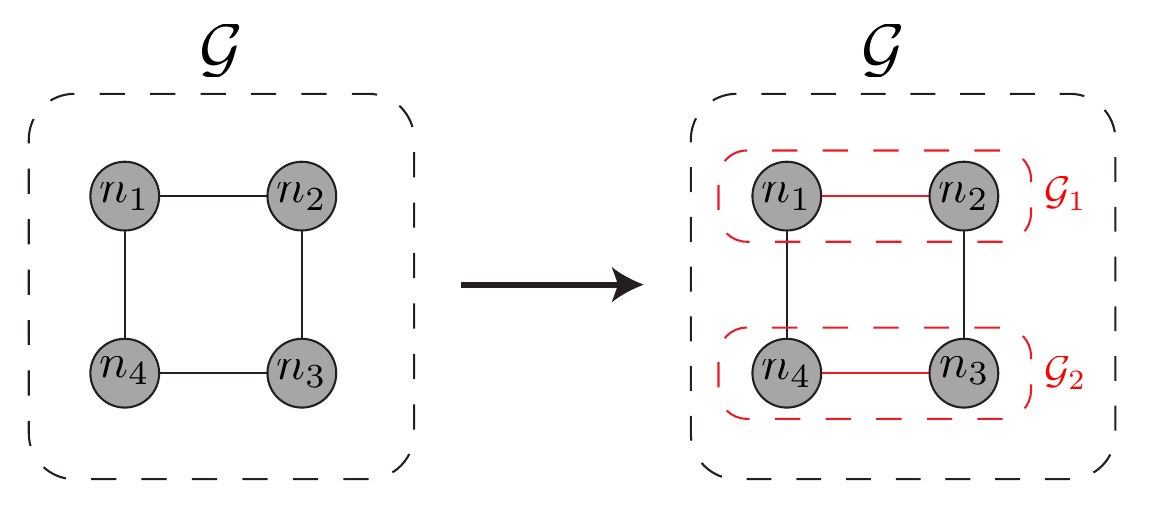}
        \vspace{-0.2in}
    \caption{Visualization of a graph and subsequent partitioning created by Code Snippet \ref{code:partitioning}.}
    \label{fig:partition_example}
\end{figure*}

Code Snippet \ref{code:partitioning} shows how partitioning is performed for an OptiGraph. The graph is defined on Lines \ref{line:part_graph_start} - \ref{line:part_graph_end}. The {\tt Plasmo.Partition} function can be called by either defining a node membership vector (Line \ref{line:node_mem_vector}) or a vector of Optinode vectors (Lines \ref{line:optinode_vec1} - \ref{line:optinode_vec2}), which are passed as arguments to the partition function along with the OptiGraph (Lines \ref{line:partition1} and \ref{line:partition2}). The partitioning object can then be used to assemble a new OptiGraph object using the {\tt assemble\_optigraph} (Line \ref{line:assemble_optigraph}) function or can be applied directly to the original OptiGraph by calling {\tt apply\_partition!} (Line \ref{line:apply_partition}). These functions result in a hierarchical graph with a couple of subgraphs, one with {\tt nodes[1]} and {\tt nodes[2]} as well as the edge between them (Line \ref{line:partition_edge1}) and the other with {\tt nodes[3]} and {\tt nodes[4]} and the edge between them (Line \ref{line:partition_edge3}). Note that the primary graph will also contain two edges (Lines \ref{line:partition_edge2} and \ref{line:part_graph_end}). This process is visualized in Figure \ref{fig:partition_example}.

\begin{figure}[!htp]
    \begin{minipage}[t]{0.7\linewidth}
        \begin{scriptsize}
        \lstset{language=Julia, breaklines = true}
        \begin{lstlisting}[label=code:partitioning, caption = {Code for partitioning a graph into subgraphs using either the {\tt assemble\_optigraph} function or the {\tt apply\_partition!} function.}] 
# Define OptiGraph and nodes
g = OptiGraph() |\label{line:part_graph_start}|
@optinode(g, nodes[1:4])

# Add variables and objective to nodes
for node in nodes
    @variable(node, 0 <= x)
    @objective(node, Min, x)
end

# Add linking constraints between nodes
@linkconstraint(g, nodes[1][:x] + nodes[2][:x] >= 1) |\label{line:partition_edge1}|
@linkconstraint(g, nodes[2][:x] + nodes[3][:x] >= 1) |\label{line:partition_edge2}|
@linkconstraint(g, nodes[3][:x] + nodes[4][:x] >= 1) |\label{line:partition_edge3}|
@linkconstraint(g, nodes[4][:x] + nodes[1][:x] >= 1) |\label{line:part_graph_end}|

# Define node-membership vector and create its accompanying partition
node_membership_vector = [1, 1, 2, 2] |\label{line:node_mem_vector}|
partition = Plasmo.Partition(g, node_membership_vector) |\label{line:partition1}|

# Build an OptiGraph with subgraphs based on the partition
partitioned_g = Plasmo.assemble_optigraph(partition) |\label{line:assemble_optigraph}|

# Define a vector of OptiNodes and create its accompanying partition
optinode_vectors = [ |\label{line:optinode_vec1}|
    [nodes[1], nodes[2]],
    [nodes[3], nodes[4]]
] |\label{line:optinode_vec2}|
partition = Plasmo.Partition(g, optinode_vectors) |\label{line:partition2}|

# Apply the partition to the original graph
apply_partition!(g, partition) |\label{line:apply_partition}|
        \end{lstlisting}
        \end{scriptsize}
    \end{minipage}
\end{figure}

\vspace{0.1in}

Aggregation is also implemented in Plasmo.jl via a couple of functions: {\tt aggregate} and \texttt{\mbox{\hspace{30pt}} \allowbreak aggregate\_to\_depth}. These functions are more specific than the mathematical definitions of aggregation above. The {\tt aggregate} function will aggregate a given subgraph into a new OptiGraph containing a single node (i.e., for an OptiGraph $\mg$, this aggregation is equivalent to setting $\mn_{agg} = \mn(\mg)$ resulting in $\mg^{agg}(\emptyset, \{n_{agg}\}, \emptyset)$). The {\tt aggregate\_to\_depth} function takes a couple of arguments: the OptiGraph to which aggregation will be applied and the level at which to apply aggregation (this level is equal to the depth minus one, as will be shown shortly). This function aggregates each subgraph at a given level into single nodes (e.g., at level 0 on graph $\mg$, each subgraph in $\msg(\mg)$ will be aggregated into a separate node). 
\\

\begin{figure*}[!ht]
    \centering
    \includegraphics[width = 0.6\textwidth]{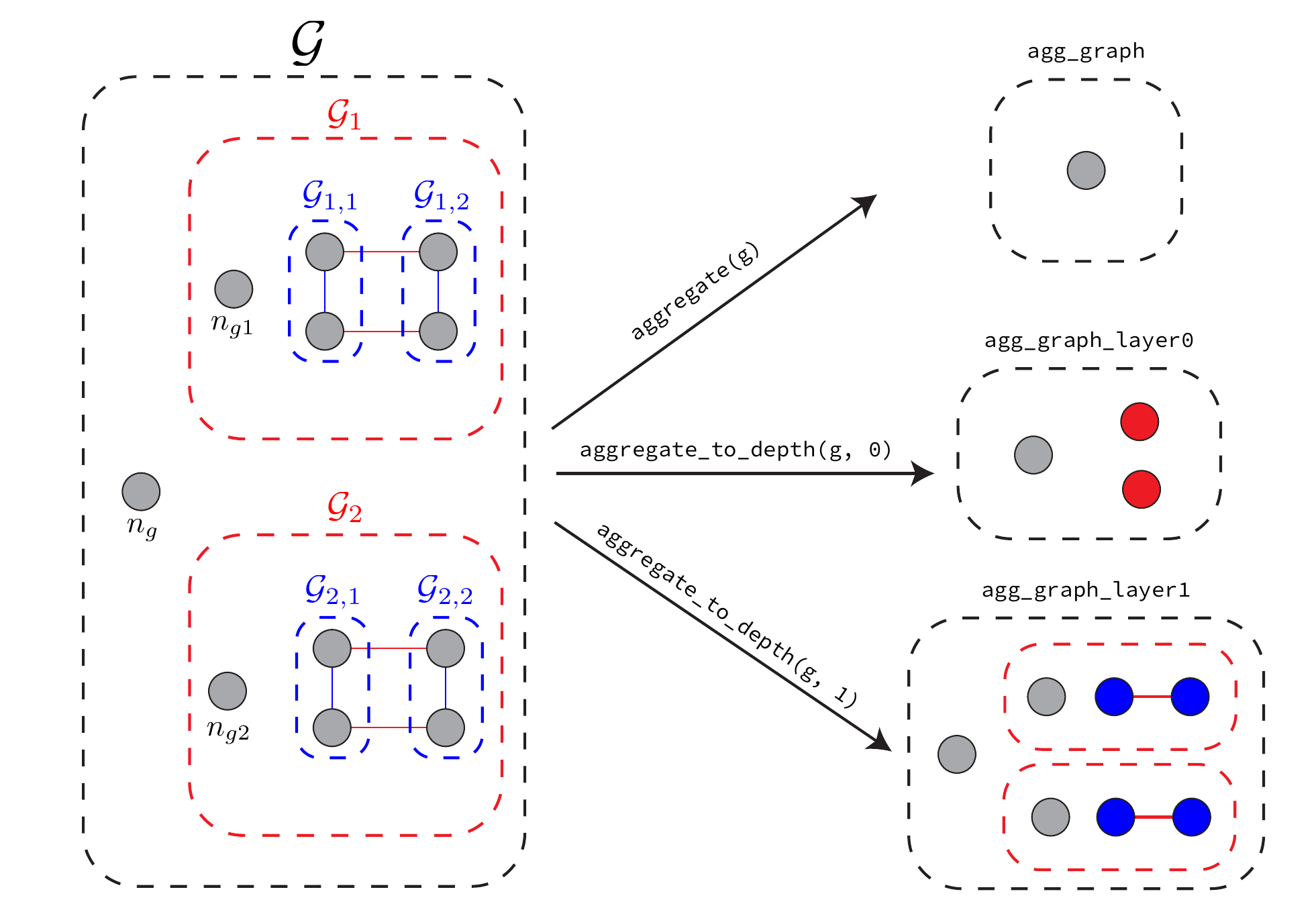}
        \vspace{-0.2in}
    \caption{An example of aggregation approaches implemented in {\tt Plasmo.jl} with the functions {\tt aggregate} and {\tt aggregate\_to\_depth}. The graph $\mg$ on the left has a depth of two, with the first layer of subgraphs colored by red (along with their owning edges) and the second layer of subgraphs colored by blue. Subgraphs can be aggregated into individual nodes, where levels of aggregation are shown on the right, where the nodes formed by aggregating the subgraphs are also color-coded based on their corresponding subgraph. }
    \label{fig:aggregation}
\end{figure*}

To highlight how aggregation functions operate, we provide an example with an accompanying code snippet and visual. We consider the graph $\mg(\{\mg_1, \mg_2\}, \{n_g \}, \emptyset)$, where $\mg_1$ can be written as $\mg_1(\{\mg_{1,1}, \mg_{1,2} \}, \{ n_{g1}\}, \me_{\mg_1})$ and $\mg_2$ can be written as $\mg_2(\{\mg_{2,1}, \mg_{2,2}\}, \{ n_{g2}\}, \me_{\mg_2})$. Code Snippet \ref{code:aggregation} constructs the graph $\mg$ and then creates three different aggregated graphs using the functions above. First, we define a function for building the graphs $\mg_1$ (containing $\mg_{1,1}$ and $\mg_{1,2}$) and $\mg_2$ (containing $\mg_{2,1}$ and $\mg_{2,2}$; Lines \ref{line:agg_function_start} - \ref{line:agg_function_end}). These graphs, along with $\mg$, are then instantiated (Lines \ref{line:agg_build_g} - \ref{line:agg_build_g2}). Graphs $\mg_1$ and $\mg_2$ are then added as subgraphs to $\mg$ (Lines \ref{line:agg_add_g1} and \ref{line:agg_add_g2}), and the nodes $n_g$, $n_{g1}$, and $n_{g2}$ are added to graphs $\mg$, $\mg_1$, and $\mg_2$ respectively. The resulting graph can be seen on the left hand side of Figure \ref{fig:aggregation}, where each layer of the subgraphs and their edges are color coded (i.e., red is at depth of one and blue is at depth of two). The {\tt aggregate} function can now be called on $\mg$ (Line \ref{line:aggregate}) to form a single node containing an equivalent optimization problem as that on $\mg$. The OptiGraph of the node can be queried by calling the function {\tt source\_graph} (Line \ref{line:aggregate_graph}; upper right of Figure \ref{fig:aggregation}). Alternatively, we can call {\tt aggregate\_to\_depth}  and set the level to zero (a depth of one; Line \ref{line:agg_to_depth1}) and aggregate $\mg_1$ and $\mg_2$ into separate nodes (right middle of Figure \ref{fig:aggregation}). Because $\mg$ has a depth of two, we can instead go one level more and call {\tt aggregate\_to\_depth} and set the level to one  (Line \ref{line:agg_to_depth2}), forming the lower right graph in Figure \ref{fig:aggregation}.
\\

There are several important observations to make about partitioning and aggregation. The aggregation functions return {\it new} subgraphs (containing equivalent optimization problems) rather than altering in place the existing subgraphs. The {\tt aggregate} and {\tt aggregate\_to\_depth} both also return a reference map (e.g., {\tt ref\_map} on Line \ref{line:aggregate}) which maps the variables, constraints, nodes, and edges of the original graph to the new, aggregated graph. While the partitioning and aggregating functions in {\tt Plasmo.jl} give flexibility to the user, the principles behind them can also be applied by a user as they implement the model. For instance, a user can decide, when constructing the original graph, how many nodes to use to build their problem, or they can decide to construct their problem in a modular way, with each module of their problem being contained on a subgraph that is then added to a higher level subgraph (effectively partitioning the problem from the initial construction). Modeling rules can also be pre-built (from modeling libraries) and embedded in a general graph structure; for example, if we have a pre-existing model of a battery, this can be easily embedded in a planning problem for a power grid. 

\begin{figure}[!htp]
    \begin{minipage}[t]{0.9\linewidth}
        \begin{scriptsize}
        \lstset{language=Julia, breaklines = true}
        \begin{lstlisting}[label=code:aggregation, caption = {Code for aggregating the graph or subgraphs into nodes using the {\tt aggregate} and {\tt aggregate\_to\_depth} functions.}] 
# Define a function for building a graph with a specific partition
function build_partitioned_graph() |\label{line:agg_function_start}|
    g = OptiGraph()
    @optinode(g, nodes[1:4])
            
    for node in nodes
        @variable(node, 0 <= x)
        @objective(node, Min, x)
    end
    @linkconstraint(g, nodes[1][:x] + nodes[2][:x] >= 1)
    @linkconstraint(g, nodes[2][:x] + nodes[3][:x] >= 1)
    @linkconstraint(g, nodes[3][:x] + nodes[4][:x] >= 1)
    @linkconstraint(g, nodes[4][:x] + nodes[1][:x] >= 1)
            
    node_membership_vector = [1, 1, 2, 2]
    partition = Plasmo.Partition(g, node_membership_vector)
            
    apply_partition!(g, partition)
    return g
end |\label{line:agg_function_end}|

# Define an overall graph and a set of partitioned subgraphs
g = OptiGraph() |\label{line:agg_build_g}|
g1 = build_partitioned_graph() |\label{line:agg_build_g1}|
g2 = build_partitioned_graph() |\label{line:agg_build_g2}|
            
# Add the partitioned subgraphs to the overall graph g
add_subgraph!(g, g1) |\label{line:agg_add_g1}|
add_subgraph!(g, g2) |\label{line:agg_add_g2}|
            
# Add extra nodes to each graph
@optinode(g, n_g)  |\label{line:agg_nodeg}|
@optinode(g1, n_g1) |\label{line:agg_nodeg1}|
@optinode(g2, n_g2) |\label{line:agg_nodeg2}|
            
# Aggregate the overall graph into a single node
agg_node, ref_map = aggregate(g) |\label{line:aggregate}|
agg_graph = source_graph(agg_node) |\label{line:aggregate_graph}|

# Alternatively aggregate the subgraphs at different depth levels
agg_graph_layer0, ref_map_layer0 = aggregate_to_depth(g, 0) |\label{line:agg_to_depth1}|
agg_graph_layer1, ref_map_layer1 = aggregate_to_depth(g, 1) |\label{line:agg_to_depth2}|
        \end{lstlisting}
        \end{scriptsize}
    \end{minipage}
\end{figure}

\section{{\tt PlasmoBenders.jl} Functionality}

Below, we give further details on {\tt PlasmoBenders.jl} functionality. Additional information can also be found in the package documentation. 

\subsubsection{Regularization}
The regularization scheme implemented in {\tt PlasmoBenders.jl} follows that outlined by \cite{pecci2024regularized}. This scheme is applied to two stage BD problems and is not supported for problems with more than two stages. The purpose of the regularization is to choose feasible solutions from the root subgraph that result in better solutions (or better cuts) when those solutions are fixed in the second stage sub-problems. Often, a challenge with BD is that the cuts generated from the root-problem's solutions are poor, and it can require many iterations (and therefore many cuts) before reaching high-quality solutions. This can be in part because using the ``optimal'' solutions of the root-problem may be extreme choices. Regularization seeks to mitigate this by choosing solutions that are on the interior of the feasible region and which can potentially generate stronger cuts after solving the second stage sub-problems. 

The regularization scheme is applied to the root-problem as follows. The root-problem has the form 
\begin{subequations}\label{eq:NBD_root}
    \begin{align}
        \min &\; c_r^\top x_r + \theta_r\\
        \textrm{s.t.} &\; A_r x_r \le b_r \\
        &\; \theta_r \ge \{ \textrm{cuts} \}. \label{eq:NBD_root_cuts}
    \end{align}
\end{subequations}
\noindent Note that this is the same form as the root graph, so the regularization scheme is easily extended to the graph-based case. Since we want solutions that are on the interior of the problem, we change the objective but add a constraint that the original objective must still be ``near'' the true objective. This is done by updating \eqref{eq:NBD_root} to be 
\begin{subequations}\label{eq:NBD_root_reg}
    \begin{align}
        \min &\; \Psi(x_r)\\
        \textrm{s.t.} &\; A_r x_r \le b_r \\
        &\; \theta_r \ge \{ \textrm{cuts} \} \\
        &\; c_r^\top x_r + \theta_r \le LB^k + \alpha (UB^k - LB^k) \label{eq:NBD_reg_con}
    \end{align}
\end{subequations}
where $\Psi$ is a criterion for choosing a feasible point, $LB^k$ and $UB^k$ are the best lower and upper bounds, respectively, up to iteration $k$, and $\alpha \in (0, 1]$ is the regularization parameter that dictates how ``close'' the chosen interior point must be to the true optimal solution. $\Psi$ can be a variety of functions \cite{pecci2024regularized}, but for {\tt PlasmoBenders.jl}, we set $\Psi(x_r) = 0$; in other words, any feasible solution to \eqref{eq:NBD_root_reg} can be returned as a solution. The additional constraint \eqref{eq:NBD_reg_con} ensures that any decision is at least less than the upper bound (note that the lower bound is in fact the optimal solution of \eqref{eq:NBD_root}). Solving \eqref{eq:NBD_root_reg} provides the feasible solutions that are then used for solving the second-stage sub-problems. Since \eqref{eq:NBD_root_reg} requires $LB^k$, \eqref{eq:NBD_root} also still must be solved at each iteration.

{\tt PlasmoBenders.jl} can implement the regularization procedure outlined by \eqref{eq:NBD_root_reg} and provides options to the user for this scheme. The BD and graph-based BD (gBD) algorithm outlined in the manuscript are not altered by the regularization scheme, but the statement in the algorithm, ``Determine root-problem decisions $\bar{x}^k_r$'' (or $\bar{x}^k_{\mathcal{G}_r}$) is where the regularization scheme is applied. In other words, \eqref{eq:NBD_root} still must be solved to get $LB^k$, but the solutions fixed in the second stage come from solving \eqref{eq:NBD_root_reg}. This regularization scheme can be applied by setting the key word argument {\tt regularize} to {\tt true} in the {\tt BendersOptimizer} constructor of {\tt PlasmoBenders.jl}. The parameter $\alpha$ can also be set via the {\tt regularize\_param} argument of the same function. Lastly, we note that the user must set the solver on the subgraphs used in {\tt PlasmoBenders.jl}, and the choice of solver and solver attributes can influence the performance of the algorithm. For instance, when using Gurobi, we found that the {\tt Method} option is best set to 2 so that the solver uses a barrier algorithm to get an interior, feasible point; using the default {\tt Method} instead used the original optimal solution to \eqref{eq:NBD_root} since it is still an optimal solution. 

\section{Case Study 1}
For the case study of the tri-level energy market operation, we offer a few additional details on the problem and solution.

\subsection{Mathematical Formulation}
To better present the mathematical definition, we define the following sets, notation, and variables. We use the sets $\Gamma$ for the set of all generators, $\Gamma_r$ for the set of renewable generators, $\Gamma_c$ for the set of conentional generators with $\Gamma_c^d \subseteq \Gamma_c$ as the set of generators committed by the DA-UC layer and $\Gamma_c^s \subseteq \Gamma_c$ as the subset of generators committed in the ST-UC layer (note that $\Gamma_c^d \cap \Gamma_c^d = \emptyset$). We also define $\Gamma_c^h = \Gamma_c^d \cup \Gamma_c^s$. We also use $*$ and $**$ to define variables, sets, or functions corresponding to a layer; here, the symbols $*$ or $**$ are exchanged for $d$, $s$, or $h$ to denote the DA-UC, ST-UC, or HA-ED layers, respectively. As each sub-problem considers different sets of times, we use $\mathcal{T}^*_i$ for the set of times (in hours) of the $i$th sub-problem of the $*$ layer. We also define the sets $\bar{\mathcal{T}}^*_i$ as the set of times without the first time point of the sub-problem. We define $\Delta^*$ as the time step for the sub-problem in hours($\Delta^d = 1$, $\Delta^s = 0.25$, and $\Delta^h = 0.25$). Note that we set $\Delta^s = \Delta^h$, and this has a small influence on the formulation of the problem presented below. The decision variables are given by:
\begin{align*}
    \begin{split}
    \bx^d_{i}  =&\; \big(x_{g, t}, s_{g, t}, z_{g, t}\big)_{\forall g \in \Gamma^d_c, t \in \mathcal{T}^d_i} \\
    \by^d_{i} =&\;  \Big((G^{+, d}_{g, t}, G^{-, d}_{g, t})_{\forall{g \in \Gamma^d_c \cup \Gamma_r}}, (F^d_{j, k, t})_{\forall (j, k)\in \mathcal{L}}  (D^{d}_{j, t}, \theta^d_{j, t})_{\forall j \in \mathcal{B}}\Big)_{\forall t \in \mathcal{T}^d_i}\\
    \color{red} \bx^s_{i} =&\; \color{red}  \big(x_{g, t}, s_{g, t}, z_{g, t}\big)_{\forall g \in \Gamma^s_c, t \in \mathcal{T}^s_i} \\
    \color{red} \by^s_{i} =&\; \color{red} \Big((G^{+, s}_{g, t}, G^{-, s}_{g, t})_{\forall{g \in \Gamma_c^h \cup \Gamma_r}}, (F^s_{j, k, t})_{\forall (j, k)\in \mathcal{L}}, (D^{s}_{j, t}, \theta^s_{j, t})_{\forall j \in \mathcal{B}}\Big)_{\forall t \in \mathcal{T}^s_i}\\
    \color{blue} \by^h_{i} =&\; \color{blue}  \Big((G^{+, h}_{g, t}, G^{-, h}_{g, t})_{\forall{g \in \Gamma_c^h \cup \Gamma_r}}, (F^h_{j, k, t})_{\forall (j, k)\in \mathcal{L}} (D^{h}_{j, t}, \theta^h_{j, t})_{\forall j \in \mathcal{B}}\Big)_{\forall t \in \mathcal{T}^s_i}
    \end{split}
\end{align*}
Symbols $\bx^d_i$ and $\by^d_i$ are binary and continuous decision variables for the $i$th sub-problem of DA-UC, $\color{red} \bx^s_i$ and $\color{red} \by^s_i$ are decision variables for the $i$th sub-problem of ST-UC, and $\color{blue} \by^h_i$ are decision variables for the $i$th sub-problem o HA-ED.  Symbols $x_{g, t}$, $s_{g, t}$, and $z_{g,t}$ are binary variables indicating for time $t$ whether generator $g$ is on/off, was turned on, or was turned off, respectively. $G^{+, *}_{g, t}$ is the power of generator $g$ consumed by the grid at time $t$ and $G^{-, *}_{g, t}$ is the power overgenerated (for conventional generators) or curtailed (for renewable generators) from $g$ at time $t$. $F^*_{j, k, t}$ is the power flow of transmission line $(j, k)$ from bus $j$ to bus $k$ during time $t$. $D^*_{j, t}$ is the amount of load shed at bus $j$ for time $t$, and $\theta^*_{j, t}$ is the bus angle for bus $j$ at time $t$. We also use black, red, and blue color to denote variables for DA-UC, ST-UC, and HA-ED, respectively.

The objective functions in the different layers are comprised of a UC part, $f_u$, and an ED part, $f_e$:
\begin{align}\label{eq:f_u}
    f_{u}(\bx^*_i) := \sum_{t \in \mathcal{T}^*_i} \sum_{g \in \Gamma^*_c} \phi_g^s s_{g, t} + \phi_g^f x_{g, t} \Delta^*
\end{align}
\begin{align}\label{eq:f_e}
    \begin{split}
        f_{e}(\by^*_i) := &\; \sum_{t \in \mathcal{T}^*_i} \Big( \sum_{g \in \Gamma^*_c} \phi^v_g G^{+,*}_{g, t} + \sum_{j \in \mathcal{B}} \Big(\phi^u_j D^{*}_{j, t} + \sum_{g \in \Gamma_j \cap \Gamma^*_c}\phi^o_g G^{-, *}_{g, t} + \sum_{g \in \Gamma_j \cap \Gamma_r} \phi^c_g G^{-, *}_{g, t}   \Big) \Big)
    \end{split}
\end{align}
The function $f_u$ accounts for the startup cost $\phi^s_g$ and the no-load cost $\phi^f_g$ for generator $g$ at time $t$. The no-load cost is multiplied by $\Delta^*$ since the DA-UC and ST-UC levels have different time resolutions. The function $f_e$ accounts for the variable cost, $\phi^v_g$, of the energy consumed by the grid, the cost of overgeneration, $\phi^o_g$, and the cost of curtailment, $\phi^c_g$, for generator $g$ at time $t$. It also accounts for the cost of unmet demand $\phi^u_j$. We next define constraints for the layers: 
\begin{align}\label{eq:c_d}
    \begin{split}
        c_d(\by_i^*) := &\; \Big( \sum_{j \in \mathcal{B}:(j, k) \in \mathcal{L}} F^*_{j, k, t} - \sum_{j \in \mathcal{B}:(k, j) \in \mathcal{L}} F^*_{k, j, t} + \sum_{g \in \Gamma_k} G^{+, *}_{g, t} + D^{*}_{k, t} - \hat{D}^*_{k, t} - \hat{R}^*_{k, t} \Big)_{k \in \mathcal{B}, t \in \mathcal{T}_i^*}
    \end{split}
\end{align}
This requires that, at each time point, the flows coming into the bus, the power consumed by the grid, and the amount of unmet demand is equal to the demand, $\hat{D}^*_{k, t}$, and the reserve requirements, $\hat{R}^*_{k, t}$ for $k \in \mathcal{B}$. The power flow equations use a DC approximation: 
\begin{align}\label{eq:c_f}
    c_f(\by^*_{i}) := \left( F^*_{j, k, t} - B_{j, k} (\theta^*_{j, t} - \theta^*_{k, t}) \right)_{(j, k) \in \mathcal{L}, t \in \mathcal{T}^*_i}
\end{align}
The renewable resources are also restricted to a specific value; this is enforced by constraint $c_r$, where $\hat{G}^*_{g, t}$ is the amount of power produced by renewable generator $g$ at time $t$: 
\begin{align}\label{eq:c_r}
    c_r(\by^*_{i}) := \left( G^{+, *}_{g, t} + G^{-, *}_{g, t} - \hat{G}^*_{g, t} \right)_{g \in \Gamma_r, t \in \mathcal{T}^*_i}
\end{align}
The ramp-up and ramp-down constraints ($c_{ru}$ and $c_{rd}$) have different forms because DA-UC and ST-UC/HA-ED have different time resolutions. Ramping constraints containing start-up and shut-down constraints are:
\begin{align}\label{eq:c_su}
    \begin{split}
    c_{su}(\by^*_{i}, \bx^{**}_j) := &\;\Big( G^{+, *}_{g, t} + G^{-, *}_{g, t} - G^{+, *}_{g, t - \Delta^*} - G^{-, *}_{g, t - \Delta^*} - (\overline{S}_g - \overline{R}_g \Delta^{**} - \underline{C}_g) s^{**}_{g, t + \Delta^{**}} - \\
    &\; \qquad \qquad (\overline{R}_g \Delta^{**} + \underline{C}_g) x^{**}_{g, t + \Delta^{**}} + \underline{C}_{g} x^{**}_{g, t} \Big)_{g \in \Gamma_c^{**}, t \in \Theta}
    \end{split}
\end{align}
\vspace{-0.3in}
\begin{align}\label{eq:c_sd}
    \begin{split}
    c_{sd}(\by^*_{i}, \bx^{**}_j) :=&\; \Big( G^{+, *}_{g, t - \Delta^*} + G^{-, *}_{g, t - \Delta^*} - G^{+, *}_{g, t} - G^{-, *}_{g, t} - (\underline{S}_g - \underline{R}_g \Delta^* - \underline{C}_g) z^{**}_{g, t + \Delta^{**}} -\\
    &\; \qquad \qquad (\underline{R}_g \Delta^{**} + \underline{C}_g) x^{**}_{g, t} +  \underline{C}_{g} x^{**}_{g, t + \Delta^{**}} \Big)_{g \in \Gamma_c^{**}, t \in \Theta}
    \end{split}
\end{align}
Here, $\overline{S}_g$ and $\underline{S}_g$ are the startup/shutdown limits for $g$, $\overline{R}_g$ and $\underline{R}_g$ are the ramp-up and ramp-down limits for $g$ as a function of time, and $\underline{C}_g$ is the minimum capacity of $g$. The set $\Theta$ is defined in this context as $\Theta = \{ t_i: t_i - \Delta^{*}= t_j - \Delta^{**}, \forall t_i \in \bar{\mathcal{T}}^*_i, t_j \in \mathcal{T}_j^{**} \}$. The model uses constraints $c_{su, 0}(\by^*_i, \by^h_k, \bx_j^{**}, \bx_{j - 1}^{**})$ and $c_{sd, 0}(\by^*_i,\by^h_k, \bx_j^{**}, \bx_{j - 1}^{**})$ to link the first time point of the $i$th sub-problem with the solutions of the previous sub-problems (this introduces complex time coupling).  We also define operating regions:
\begin{align}
    \mathcal{X}_i^* &:=  \left\{ \bx^*_i | (x_{g,t}, s_{g, t}, z_{g, t}) \in \{0, 1\}^3, \quad \substack{\forall g \in \Gamma_c^*,\\ t \in \mathcal{T}_i^*} \right\}\label{eq:x}\\    
    \mathcal{Y}^* &:= \left\{ \begin{array}{c|l} & \underline{\theta}_j \le \theta_{j, t}^d \le \overline{\theta}_j, \quad \forall j \in \mathcal{B}, t \in \mathcal{T}_i^d\\
        & \underline{F}_{j, k} \le F^*_{j, k, t} \le \overline{F}_{j, k}, \quad \forall (j, k) \in \mathcal{L}, t \in \mathcal{T}_i^d\\
        \by^*_i  & (G^{+, *}_{g, t}, G^{-, *}_{g, t}) \in \mathbb{R}_{+}^2, \quad \forall g \in \Gamma_c^d \cup \Gamma_r, t \in \mathcal{T}_i^d \\
        & F^*_{j, k, t} \in \mathbb{R}, \quad \forall (j, k) \in \mathcal{L}, t \in \mathcal{T}^d_i\\
        & \theta^*_{j, t} \in \mathbb{R},  D^*_{j, t} \in \mathbb{R}_+, \quad \forall j \in \mathcal{B}, t \in \mathcal{T}^*_i
        \end{array}\right\}\label{eq:y}
\end{align}
The $i$th DA-UC sub-problem is given by \eqref{eq:dauc}
\begin{subequations}\label{eq:dauc}
\begin{align}
        \min &\; f_{u}(\bx^d_i) + f_{e}(\by_i^d) \Delta^d \\
        \textrm{s.t.} &\; c_d(\by^d_i) = 0, \quad c_f(\by^d_i) = 0, \quad c_r(\by^d_i) = 0\\
        &\; c_{su}(\by^d_i, \bx^d_i) \le 0, \quad c_{su, 0}(\by^d_i, \textcolor{blue}{\by^h_k}, \bx^d_i, \bx^d_{i -1}) \le 0 \label{eq:dauc_link1}\\
        &\;  c_{sd}(\by^d_i, \bx^d_i) \le 0, \quad c_{sd, 0}(\by^d_i, \textcolor{blue}{\by^h_k}, \bx^d_i, \bx^d_{i -1}) \le 0 \label{eq:dauc_link2}\\
        &\; x_{g, t} - x_{g, t - 1} = s_{g, t} - z_{g, t} \quad \forall g \in \Gamma_c^d, t \in \mathcal{T}^{d}_i\label{eq:dauc_onoff}\\
        &\; \sum_{j \in {\mathcal{U}^d_{g, t}}} s_{g, j} \le x_{g, t} \quad \forall g \in \Gamma_c^d, t \in \mathcal{T}^d_i \label{eq:dauc_uptime}\\
        &\; \sum_{j \in \mathcal{D}^d_{g, t}}^t z_{g, j} \le 1 - x_{g, t} \quad \forall g \in \Gamma_c^d, t \in \mathcal{T}_i^d \label{eq:dauc_downtime}\\
        &\; \underline{C}_g x_{g, t} \le G^{+, d}_{g, t} + G^{-, d}_{g, t} \le \overline{C}_g x_{g, t}, \quad \substack{\forall g \in \Gamma_c^d, \\t \in \mathcal{T}_i^d} \label{eq:dauc_cap}
\end{align}
\end{subequations}
The $i$th sub-problem of ST-UC is given by \eqref{eq:stuc}. This has a similar structure as DA-UC; however, there are now DA-UC variables that are incorporated into this lower layer solution through \eqref{eq:stuc_su}, \eqref{eq:stuc_sd}, \eqref{eq:stuc_cap}, and \eqref{eq:stuc_link}. The last constraint ensures that the generation amounts for the DA-UC generators in ST-UC are within a certain bound ($\epsilon^s_g$) of the DA-UC sub-problem solutions. This helps avoid myopic solutions, since the ST-UC time horizon is shorter than that of DA-UC. 
\begin{subequations}\label{eq:stuc}
\begin{align}
        \color{red} \min &\; \color{red} f_{u}(\bx^s_i) + f_{e}(\by_i^s) \Delta^s \\
        \color{red} \textrm{s.t.} &\; \color{red} c_d(\by^s_i) = 0, \quad \color{red} c_f(\by^s_i) = 0, \quad \color{red} c_r(\by^s_i) = 0\\
        &\; \color{red} c_{su}(\by^s_i, \textcolor{black}{\bx^d_j}) \le 0, \quad c_{su, 0}(\by^s_i, \textcolor{blue}{\by^h_k}, \textcolor{black}{\bx^d_j}, \textcolor{black}{\bx^d_{j -1}}) \le 0 \label{eq:stuc_su}\\
        &\;  \color{red} c_{sd}(\by^s_i, \textcolor{black}{\bx^d_j}) \le 0, \quad c_{sd, 0}(\by^s_i, \textcolor{blue}{\by^h_k}, \textcolor{black}{\bx^d_j}, \textcolor{black}{\bx^d_{j -1}}) \le 0 \label{eq:stuc_sd} \\
        &\; \color{red} c_{su}(\by^s_i, \bx^s_i) \le 0, \quad c_{su, 0}(\by^s_i, \textcolor{blue}{\by^h_k}, \bx^s_i, \bx^s_{i - 1}) \le 0\\
        &\; \color{red} c_{sd}(\by^s_i, \bx^s_i) \le 0, \quad c_{sd, 0}(\by^s_i, \textcolor{blue}{\by^h_k}, \bx^s_i, \bx^s_{i - 1}) \le 0 \label{eq:stuc_link1}\\
        &\; \color{red} c_{ru}(\by^s_i) \le 0, \quad \quad \color{red} c_{rd}(\by^s_i) \le 0 \\
        &\; \color{red} x_{g, t} - x_{g, t - 1} = s_{g, t} - z_{g, t} \quad \forall g \in \Gamma_c^s, t \in \mathcal{T}^{s}_i\label{eq:stuc_binary}\\
        &\; \color{red} \sum_{j \in \mathcal{U}^s_{g, t}} s_{g, j} \le x_{g, t} \quad \forall g \in \Gamma_c^s, t \in \mathcal{T}^s_i \label{eq:stuc_ut}\\
        &\; \color{red} \sum_{j \in \mathcal{D}^s_{g, t}} z_{g, j} \le 1 - x_{g, t} \quad \forall g \in \Gamma_c^s, t \in \mathcal{T}_i^s \label{eq:stuc_dt} \\
        &\; \color{red} \underline{C}_g \textcolor{black}{x_{g, \lceil t \rceil}} \le G^{+, s}_{g, t} + G^{-, s}_{g, t} \le \overline{C}_g\textcolor{black}{x_{g, \lceil t \rceil }} \quad \substack{ \forall g \in \textcolor{black}{\Gamma_c^d}, \\ t \in \mathcal{T}_i^s} \label{eq:stuc_cap} \\
        &\; \color{red} \underline{C}_g x_{g, t} \le G^{+, s}_{g, t} + G^{-, s}_{g, t} \le \overline{C}_g x_{g, t} \quad \substack{\forall g \in \Gamma_c^s,\\ t \in \mathcal{T}_i^s}\\
        &\; \color{red} |G^{+, s}_{g, t} + G^{-, s}_{g, t} - \textcolor{black}{G^{+, d}_{g, \lceil t \rceil}} - \textcolor{black}{G^{-, d}_{g, \lceil t \rceil}}| \le \epsilon^s_g, \quad \substack{\forall g \in \textcolor{black}{\Gamma_c^d}, \\t \in \mathcal{T}_i^s} \label{eq:stuc_link}
\end{align}
\end{subequations}
The $i$th sub-problem of HA-ED is given by \eqref{eq:haed}. This formulation is similar to ST-UC but without the binary variables and their accompanying constraints and objective function. In addition, there are now links between {\it both} the ST-UC and DA-UC layers. 
\begin{subequations}\label{eq:haed}
\begin{align}
        \color{blue} \min &\; \color{blue} f_{e}(\by_i^h) \Delta^h \\
        \color{blue} \textrm{s.t.} &\; \color{blue} c_d(\by^h_i) = 0, \quad c_f(\by^h_i) = 0, \quad  c_r(\by^h_i) = 0\\
        &\; \color{blue} c_{su}(\by^h_i, \textcolor{black}{\bx^d_j}) \le 0, \quad c_{su, 0}(\by^h_i, \by^h_{i - 1}, \textcolor{black}{\bx^d_j}, \textcolor{black}{\bx^d_{j - 1}}) \label{eq:haed_link1}\\
        &\; \color{blue} c_{sd}(\by^h_i, \textcolor{black}{\bx^d_j}) \le 0, \quad c_{sd, 0}(\by^h_i, \by^h_{i - 1}, \textcolor{black}{\bx^d_j}, \textcolor{black}{\bx^d_{j - 1}})\\
        &\; \color{blue} c_{su}(\by^h_i, \textcolor{red}{\bx^s_i}) \le 0, \quad c_{su, 0}(\by^h_i, \by^h_{i - 1}, \textcolor{red}{\bx^d_j}, \textcolor{red}{\bx^d_{j - 1}})\label{eq:haed_stuc_su} \\
        &\; \color{blue} c_{sd}(\by^h_i, \textcolor{red}{\bx^s_i}) \le 0, \quad c_{sd, 0}(\by^h_i, \by^h_{i - 1}, \textcolor{red}{\bx^d_j}, \textcolor{red}{\bx^d_{j - 1}})\label{eq:haed_stuc_sd} \\
        &\; \color{blue} c_{ru}(\by^h_i) \le 0, \quad \quad c_{rd}(\by^h_i) \le 0 \\
        &\; \color{blue} \underline{C}_g \textcolor{black}{x_{g, \lceil t \rceil}} \le G^{+, h}_{g, t} + G^{-, h}_{g, t} \le \overline{C}_g \textcolor{black}{x_{g, \lceil t \rceil}} \quad \substack{\forall g \in \textcolor{black}{\Gamma_c^d}, \\t \in \mathcal{T}_i^h}\\
        &\; \color{blue} \underline{C}_g \textcolor{red}{x_{g, t}} \le G^{+, h}_{g, t} + G^{-, h}_{g, t} \le \overline{C}_g \textcolor{red}{x_{g, t}} \quad \substack{\forall g \in \textcolor{red}{\Gamma_c^s}, \\t \in \mathcal{T}_i^h}\label{eq:haed_stuc_cap} \\
        &\; \color{blue} |G^{+, h}_{g, t} + G^{-, h}_{g, t} - \textcolor{black}{G^{+, d}_{g, \lceil t \rceil}} - \textcolor{black}{G^{-, d}_{g, \lceil t \rceil}}| \le \epsilon_g, \quad \substack{\forall g \in \textcolor{black}{\Gamma_c^d}, \\t \in \mathcal{T}_i^h } \label{eq:haed_dauclinks}\\
        &\; \color{blue} |G^{+, h}_{g, t} + G^{-, h}_{g, t} - \textcolor{red}{G^{+, s}_{g, \lceil t \rceil}} - \textcolor{red}{G^{-, s}_{g, t}}| \le \epsilon_g, \quad \forall \substack{g \in \textcolor{red}{\Gamma_c^s}, \\t \in \mathcal{T}_i^h} \label{eq:haed_stuc_links}.
\end{align}
\end{subequations}

\subsubsection{Graph Representation}

\begin{figure}[!htp]
    \centering
    \vspace{-0.2in}
    \includegraphics[width=0.35\textwidth]{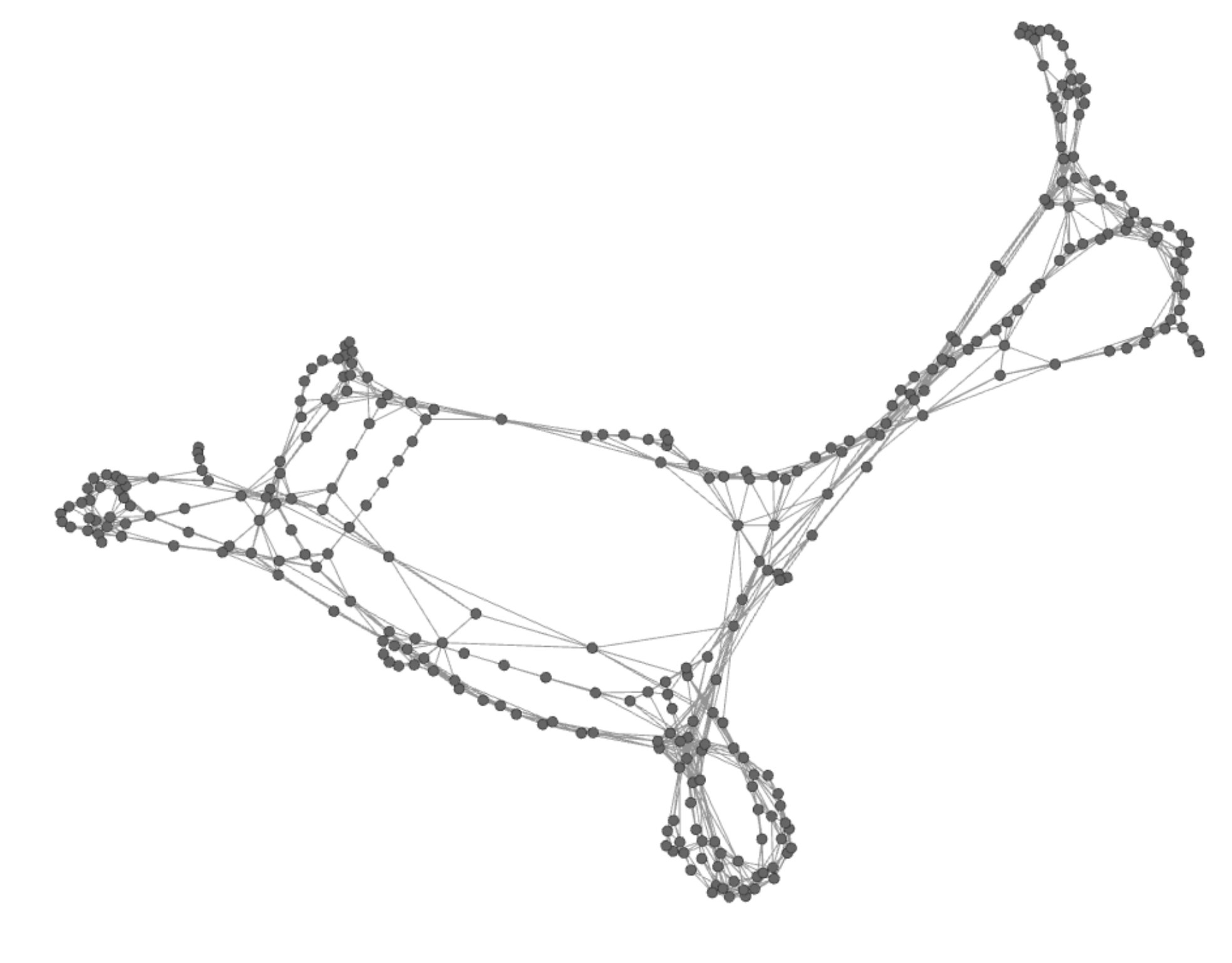}
    \vspace{-0.1in}
    \caption{A graph visualization of the 118-bus system, where nodes correspond to buses and to transmission lines.}
    \label{fig:118_bus_graph}    
\end{figure}

\begin{figure}[!htp]
    \centering
    \vspace{-0.2in}
    \includegraphics[width = 0.4\textwidth]{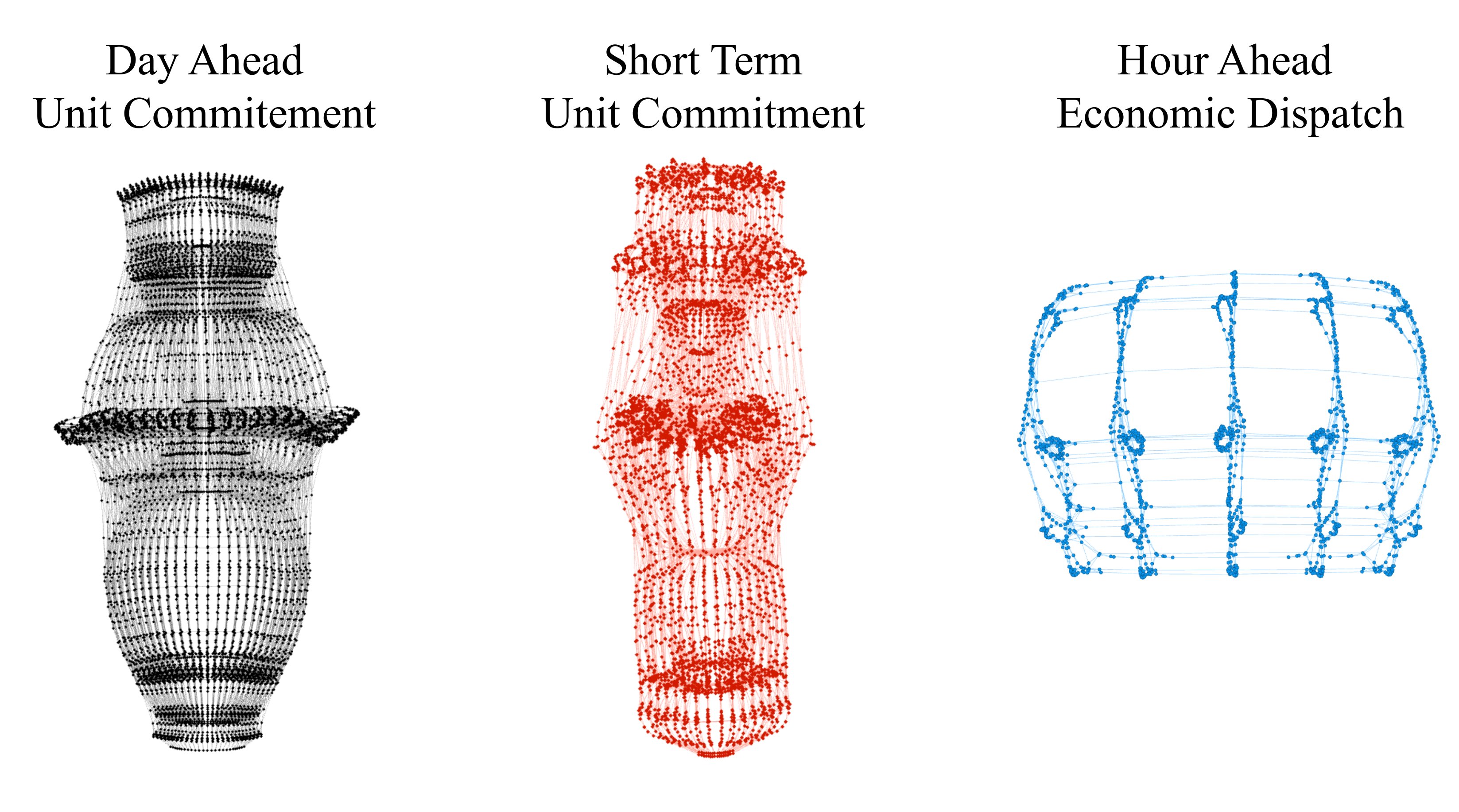}
        \vspace{-0.1in}
    \caption{Representation of a single sub-problem subgraph for the DA-UC, ST-UC and HA-ED sub-problems for the tri-level hierarchical problem of \cite{atakan2022towards}.}
    \label{fig:subproblems}
\end{figure}

We next outline how we represent the tri-level hierarchical architecture using the OptiGraph abstraction in {\tt Plasmo.jl}. The overall problem for a single day is represented by a graph $\mg(\{\mg_i\}_{i \in P}, \emptyset, \me_\mg)$, which contains a subgraph for each DA-UC, ST-UC, and HA-ED sub-problem, creating 105 total subgraphs (i.e., $P = \{DA1, ST1, ST2,..., ST8, HA1, HA2, ...,  HA96\}$, where an index encodes both number and layer). Within each subgraph, each time point $t$ is represented as an additional subgraph, and nodes are placed on these time-point-subgraphs for each bus and each transmission line in the 118-bus system \cite{pena2017extended}. The nodes corresponding to buses contain the variables $D^*_{k, t}$, $\theta^*_{k, t}$, $G^{+, *}_{g, t}$, and $G^{-, *}_{g, t}, \forall g \in \Gamma_k \cap \{ \Gamma_c^* \cup \Gamma_r \}$. In the case of DA-UC and ST-UC problems, the bus nodes also contain $x_{g, t}, s_{g, t}$, and  $z_{g, t} \: \forall g \in \Gamma_k \cap \Gamma_c^*$ and any constraints for these variables. The nodes corresponding to transmission lines contain variables $F^*_{i, j, t}, \forall (i, j) \in \mathcal{L}$. For each node representing line $(i, j) \in \mathcal{L}$, edges (linking constraints) are also placed connecting to bus $i$ and to bus $j$. The resulting subgraph for any time point is shown in Figure \ref{fig:118_bus_graph}. Each DA-UC subgraph had 25 time point subgraphs (i.e., 25 replicates of the network shown in Figure \ref{fig:118_bus_graph}; note that we used 25 rather than 24 time points in this problem to simplify modeling) each representing one hour, each ST-UC subgraph had 16 time point subgraphs with each representing 15 minutes (four hours total), and each HA-ED subgraph had 5 time point subgraphs with each representing 15 minutes (75 minutes total). Linking constraints were also placed between time point subgraphs where applicable, such as for $c_{ru}$, $c_{rd}$, \eqref{eq:dauc_onoff} - \eqref{eq:dauc_downtime}, or \eqref{eq:stuc_binary} - \eqref{eq:stuc_dt}. An example of each sub-problem graph is shown in Figure \ref{fig:subproblems}

\begin{figure}[!ht]
    \centering
        \vspace{-0.2in}
    \includegraphics[width = 0.3\textwidth]{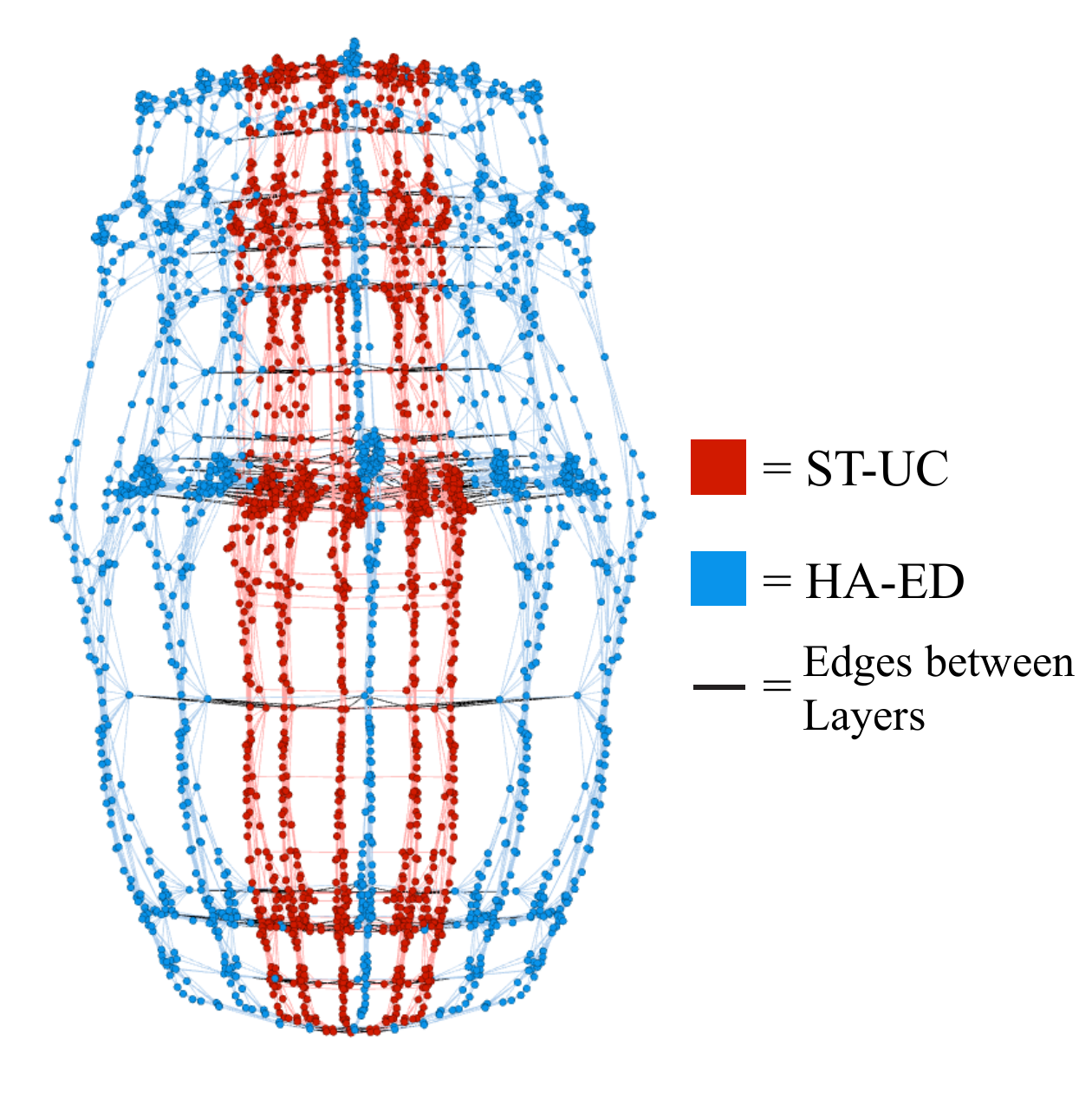}
        \vspace{-0.1in}
    \caption{An example of the hierarchical linking between sub-problems. Five time points from an ST-UC sub-problem are shown linked to a full HA-ED sub-problem, with linking constraints highlighted in black.}
    \label{fig:linking_constraints}
\end{figure}

\begin{figure}[!ht]
    \centering
        \vspace{-0.2in}
    \includegraphics[width = 0.9\textwidth]{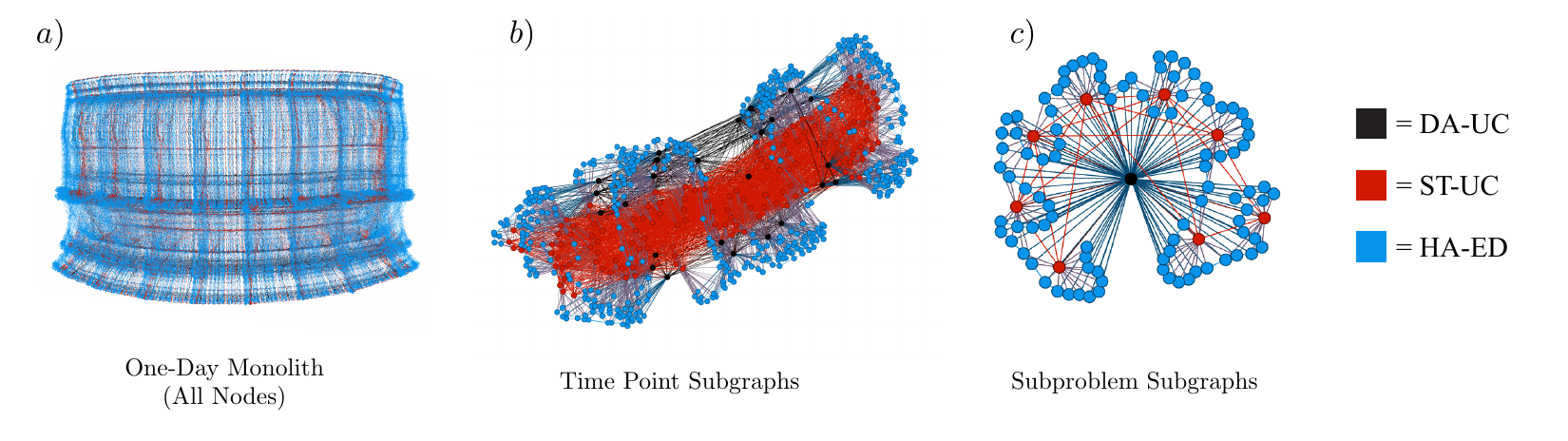}
        \vspace{-0.1in}
    \caption{Representation of (a) the 1-day monolithic graph highlighted by the sub-problems (DA-UC, ST-UC, HA-ED), and equivalent representations where (b) individual time point subgraphs are represented as nodes and where (c) sub-problem subgraphs are represented as nodes.}
    \label{fig:CS1_subgraph_representations}
\end{figure}

The hierarchical structure in this problem is captured by the different subgraphs and by edges, $\me_\mg$, placed between subgraphs. These edges contain the linking constraints in \eqref{eq:dauc}, \eqref{eq:stuc}, and \eqref{eq:haed} which can be identified by the constraints containing variables from different sub-problems (e.g., \eqref{eq:dauc_link1}-\eqref{eq:dauc_link2}, \eqref{eq:stuc_su}-\eqref{eq:stuc_link1}, \eqref{eq:haed_link1}-\eqref{eq:haed_stuc_sd}). An example of the linking between layers is shown in Figure \ref{fig:linking_constraints}, where five time points from an ST-UC sub-problem are shown with links to a full HA-ED sub-problem, where links are highlighted in black. The overall structure can be further visualized by representing nodes and subgraphs. The overall graph $\mg$ with all of its nodes is shown in Figure \ref{fig:CS1_subgraph_representations}a, where nodes are colored by the sub-problem type to which they belong (DA-UC, ST-UC, and HA-ED). Equivalently, we can visualize each time point subgraph as a single node (i.e., aggregating these subgraphs into nodes) as shown in Figure \ref{fig:CS1_subgraph_representations}b, resulting in a graph with fewer total nodes. Finally, we can also represent each sub-problem subgraph (i.e., $\msg(\mg)$) as a single node as shown in Figure \ref{fig:CS1_subgraph_representations}c. These visualizations highlight the hierarchy that exists in the problem between layers, and they show how the subgraphs form partitions in the problem. For instance, Figures \ref{fig:CS1_subgraph_representations}a and \ref{fig:CS1_subgraph_representations}c are valid representations of the same system, but Figure \ref{fig:CS1_subgraph_representations}c helps visualize the structure in a more intuitive way. 

\begin{figure}[!ht]
    \centering
            \vspace{-0.1in}
    \includegraphics[width=0.5\textwidth]{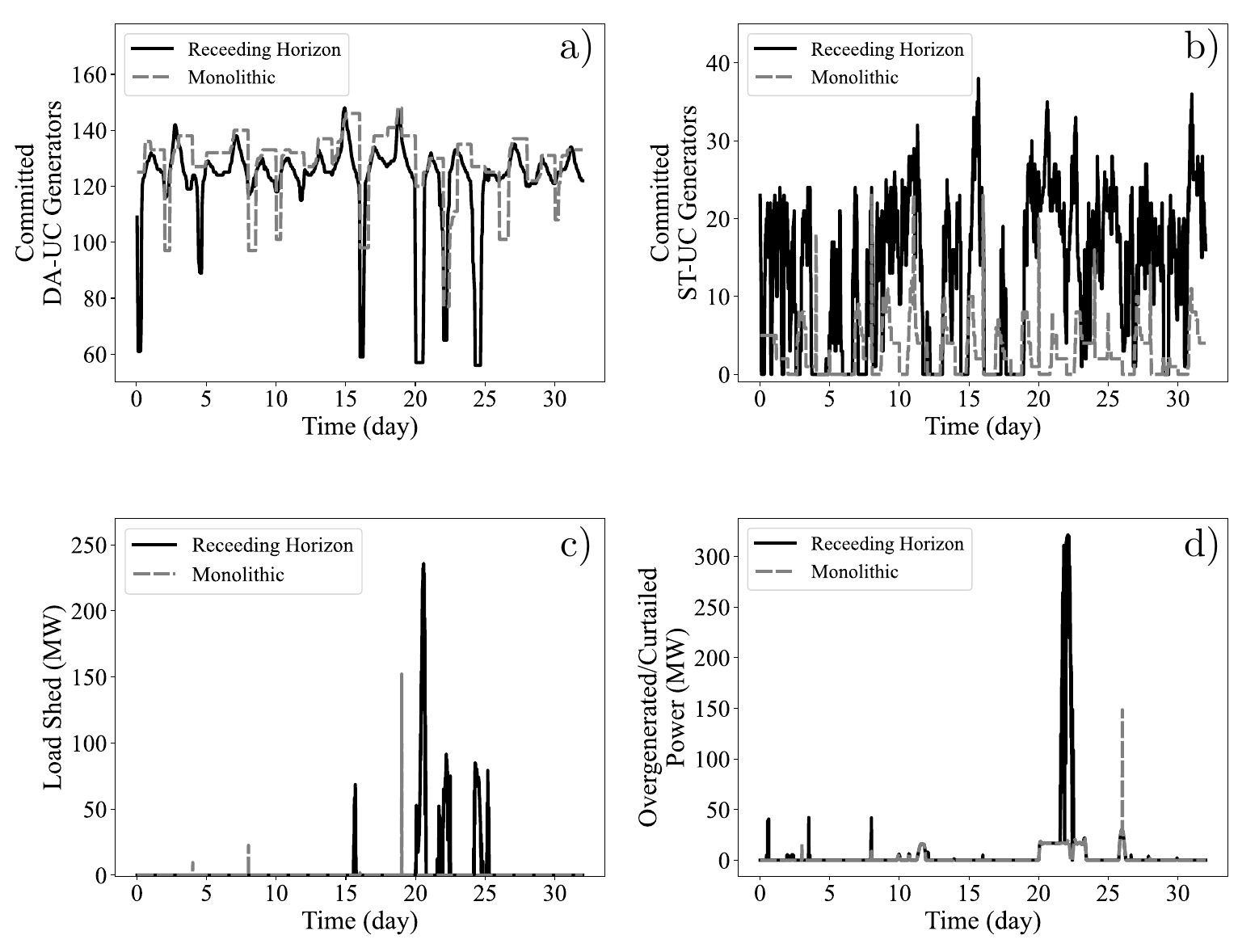}
        \vspace{-0.1in}
    \caption{Results of the tri-level problem using a "receding horizon" type approach and a "monolothic" approach. a) number of committed DA-UC generators over time, b) number of committed ST-UC generators over time, c) the load shedding, and d) the overgenerated or curtailed power.}
    \label{fig:results}    
\end{figure}

\subsection{Data and Solution Details}
\begin{figure}[!htp]
    \centering
    \includegraphics[width=0.5\linewidth]{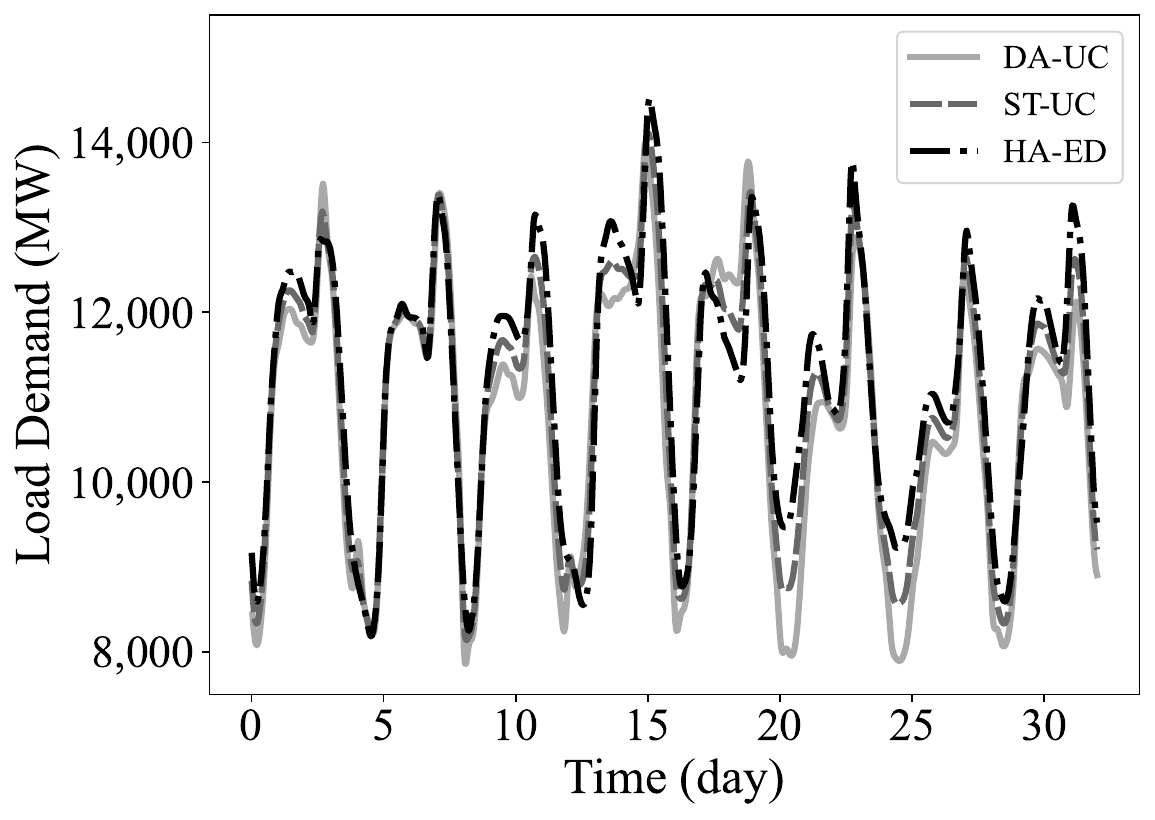}
    \caption{Load demand used in the DA-UC, ST-UC, and HA-ED sub-problems for the 32-day horizon.}
    \label{fig:demands}
\end{figure}

First, the set of demands can be seen in Figure \ref{fig:demands} for each level of the model. The demand data used in \cite{atakan2022towards} only included hourly demands for the day-ahead unit commitment (DA-UC) and hour-ahead economic dispatch (HA-ED) layers. For the short-term unit commitment (ST-UC) layer, we used the average of the other two layers. In addition, to get the data on a fifteen minute resolution for the HA-ED and ST-UC layers, we interpolated between hourly data points. 

Reserve requirements were also included in the model. We used the factors suggested by \cite{atakan2022towards} corresponding to the ``low reserve requirements'' scenario, which were set to 10\% of the demands for the UC problems and 2.5\% of the demand for the HA-ED problem. 

We also make some additional notes about solving the monolithic and receding horizon problems. First, the monolithic problem was solved with a MIP gap of 1\%. There was an additional time limit set of 2 hours for each day's simulation, but this time limit was never reached. For the receding horizon problems, we used a MIP gap of 0.3\% or less, though we did enforce a time limit of 15 minutes for the DA-UC problems and 5 minutes for the ST-UC problems. This time limit was reached for several problems. We also noticed that there were some numerical challenges in the modeling where a sub-problem in the ST-UC or HA-ED layers would return a termination status of infeasible. To avoid this problem, we added the option of using slack variables in the receding horizon approach, such that many linking constraint between time points (e.g., ramping constraints or capacity constraints) can be relaxed at the cost of having heavily penalization in the objective function. There were four infeasible problems in the HA-ED layer on day 25, but the slack variables had nonzero values in two different ramping constraints, and had very small values (e.g., < 1e-5) which is why we suspect numerical issues with the data and with solutions being passed between subproblems.

Finally, we note that an area of future work could be better facilitating passing solutions between subgraphs. While the receding horizon and monolithic problems used very similar code to model the problem, we did have to create two sets of some of the functions. This is because the receding horizon sub-problems used the optimal values of some other sub-problems for certain constraints, whereas the monolithic problem could directly link {\it variables} of other subproblems. 

\section{Case Study 3}
Below, we provide time and memory statistics from solving the 5-bus system and the reliability test system (RTS) with graph-based NBD. The purpose of this case study was not to show if these decomposition schemes are superior to other solvers but rather to show how the graphs facilitate the decomposition schemes, so these statistics are not included in the main manuscript. However, we include them here for reference. The memory shown here was the maximum resident set size reported by the timing function, and the solution times do not include the time to build the model or the time to initialize the {\tt BendersOptimizer}, only the time to actually solve the model or to perform both the forward and backward passes. Here, Gurobi with default parameters was set as the solver for both sub-problems and the non-decomposed, full problem. The results for the 5-bus system and RTS are shown in Tables \ref{tab:5bus} and \ref{tab:RTS}.

\begin{table}[ht] 
\caption{Reported information from solving the 5-bus system using the gBD scheme in {\tt PlasmoBenders.jl} compared with solving without decomposition (``No Decomp''). Results are reported based on the number of partitions (``Num Partitions'') used for gBD with $\mathcal{G}_r = \mathcal{G}_2$ indicating that the second subgraph (rather than the first) was used as the root subgraph.}
\centering
\label{tab:5bus}
\begin{tabular}{c|c|c|c|c|c}
\hline
& \multicolumn{4}{|c|}{gBD} & \multicolumn{1}{|c}{No Decomp} \\
\cline{1-6}
Num Partitions & 3 & 3 $\mathcal{G}_r = \mathcal{G}_2$  & 5 & 10 & ---- \\
\hline
Best Solution (USD) & 4.861e5 & 4.861e5 & 4.861e5 & 4.861e5 & 4.861e5 \\
\hline
Gap & -1.2e-14\% & -2.5e-9 \% & 0.4\% & 1.7\% & 0.097\% \\
\hline
Time (min) & 5.9 & 12.4 & 5.2 & 8.0 & 64.2 \\
\hline
Required Memory (GB) & 7.8 & 10.0  & 7.5 & 4.5 & 26.8 \\
\hline
\end{tabular}
\end{table}
s

\begin{table}[ht] 
\caption{Reported information from solving the RTS using the gBD scheme in {\tt PlasmoBenders.jl} compared with solving without decomposition (``No Decomp''). Results are reported based on the number of partitions (``Num Partitions'') used for gBD.}
\centering
\label{tab:RTS}
\begin{tabular}{c|c|c|c|c}
\hline
& \multicolumn{3}{|c|}{gBD} & \multicolumn{1}{|c}{No Decomp} \\
\cline{1-5}
Num Partitions & 3 & 5 & 8 &  ---- \\
\hline
Best Solution (USD) & 6.290e6 & 6.305e6 & 6.362e6 & 6.275e6\\
\hline
Gap & 0.69\% & 1.04\% & 3.30\% & 0.29\%  \\
\hline
Time (hr) & 69.4 & 58.4 & 27.5 &  3.1 \\
\hline
Required Memory (GB) & 43.5 & 40.8 & 39.2 &  46.7 \\
\hline
\end{tabular}
\end{table}

The gBD scheme was efficient for solving the 5-bus system (Table \ref{tab:5bus}) but was slow for solving the RTS. Using 3 partitions, the gBD approach was about an order of magnitude faster for the 5-bus system compared with solving without decomposition. In contrast, for the RTS system, it was an order of magnitude slower. Both cases did provide decreases in required maximum memory however. Some problems may be better suited to decomposition approaches than others, and the 5-bus system may be efficient in part because it can use a longer time horizon within the sub problems (since the system itself is smaller). Future work could include exploring ways to speed up these algorithms. 


\bibliography{./PlasmoDecompositions}